\documentclass{amsproc}

\usepackage[T1]{fontenc}
\usepackage[utf8]{inputenc}
\usepackage[english]{babel}

\usepackage[dvipsnames]{xcolor}
\usepackage[colorlinks=true,hyperindex, linkcolor=magenta, pagebackref=false, citecolor=cyan,pdfpagelabels]{hyperref}
\hypersetup{linktocpage}

\hypersetup{
    colorlinks, linkcolor=BurntOrange, urlcolor=BurntOrange, citecolor=NavyBlue
}

\usepackage{amssymb}
\usepackage{amsmath}
\usepackage[matrix,arrow,curve,cmtip]{xy}
\entrymodifiers={+!!<0pt,\fontdimen22\textfont2>} 
\usepackage{units}
\usepackage{enumitem}
\setenumerate[1]{leftmargin=2em}
\usepackage{eucal}
\usepackage[nameinlink,capitalise,noabbrev]{cleveref} 
 
\numberwithin{equation}{section}

\theoremstyle{plain}   

\newtheorem{theorem}[equation]{Theorem}  
\newtheorem{corollary}[equation]{Corollary}     
\newtheorem{lemma}[equation]{Lemma}         
\newtheorem{proposition}[equation]{Proposition} 
\newtheorem{addendum}[equation]{Addendum} 

\theoremstyle{definition}
\newtheorem{definition}[equation]{Definition}
\newtheorem{construction}[equation]{Construction}
\newtheorem*{notation}{Notation}
\newtheorem*{terminology}{Terminology}
\newtheorem*{acknowledgements}{Acknowledgements}

\theoremstyle{remark}
\newtheorem{remark}[equation]{Remark}
\newtheorem{example}[equation]{Example}
\newtheorem{warning}[equation]{Warning}
\newtheorem{question}[equation]{Question}

\newcommand{\spaces}{\mathcal{S}} 

\DeclareMathOperator{\Spec}{Spec}

\DeclareMathOperator{\Spf}{Spf}
\DeclareMathOperator{\SubExt}{SubExt}
\DeclareMathOperator{\Hom}{Hom}
\DeclareMathOperator{\End}{End}
\DeclareMathOperator{\Aut}{Aut}

\DeclareMathOperator{\Mod}{Mod}
\DeclareMathOperator{\LMod}{LMod}
\DeclareMathOperator{\CAlg}{CAlg}

\DeclareMathOperator{\Fun}{Fun}
\DeclareMathOperator{\Aff}{Aff}
\DeclareMathOperator{\op}{op}
\DeclareMathOperator{\cp}{cp}
\DeclareMathOperator{\an}{an}
\DeclareMathOperator{\Tor}{Tor}
\DeclareMathOperator{\QCoh}{QCoh}
\DeclareMathOperator{\Cat}{Cat}
\DeclareMathOperator{\Sch}{Sch}
\DeclareMathOperator{\verylargecat}{\widehat{\operatorname{Cat}}}
\newcommand{\LPr}{\operatorname{LPr}}
\newcommand{\presheaves}{\mathcal{P}}
\DeclareMathOperator{\sheaves}{Shv}
\DeclareMathOperator{\Set}{Set}
\DeclareMathOperator{\Map}{Map}
\DeclareMathOperator{\map}{map}
\DeclareMathOperator{\id}{id}
\DeclareMathOperator{\Sp}{Sp}
\DeclareMathOperator{\Ab}{Ab}
\DeclareMathOperator{\Lat}{Lat}
\DeclareMathOperator{\Lie}{Lie}
\DeclareMathOperator{\fpqc}{fpqc}
\DeclareMathOperator{\Groth}{Groth}
\DeclareMathOperator{\lex}{lex}
\DeclareMathOperator{\Sym}{Sym}
\DeclareMathOperator{\Hyp}{Hyp}
\DeclareMathOperator{\FGroup}{FGroup}
\DeclareMathOperator{\gr}{gr}
\DeclareMathOperator{\Fil}{Fil}
\newcommand{\affinespace}{\hat{\mathbb{A}}}
\newcommand{\additivegroup}{\hat{\mathbb{G}}_a}

\newcommand{\ff}{\operatorname{f}\!\operatorname{f}}
\DeclareMathOperator{\pr}{pr}

\DeclareMathOperator{\Ind}{Ind}
\DeclareMathOperator{\Nil}{Nil}
\DeclareMathOperator{\Loc}{Loc}
\DeclareMathOperator{\Cpl}{Cpl}
\DeclareMathOperator{\cof}{cof}
\DeclareMathOperator{\Vect}{Vect}
\DeclareMathOperator{\Grp}{Grp}
\DeclareMathOperator{\geom}{geom}
\DeclareMathOperator{\THH}{THH}
\DeclareMathOperator{\aff}{aff}
\DeclareMathOperator{\loc}{loc}
\DeclareMathOperator{\Pro}{Pro}

\begin{document}

\title{Dirac geometry II: Coherent cohomology}

\author{Lars Hesselholt}
\address{Nagoya University, Japan, and University of Copenhagen, Denmark}
\email{larsh@math.nagoya-u.ac.jp}

\author{Piotr Pstr\k{a}gowski}
\address{Harvard University, USA, and Institute for Advanced Study, USA}
\email{pstragowski.piotr@gmail.com}

\thanks{Both authors were supported by the Danish National Research Foundation through the Copenhagen Center for Geometry and Topology (DNRF151). Lars Hesselholt received support from JSPS Grant-in-Aid for Scientific Research number 21K03161 and from a grant from the Institute for Advanced Study School of Mathematics, and Piotr Pstr\k{a}gowski received support from NSF Grant \#DMS-1926686 and Deutsche Forschungsgemeinschaft \#EXC-2047/1 – 390685813.}

\maketitle

\begin{abstract}
  Dirac rings are commutative algebras in the symmetric monoidal
  category of $\mathbb{Z}$-graded abelian groups with the Koszul sign
  in the symmetry isomorphism. In the prequel to this paper, we
  developed the commutative algebra of Dirac rings and defined the
  category of Dirac schemes. Here, we embed this category in the
  larger $\infty$-category of Dirac stacks, which also contains formal
  Dirac schemes, and develop the coherent cohomology of Dirac
  stacks. We apply the general theory to stable homotopy theory and
  use Quillen's theorem on complex cobordism and Milnor's theorem on
  the dual Steenrod algebra to identify the Dirac stacks corresponding
  to $\operatorname{MU}$ and $\mathbb{F}_p$ in terms of their functors
  of points. Finally, in an appendix, we develop a rudimentary theory
  of accessible presheaves.
\end{abstract}

\setcounter{tocdepth}{2}
\tableofcontents

\section{Introduction}

Whatever it is that animates anima and breathes life into higher algebra, this something leaves its trace in the structure of a Dirac ring on the homotopy groups of a commutative algebra in spectra. Dirac geometry is built from Dirac rings in the same way as algebraic geometry is built from commutative rings. 
In the prequel to this paper, we developed the commutative algebra of Dirac rings and defined the category of Dirac schemes. Here, we first embed this category in the larger $\infty$-category of Dirac stacks, which also contains formal Dirac schemes. We next develop the coherent cohomology of Dirac stacks, which amounts to a functor that to a Dirac stack $X$ assigns a presentably symmetric monoidal stable $\infty$-category $\QCoh(X)$ of quasi-coherent $\mathcal{O}_X$-modules together with a symmetric monoidal subcategory $\QCoh(X)_{\geq 0} \subset \QCoh(X)$, which is the connective part of a $t$-structure. As applications of the general theory to stable homotopy theory, we use Quillen's theorem on complex cobordism and Milnor's theorem on the dual Steenrod algebra to identify the Dirac stacks corresponding to $\operatorname{MU}$ and $\mathbb{F}_p$ in terms of their functors of points. Finally, in an appendix, we develop a rudimentary theory of accessible presheaves of anima on coaccessible $\infty$-categories.

To develop the theory of Dirac stacks and their coherent cohomology, we follow Lurie in~\cite{luriesag}. The category $\Aff$ of affine Dirac schemes is the opposite of the category $\CAlg(\Ab)$ of Dirac rings. The category of Dirac rings $\CAlg(\Ab)$ and the $\infty$-category of anima $\spaces$ are both presentable, and we define the $\infty$-category of Dirac prestacks to be the $\infty$-category of accessible presheaves of anima on $\Aff$,
$$\presheaves(\Aff) \simeq \Fun_{\mathcal{A}}(\CAlg(\Ab),\spaces).$$
We show in \cref{appendix:accessible_sheaves} that this $\infty$-category has excellent categorical properties and that it admits a number of equivalent descriptions. If $\kappa$ is a small regular cardinal, then we let $\Aff_{\kappa} \subset \Aff$ be the full subcategory spanned by the affine Dirac schemes of cardinality $< \kappa$. Since it is essentially small, the $\infty$-category $\presheaves(\Aff_{\kappa})$ of presheaves of anima on $\Aff_{\kappa}$ is presentable, and we show that
$$\textstyle{ \presheaves(\Aff) \simeq \varinjlim_{\kappa} \presheaves(\Aff_{\kappa}). }$$
It follows that the $\infty$-category $\presheaves(\Aff)$ admits all small colimits and limits, both of which are calculated pointwise. We also show that the Yoneda embedding restricts to a fully faithful functor
$$\xymatrix@C=10mm{
{ \Aff } \ar[r]^-{\phantom{,}h\phantom{,}} &
{ \presheaves(\Aff) } \cr
}$$
and that this functor exhibits the target as the cocompletion of the source.

The flat topology on $\Aff$ is the Grothendieck topology, where a sieve $j \colon U \to h(S)$ is a covering sieve if there exists a finite family of maps $(T_i \to S)_{i \in I}$ such that the induced map $T \simeq \coprod_{i \in I} T_i \to S$ is faithfully flat and such that each $h(T_i) \to h(S)$ factors through $j$.

\begin{definition}
\label{definition:dirac_stack}The $\infty$-category of Dirac stacks is the full subcategory
$$\sheaves(\Aff) \subset \presheaves(\Aff)$$
spanned by the sheaves for the flat topology.
\end{definition}

It is not immediately clear that there exists a ``sheafification'' functor, left adjoint to the canonical inclusion of sheaves in presheaves. That this is nevertheless the case is a consequence of a theorem of Waterhouse~\cite[Theorem~5.1]{waterhouse1975basically}, which is perhaps not as well-known as it deserves. It shows that the flat topology is well-behaved. What is not well-behaved is the $\infty$-category of all presheaves on $\Aff$; once we limit ourselves to considering accessible presheaves only, the supposed problems with the flat topology disappear. In the Dirac context, Waterhouse's theorem is the following statement, and the essential ingredient in its proof is the Dirac analogue of the equational criterion for faithful flatness, which we proved in \cite[Addendum~3.9]{diracgeometry1}.

\begin{theorem}[{\hyperref[theorem:waterhouse_preservation_of_flat_sheaves_under_lke]{Waterhouse's theorem}}]\label{theorem:waterhouse}Let $\kappa < \kappa'$ be small regular cardinals. In this situation, the left Kan extension along the canonical inclusion
$$\xymatrix@C=9mm{
{ \presheaves(\Aff_{\kappa}) } \ar[r]^-{i_!} &
{ \presheaves(\Aff_{\kappa'}) } \cr
}$$
preserves sheaves with respect to the flat topology.
\end{theorem}

As we have already mentioned, an important consequence of \cref{theorem:waterhouse} is the existence of a left exact ``sheafification'' functor
$$\xymatrix@C=10mm{
{ \presheaves(\Aff) } \ar@<.7ex>[r]^-{\phantom{,}L\phantom{,}} &
{ \sheaves(\Aff), } \ar@<.7ex>[l]^-{\phantom{!}\iota\phantom{!}} \cr
}$$
left adjoint to the canonical inclusion functor. Indeed, if we replace $\Aff$ by the essentially small $\Aff_{\kappa}$, then Lurie has proved in~\cite[Proposition~6.2.2.7]{luriehtt} that a left exact sheafification functor exists, and  \cref{theorem:waterhouse} and \cite[\href{https://kerodon.net/tag/02FV}{Tag 02FV}]{kerodon} show that these left adjoint functors assemble to give the stated left adjoint functor.

We can now deduce from~\cite[Theorem~6.1.0.6]{luriehtt} that the $\infty$-category of Dirac stacks satisfies the $\infty$-categorical Giraud axioms, except that it not presentable. 

\begin{theorem}[{\hyperref[theorem:giraud_axioms_for_dirac_stacks]{Giraud's axioms for Dirac stacks}}]\label{theorem:almost_topos}The $\infty$-category $\,\sheaves(\Aff)$ of Dirac stacks has the following properties:
  \begin{enumerate}[leftmargin=8mm]
\item[{\rm (i)}]It is complete and cocomplete, and it is generated under small colimits by the essential image of the Yoneda embedding $h \colon \Aff \to \sheaves(\Aff)$, which consists of $\omega_1$-compact objects and is coaccessible.
\item[{\rm (ii)}]Colimits in $\sheaves(\Aff)$ are universal.
\item[{\rm (iii)}]Coproducts in $\sheaves(\Aff)$ are disjoint.
\item[{\rm (iv)}]Every groupoid in $\sheaves(\Aff)$ is       effective.
\end{enumerate}
\end{theorem}

The $\infty$-category of Dirac stacks is close enough to being an $\infty$-topos that it retains the important property of $\infty$-topoi that the contravariant functor
$$\xymatrix{
{ \sheaves(\Aff)^{\op} } \ar[r] &
{ \verylargecat_{\infty} } \cr
}$$
that to $X$ assigns the slice $\infty$-category $\sheaves(\Aff)_{/X}$ takes colimits of Dirac stacks to limits of $\infty$-categories. This statement is the source of all descent statements in the paper. In typical situations, we wish to show that flat descent holds for the full subcategory
$$\sheaves(\Aff)_{/X}^P \subset \sheaves(\Aff)_{/X}$$
spanned by the maps $f \colon Y \to X$ that have some property $P$. But this follows from the fundamental descent statement, once we prove that the property $P$ satisfies descent for the flat topology in the sense that for every effective epimorphism $f \colon X' \to X$, the diagram of $\infty$-categories
$$\xymatrix{
{ \sheaves(\Aff)_{/X}^P } \ar[r] \ar[d] &
{ \sheaves(\Aff)_{/X'}^P } \ar[d] \cr
{ \sheaves(\Aff)_{/X} } \ar[r] &
{ \sheaves(\Aff)_{/X'} } \cr
}$$
is cartesian. For example, we define a functor that to a Dirac stack $X$ assigns the $\infty$-category $\FGroup(X)$ of formal groups over $X$ and use this strategy to show that it descends along effective epimorphisms.

It follows from Grothendieck's faithfully flat descent for graded modules, which we proved in~\cite[Theorem~1.5]{diracgeometry1}, that the flat topology on $\Aff$ is subcanonical. Here, we prove that, more generally, the category of Dirac schemes, which we defined in \cite[Definition~2.29]{diracgeometry1}, embeds fully faithfully into the $\infty$-category of Dirac stacks.

\begin{theorem}[{\hyperref[theorem:presheaf_represented_by_a_scheme_is_a_dirac_stack_and_restricted_yoneda_is_fully_faithful]{Schemes are stacks}}]Dirac schemes are Dirac stacks:
\begin{enumerate}
    \item[{\rm (1)}]If $X$ is a Dirac scheme, then the functor $h(X) \colon \Aff^{\op} \to \spaces$ that to an affine scheme $S$ assigns $\Map(S,X)$ is accessible and a sheaf for the flat topology.
    \item[{\rm (2)}]The resulting functor $h \colon \Sch \to \sheaves(\Aff)$ is fully faithful.
\end{enumerate}
\end{theorem}

This result is the starting point for the definition following To\"{e}n--Vezzosi~\cite{toenvezzosi} and Lurie~\cite{luriesagtemp} of the property of a map of Dirac stacks $f \colon Y \to X$ of being geometric. Informally, the geometric maps span the smallest full subcategory
$$\sheaves(\Aff)_{/X}^{\geom} \subset \sheaves(\Aff)_{/X}$$
that contains the maps $f \colon Y \to X$ that, locally on $X$, are equivalent to maps of Dirac schemes, and that is closed under the formation of the geometric realization of groupoids with flat face maps.\footnote{\,We recall that, in $\infty$-category theory, geometric realizations of groupoids play the same role as do quotients by equivalence relations in $1$-category theory; see~\cite[Section~6.1.2]{luriehtt}.} However, the formal definition, which we give in \cref{subsection:geometric_maps} below is more complicated, because the notions of being geometric and being flat must be defined recursively together. The property of being geometric is well-behaved: It is preserved under composition and base-change; it is local on the target; and maps in $\sheaves(\Aff)_{/X}^{\geom}$ are automatically geometric.

We now explain the theory of coherent cohomology of Dirac stacks, which we develop following Lurie~\cite[Chapter~6]{luriesag}. It consists of the following data:
\begin{enumerate}[leftmargin=10mm]
\item[(Q1)]A functor $\QCoh \colon \sheaves(\Aff)^{\op} \to \CAlg(\LPr)$; it assigns to map of Dirac stacks $f \colon Y \to X$ a symmetric monoidal adjunction
$$\xymatrix@C=10mm{
{ \QCoh(X) } \ar@<.7ex>[r]^-{f^*} &
{ \QCoh(Y) } \ar@<.7ex>[l]^-{f_*} \cr
}$$
between presentably symmetric monoidal stable $\infty$-categories.
\item[(Q2)]For every Dirac stack $X$, a symmetric monoidal adjunction
$$\xymatrix@C=10mm{
{ \QCoh(X)_{\geq 0} } \ar@<.7ex>[r]^-{i} &
{ \QCoh(X) } \ar@<.7ex>[l]^-{r} \cr
}$$
such that $i$ is the fully faithful inclusion of the connective part of a $t$-structure, and such that for every map of Dirac stacks $f \colon Y \to X$, the functor $f^*$ is right $t$-exact, or equivalently, the functor $f_*$ is left $t$-exact.
\end{enumerate}
First, to produce the functor~(Q1), we begin with the functor that to an affine Dirac scheme $S \simeq \Spec(A)$ assigns the presentably symmetric monoidal category
$$\QCoh(S)^{\heartsuit} \simeq \Mod_A(\Ab)$$
of graded $A$-modules. It promotes via animation to a functor that to $S$ assigns the presentably symmetric monoidal $\infty$-category
$$\QCoh(S)_{\geq 0} \simeq \Mod_A(\Ab)^{\an}$$
of animated graded $A$-modules, which, in turn, promotes via stabilization to a functor that to $S$ assigns the presentably symmetric monoidal $\infty$-category
$$\QCoh(S) \simeq \Sp(\Mod_A(\Ab)^{\an})$$
of spectra in animated graded $A$-modules. Following Grothendieck, we show that this functor is a sheaf for the flat topology on $\Aff$, a fact, which, in homotopy theory, was first observed and exploited by Hopkins~\cite{hopkins1999complex}. This gives us~(Q1):

\begin{theorem}[{\hyperref[thm:qcohstackification]{Faithfully flat descent for quasi-coherent modules}}]
\label{theorem:flat_descent_for_qcoh}The right Kan extension of the functor $\QCoh \colon \Aff^{\op} \to \CAlg(\LPr)$ along $h \colon \Aff \to \presheaves(\Aff)$ admits a unique factorization
$$\xymatrix@C=12mm{
{ \sheaves(\Aff)^{\op} } \ar[r]^-{\QCoh} &
{ \CAlg(\LPr) } \cr
}$$
through the sheafification functor.
\end{theorem}

To spell out this definition of $\QCoh(X)$, if we choose a regular cardinal $\kappa$ such that $X \in \sheaves(\Aff)$ is the left Kan extension of $X_{\kappa} \in \sheaves(\Aff_{\kappa})$, then
$$\textstyle{
\QCoh(X) \simeq \varprojlim_{\;\eta \colon \!S \to X}\QCoh(S), }$$
where the limit is indexed by $((\Aff_{\kappa})_{/X})^{\op}$ and is calculated in $\CAlg(\LPr)$.\footnote{The forgetful functors
$\CAlg(\LPr) \to \LPr \to \verylargecat_{\infty}$ preserve and reflect limits} 

Second, to produce the symmetric monoidal adjunction~(Q2), we show that for an affine Dirac scheme $S$, the $\infty$-category $\QCoh(S)_{\geq 0}$ is Grothendieck prestable in the sense of~\cite[Definition~C.1.4.2]{luriesag}. Thus, the symmetric monoidal adjunction
$$\xymatrix@C=10mm{
{ \QCoh(S)_{\geq 0} } \ar@<.7ex>[r]^-{\Sigma^{\infty}} &
{ \QCoh(S) } \ar@<.7ex>[l]^-{\Omega^{\infty}} \cr
}$$
is the connective part of a $t$-structure, and we obtain~(Q2) by right Kan extension. More concretely, if $X$ is a Dirac stack, then $\mathcal{F} \in \QCoh(X)$ is connective if and only if $\eta^*(\mathcal{F}) \in \QCoh(S)$ is connective for every map $\eta \colon S \to X$ with $S$ affine. It is also true that if $\eta^*(\mathcal{F}) \in \QCoh(S)$ is coconnective for every map $\eta \colon S \to X$ with $S$ affine, then $\mathcal{F} \in \QCoh(X)$ is coconnective, but the converse is false. In particular, the tensor unit $\mathcal{O}_X$ belongs to the heart of the $t$-structure.

In general, the $t$-structure on $\QCoh(X)$ is not left complete or right complete, and the coconnective part is next to impossible to understand. However, for geometric Dirac stacks, including Dirac schemes, the situation improves:

\begin{theorem}[{\hyperref[theorem:qcoh_on_geometric_dirac_stack_is_grothendieck_prestable]{Geometric stacks and $t$-structure}}]
\label{proposition:t-structure_for_geometric_stacks}
If $X$ is a geometric Dirac stack, then the $t$-structure on $\QCoh(X)$ is left and right complete and the coconnective part is closed under filtered colimits. Moreover, if $\eta \colon S \to X$ is submersive with $S$ affine, then $\mathcal{F} \in \QCoh(X)$ is coconnective if and only if $\eta^*(\mathcal{F}) \in \QCoh(S)$ is so.
\end{theorem}

We let $\pi_i^{\heartsuit}(\mathcal{F}) \in \QCoh(X)^{\heartsuit}$
be the $i$th homotopy $\mathcal{O}_X$-module of $\mathcal{F} \in \QCoh(X)$ with respect to the canonical $t$-structure. We refer to $i$ as the ``animated'' degree to separate it from the ``spin'' degree intrinsic to Dirac geometry. 
If $f \colon Y \to X$ is a map of Dirac stacks, then the ``correct'' relative coherent cohomology $\mathcal{O}_X$-modules of the $\mathcal{O}_Y$-module $\mathcal{G} \in \QCoh(Y)$ are given by
$$R^{\hspace{.6pt}i\hspace{-.6pt}}f_*(\mathcal{G}) := \pi_{-i}^{\heartsuit}(f_*(\mathcal{G})).$$
This functor can generally not be calculated as the derived functor of the left-exact functor between abelian categories $f_*^{\heartsuit}$ induced by the left $t$-exact functor $f_*$.

Another special class of Dirac stacks that we consider are the formal Dirac schemes. If $\eta \colon S \to X$ is a closed immersion of Dirac schemes, whose defining quasi-coherent ideal $\mathcal{I} \subset \mathcal{O}_X$ is of finite type, then we define the formal completion of $X$ along $S$ to be the colimit in Dirac stacks
$$\xymatrix{
{ Y \simeq \varinjlim_{\,m} S^{(m)} } \ar[r]^-{j} &
{ X } \cr
}$$
of the infinitesimal thickenings $\eta^{(m)} \colon S^{(m)} \to X$. The map $j$ is not geometric, but if the closed immersion $\eta \colon S \to X$ is regular in the sense that, locally on $X$, it is defined by a regular sequence, then it behaves as an open and affine immersion of a tubular neighborhood of $S$ in $X$, as envisioned by Grothendieck.

\begin{theorem}[{\hyperref[thm:recollement]{Formal completion and recollement}}]
\label{theorem:formal_completion_qcoh}Let $X$ be a Dirac scheme, and let $\eta \colon S \to X$ be a regular closed immersion. Let $j \colon Y \to X$ be the formal completion of $X$ along $S$, and let $i \colon U \simeq X \smallsetminus S \to X$ be the inclusion of the open complement of $S$. In this situation, there is a stable recollement
\vspace{-1mm}
$$\begin{xy}
(0,0)*+{ \QCoh(U) };
(30.5,0)*+{ \QCoh(X) };
(61,0)*+{ \QCoh(Y) };
(7,0)*+{}="1";
(23,0)*+{}="2";
(38,0)*+{}="3";
(54,0)*+{}="4";
(0,-5)*+{};
{ \ar@<.5ex>@/^.9pc/_-{i^!} "1";"2";};
{ \ar^-{i_* \simeq\, i_!} "2";"1";};
{ \ar@<-1ex>@/_.9pc/_-{i^*} "1";"2";};  
{ \ar@<.5ex>@/^.9pc/_-{j_*} "3";"4";};
{ \ar^-{j^! \simeq\, j^*} "4";"3";};
{ \ar@<-1ex>@/_.9pc/_-{j_!} "3";"4";};  
\end{xy}$$
and, in addition, the functor $j_*$ is $t$-exact.
\end{theorem}

The formal schemes that we will consider are the formal affine spaces, defined as follows. If $S$ is a Dirac scheme and $\mathcal{E} \in \QCoh(S)^{\heartsuit}$ an $\mathcal{O}_S$-module locally free of finite rank, then the affine space associated with $\mathcal{E}$ is the affine map
$$\xymatrix{
{ \mathbb{A}_S(\mathcal{E}) \simeq \Spec(\Sym_{\mathcal{O}_S}(\mathcal{E})) } \ar[r]^-{p} &
{ S, } \cr
}$$
and the formal affine space associated with $\mathcal{E}$ is its formal completion
$$\xymatrix{
{ \affinespace_S(\mathcal{E}) \simeq \Spf(\Sym_{\mathcal{O}_S}(\mathcal{E})) } \ar[r]^-{q} &
{ S } \cr
}$$
along the zero section. We define a map of Dirac stacks $q \colon Y \to X$ to be a formal hyperplane if its base change along any map $\eta \colon S \to X$ from an affine Dirac scheme is equivalent to a formal affine space. This property is stable under base change, but it does not descend along effective epimorphisms in general, even if $X$ is a Dirac scheme.\footnote{\,For a counterexample, see~\cite[Warning~1.1.22]{lurieell2}.} However, if we fix a section, then it does. We consider a formal affine space to be pointed by the zero section.

\begin{definition}
\label{definition:formal_hyperplane}
A pointed Dirac stack
$$\xymatrix@C=6mm{
{} &
{ Y } \ar[dr]^-{q} &
{} \cr
{ X } \ar[rr]^-{\id_X} \ar[ur]^-{\eta} &&
{ X } \cr
}$$
is a pointed formal hyperplane if its base change along any map $\eta \colon S \to X$ from an affine Dirac scheme is equivalent to a pointed formal affine space.
\end{definition}

By definition, the $\infty$-category of pointed formal hyperplanes over a Dirac stack is the full subcategory of the $\infty$-category of pointed Dirac stacks
$$\Hyp_*(X) \subset (\sheaves(\Aff)_{/X})_*$$
spanned by the pointed formal hyperplanes. We show that it is equivalent to a $1$-category and that the property of being a pointed formal hyperplane descends along effective epimorphisms of Dirac stacks $f \colon X' \to X$ in the sense that
$$\xymatrix@C=10mm{
{ \Hyp_*(X) } \ar[r]^-{f^*} \ar[d] &
{ \Hyp_*(X') } \ar[d] \cr
{ (\sheaves(\Aff)_{/X})_* } \ar[r]^-{f^*} &
{ (\sheaves(\Aff)_{/X'})_* } \cr
}$$
is a cartesian diagram of $\infty$-categories. So by \cref{theorem:almost_topos}, the $\infty$-category of pointed formal hyperplanes satisfies descent along effective epimorphisms:

\begin{theorem}[{\hyperref[cor:formal_hyperplanes_descent_along_effective_epis]{Flat descent for pointed formal hyperplanes}}]
\label{theorem:descent_for_pointed_formal_hyperplanes}If $f \colon X' \to X$ is an effective epimorphism of Dirac stacks, then the canonical maps
$$\xymatrix{
{ \Hyp_*(X) } \ar[r] &
{ \varprojlim_{\,[n] \in \Delta} \Hyp_*(X'^{\times_X[n]}) } \ar[r] &
{ \varprojlim_{\,[n] \in \Delta_{\leq 2}} \Hyp_*(X'^{\times_X[n]}) } \cr
}$$
are equivalences of $\infty$-categories.
\end{theorem}

We recall from~\cite[Proposition~A.1]{diracgeometry1} that the right-hand map in \cref{theorem:descent_for_pointed_formal_hyperplanes} is an equivalence, because the $\infty$-category of pointed formal hyperplanes is equivalent to a $1$-category, and moreover, the right-hand term is equivalent to the $1$-category of pointed formal hyperplanes over $X'$ with descent data along $f \colon X' \to X$. 

Let $\mathcal{C}$ be an $\infty$-category which admits finite products. Following Lawvere, we let $\Lat$ be the category of finitely generated free abelian groups and define the $\infty$-category $\Ab(\mathcal{C})$ of abelian group objects in the $\infty$-category $\mathcal{C}$ to be the full subcategory
$$\Ab(\mathcal{C}) \simeq \Fun^{\Sigma}(\Lat^{\op},\mathcal{C}) \subset \Fun(\Lat^{\op},\mathcal{C})$$
spanned by the functors that preserve finite products.

\begin{definition}
\label{definition:formal_groups_introduction}The $\infty$-category of formal groups over a Dirac stack $X$ is the $\infty$-category
$$\FGroup(X) \simeq \Ab(\Hyp(X))$$
of abelian group objects in the $\infty$-category of formal hyperplanes over $X$.
\end{definition}

The $\infty$-category $\FGroup(X)$ is equivalent to a $1$-category. If $\mathcal{G}$ is a formal group over $X$, then the underlying hyperplane $\mathcal{G}(\mathbb{Z})$ is pointed by the zero section $\mathcal{G}(0) \to \mathcal{G}(\mathbb{Z})$. Therefore, by \cref{theorem:descent_for_pointed_formal_hyperplanes}, if $f \colon X' \to X$ is an effective epimorphism of Dirac stacks, then the canonical map
$$\xymatrix{
{ \FGroup(X) } \ar[r] &
{ \varprojlim_{\,[n] \in \Delta} \FGroup(X'^{\times_X[n]}) } \cr
}$$
is an equivalence of $\infty$-categories, and the target $\infty$-category is equivalent to the $1$-category of formal groups over $X'$ with descent data along $f \colon X' \to X$. Moreover, if $\mathcal{G}$ is a formal group over $X$, then we define its Lie algebra to be the vector bundle
$$\Lie(\mathcal{G}) \simeq V_X(\mathcal{N}_{X/Y}) \in \Vect(X)$$
associated with the conormal sheaf at the zero section. Since $\mathcal{G}$ is abelian, the Lie bracket on this vector bundle will be zero, so we do not take the trouble to define it.

We end with a number of applications of Dirac geometry to algebraic topology, which we encode as a six-functor formalism on the $\infty$-category of anima following the general theory of Liu--Zheng~\cite{liuzheng} and Mann~\cite[Appendix~A.5]{mann}. The definition of this six-functor formalism is informed by the fact that anima are ``discrete'' objects, as opposed to ``continuous'' objects. So we declare that every map of anima $f \colon T \to S$ is a local isomorphism and that a map of anima $f \colon T \to S$ is proper if its fibers are equivalent to finite sets, including the empty set. By~\cite[Proposition~A.5.10]{mann}, this determines a six-functor formalism on $\mathcal{S}$ that to an anima $S$ assigns the functor $\infty$-category
$$\Sp^S \simeq \Fun(S,\Sp)$$
and that to a map of anima $f \colon T \to S$ assigns the adjoint functors
\vspace{-1mm}
$$\begin{xy}
(0,0)*+{ \Sp^T };
(23,0)*+{ \Sp^S, };
(3,0)*+{}="1";
(20,0)*+{}="2";
(0,-5)*+{};
{ \ar@<.9ex>@/^.9pc/^-{f_!} "2";"1";};
{ \ar_-{f^! \simeq\, f^*} "1";"2";};
{ \ar@<-.5ex>@/_.9pc/^-{f_*} "2";"1";};
\end{xy}$$
where $f^! \simeq f^*$ is the restriction along $f$, and  $f_!$ and $f_*$ are the left and right Kan extensions along $f$.

Given a commutative algebra in spectra $E$ and a group in anima $G$, we apply the general theory to construct a formal Dirac scheme
$$\xymatrix{
{ Y \simeq Y_{E,G} } \ar[r]^-{q} &
{ S \simeq S_E }
}$$
over $S \simeq \Spec(R)$ with $R \simeq \pi_*(E)$, and we will consider situations in which $q$ is a formal hyperplane. By varying $G$, we will obtain formal groups, and by varying $E$, we will obtain formal groups over geometric Dirac stacks, which we characterize by exhibiting their functors of points.

To this end, we recall from \cite[Theorem~5.6.2.10]{lurieha} that every group in anima $G$ is the loop group $\Omega BG$ of a pointed $1$-connective anima
$$\xymatrix{
{ 1 } \ar[r]^-{s} &
{ BG, } \cr
}$$
and we apply the general theory to the unique map
$$\xymatrix{
{ BG } \ar[r]^-{p} &
{ 1. } \cr
}$$
The $\infty$-category $\Sp^{BG}$ is the $\infty$-category of spectra with $G$-action, and $p_!$ and $p_*$ are the functors that to a spectrum with $G$-action $Y$ assign its homotopy orbit spectrum $p_!(Y) \simeq Y_{hG}$ and its homotopy fixed point spectrum $p_*(Y) \simeq Y^{hG}$, respectively, whereas $p^*$ is the functor that to a spectrum $X$ assigns the spectrum with trivial $G$-action $Y \simeq p^*(X)$. We only consider spectra with trivial $G$-action.

Given a commutative algebra in spectra $E$, the spectrum 
$$p_*p^*(E) \simeq \map(\Sigma_+^{\infty}BG,E)$$
promotes to a commutative algebra in spectra, and the unit maps
$$\xymatrix{
{ E } \ar[r] &
{ p_*p^*(E) } \ar[r] &
{ p_*s_*s^*p^*(E) \simeq E } \cr
}$$
promote to maps of commutative algebras in spectra. Hence, if $R \simeq \pi_*(E)$ and $A \simeq \pi_*(p_*p^*(E))$ are the respective Dirac rings of homotopy groups, then these maps give $A$ the structure of an augmented $R$-algebra. We now define
$$\xymatrix{
{ Y \simeq Y_{E,G} \simeq \Spf(A) } \ar[r]^-{q} &
{ S \simeq S_E \simeq \Spec(R) } \cr
}$$
to be the formal completion of $X \simeq \Spec(A)$ along the closed immersion
$$\xymatrix{
{ S } \ar[r]^-{\eta} &
{ X } \cr
}$$
defined by the augmentation.

First, we let $G$ be the Pontryagin dual $\widehat{L} \simeq \Hom(L,U(1))$ of $L \in \Lat$. If $E$ is complex orientable, then the formal scheme
$$\xymatrix{
{ Y \simeq Y_{E,\widehat{L}} } \ar[r]^-{q} &
{ S \simeq S_E } \cr
}$$
is a formal hyperplane, and as $L$ varies, this defines the Quillen formal group
$$\xymatrix@C=10mm{
{ \Lat^{\op} } \ar[r]^-{\mathcal{G}_E^Q} &
{ \Hyp(S) } \cr
}$$
associated with $E$. Moreover, by considering the Postnikov filtration of $E$, we obtain a canonical isomorphism of line bundles
$$\xymatrix@C=10mm{
{ \mathbb{A}_S(1) } \ar[r]^-{\phi_E^Q} &
{ \Lie(\mathcal{G}_E^Q) } \cr
}$$
from the spin-$1$ affine line over $S$ to the Lie algebra of the Quillen formal group. We now let $E$ vary through the tensor powers of the commutative algebra in spectra $\operatorname{MU}$ representing complex cobordism. In this way, we obtain a simplicial formal group
$$\xymatrix@C=10mm{
{ \Delta^{\op} } \ar[r] &
{ \FGroup } \cr
}$$
that to $[n]$ assigns the Quillen formal group associated with $\operatorname{MU}^{\otimes[n]}$. Its geometric realization is a formal group $\mathcal{G}^Q$ over the geometric Dirac stack
$$\textstyle{ X^Q \simeq \varinjlim_{\,[n] \in \Delta^{\op}} S_{\,\operatorname{MU}^{\otimes[n]}} }$$
equipped with a trivialization $\phi^Q$ of its Lie algebra. We show that Quillen's theorem on complex cobordism, \cite[Theorem~2]{quillen6}, implies and is implied by the statement that $(\mathcal{G}^Q,\phi^Q)$ is the universal $1$-dimensional spin-$1$ formal group equipped with a trivialization of its Lie algebra.

\begin{theorem}[{\hyperref[thm:quillen]{Quillen's theorem}}]
\label{theorem:quillen's_theorem}Let $T$ be a Dirac stack. The functor that to $f \colon T \to X^Q$ assigns the pair $(f^*\mathcal{G}^Q,f^*\phi^Q)$ is an equivalence of anima from $X^Q(T)$ to the groupoid of pairs $(\mathcal{G},\phi)$ of a formal group $\mathcal{G}$ over $T$ and an isomorphism
$$\xymatrix{
{ \mathbb{A}_T(1) } \ar[r]^-{\phi} &
{ \Lie(\mathcal{G}) } \cr
}$$
of line bundles over $T$ from the spin-$1$ affine line.
\end{theorem}

We also show that if $T$ is an affine Dirac scheme and $(\mathcal{G},\phi) \in X^Q(T)$, then the trivialization $\phi$ lifts to a spin-$1$ global coordinate on the formal group $\mathcal{G}$.

Second, we let $p$ be an odd prime number, and let $G$ be the $p$-torsion subgroup
$$\xymatrix{
{ \widehat{L}[p] \simeq \Hom(L,C_p) } \ar[r]^-{i} &
{ \widehat{L} \simeq \Hom(L,U(1)) } \cr
}$$
of the Pontryagin dual of $L \in \Lat$. If $k \simeq \mathbb{F}_p$, then for every commutative $k$-algebra in spectra $E$,
$$\xymatrix{
{ Y \simeq Y_{E,\widehat{L}[p]} } \ar[r]^-{q} &
{ S \simeq S_E } \cr
}$$
is a formal hyperplane, and as $L$ varies, this defines the Milnor formal group
$$\xymatrix@C=10mm{
{ \Lat^{\op} } \ar[r]^-{\mathcal{G}_E^M} &
{ \Hyp(S) } \cr
}$$
associated with $E$. The canonical inclusion $i$ induces a map of formal groups
$$\xymatrix{
{ \mathcal{G}_E^M } \ar[r] &
{ \mathcal{G}_E^Q } \cr
}$$
to the Quillen formal group. It is surjective and its kernel
$$\xymatrix{
{ \Fil^1\mathcal{G}_E^M } \ar[r] &
{ \Fil^0\mathcal{G}_E^M \simeq \mathcal{G}_E^M } \cr
}$$
is again a formal group over $S$. This filtration gives rise to a grading of the Lie algebra $\Lie(\mathcal{G}_E^M)$, which we refer to as the ``charge'' grading.\footnote{\,This grading encodes the residual $U(1) \simeq U(1)/C_p$-action on $BC_p$.} Moreover, from the Postnikov filtration of $E$, we obtain a canonical isomorphism of graded vector bundles over $S$,
$$\xymatrix@C=10mm{
{ \mathbb{A}_S(e,\gamma) } \ar[r]^-{\phi_E^M} &
{ \gr\Lie(\mathcal{G}_E^M), } \cr
}$$
where $e$ has spin $\nicefrac{1}{2}$ and charge $1$, and where $\gamma$ has spin $1$ and charge $0$.

Letting $E$ vary through the tensor powers of $k$, we obtain a simplicial filtered formal group that to $[n] \in \Delta^{\op}$ assigns the Milnor formal group associated with $k^{\otimes [n]}$ with the filtration defined by the map to the Quillen formal group. Its geometric realization is a filtered formal group $\Fil\mathcal{G}^M$ over the geometric Dirac stack
$$\textstyle{ X^M \simeq \varinjlim_{\,[n] \in \Delta^{\op}} S_{k^{\otimes[n]}} }$$
equipped with a trivialization $\phi^M$ of its graded Lie algebra. We show that Milnor's theorem~\cite{milnor} on the structure of the dual Steenrod algebra implies and is implied by the following statement.

\begin{theorem}[{\hyperref[thm:milnor]{Milnor's theorem}}]
\label{theorem:milnor's_theorem}Let $S \simeq \Spec(\mathbb{F}_p)$ with $p$ odd, and let $T$ be an $S$-Dirac stack. The functor that to a map of $S$-Dirac stacks $f \colon T \to X^M$ assigns the pair $(f^*\Fil\mathcal{G}^M,f^*\phi^M)$ is an equivalence from $X^M(T)$ to the groupoid of pairs $(\Fil\mathcal{G},\phi)$ of a filtered formal group $\Fil\mathcal{G}$ over $T$, locally isomorphic to a filtered additive formal group, and an isomorphism of graded vector bundles over $T$,
$$\xymatrix{
{ \mathbb{A}_T(e,\gamma) } \ar[r]^-{\phi} &
{ \gr\Lie(\mathcal{G}), } \cr
}$$
where $e$ has spin $\nicefrac{1}{2}$ and charge $1$, and where $\gamma$ has spin $1$ and charge $0$.
\end{theorem}

Since $k \simeq \mathbb{F}_p$ is a prime field, the face maps $d_0,d_1 \colon S_{k^{\otimes[1]}} \to S_{k^{\otimes[0]}}$ are equal, so the groupoid Dirac scheme $S_{k^{\otimes[-]}}$ is in fact a group Dirac scheme $G^M$ and
$$X^M \simeq BG^M$$
is its classifying Dirac stack. The dual Steenrod algebra is the Hopf algebra in Dirac $k$-vector spaces corresponding to $G^M$.

Finally, we express the method of descent, which goes back to Adams~\cite{adams2}, in Dirac geometric terms. Given a map of commutative algebras in spectra $\phi \colon k \to E$, we may form the Dirac stack
$$\textstyle{ X \simeq \varinjlim_{[n] \in \Delta^{\op}} \Spec(\pi_*(E^{\otimes_k[n]})). }$$
This Dirac stack and its coherent cohomology is particularly useful if the face maps in the simplicial Dirac scheme are flat, in which case, we construct a functor
$$\xymatrix{
{ \Mod_k(\Sp) } \ar[r]^-{\mathcal{F}} &
{ \QCoh(X)^{\heartsuit} } \cr
}$$
that to a $k$-module in spectra $V$ assigns a quasi-coherent $\mathcal{O}_X$-module $\mathcal{F}(V)$. This functor gives rise to an Adams filtration of $V$, and the associated spectral sequence is the Adams or descent spectral sequence
$$E_{i,j}^2 = H^{-i}(X,\mathcal{F}(V)(\nicefrac{j}{2})) \Rightarrow \pi_{i+j}(\widehat{V}).$$
It starts from the coherent cohomology of $X$ with coefficients in the half-integer Serre twists of $\mathcal{F}(V)$, and it converges conditionally to the homotopy groups of the completion of $V$ with respect to the Adams filtration. A particularly understandable application of this method was given by Liu--Wang~\cite{liuwang} in the situation of the base-change map in topological Hochschild homology
$$\xymatrix{
{ k \simeq \THH(\mathcal{O}_K/\,\mathbb{S}_W) } \ar[r]^-{\phi} &
{ E \simeq \THH(\mathcal{O}_K/\,\mathbb{S}_W[z]), } \cr
}$$
where the descent spectral sequence is concentrated on the lines $-1 \leq i \leq 0$, and therefore, degenerates for degree-reasons.

\begin{terminology}We write $\Ab$  for the symmetric monoidal category of $\mathbb{Z}$-graded abelian groups with the Koszul sign in the symmetry isomorphism, and we define the category of Dirac rings to be the category $\CAlg(\Ab)$ of commutative algebras therein.
We use the terminology of Clausen--Jansen~\cite[Theorem~2.19]{clausenjansen} and say that a map $f \colon L \to K$ between small $\infty$-categories is a $\smash{ \varinjlim }$-equivalence if for every diagram $X \colon K \to \mathcal{C}$, the induced map
$$\xymatrix{
{ \varinjlim_Lf^*(X) } \ar[r] &
{ \varinjlim_KX } \cr
}$$
is an equivalence. These maps are referred to as cofinal in \cite[Section~4.1]{luriehtt}.
\end{terminology}

\begin{acknowledgements}It is a pleasure to thank Maxime Ramzi for
  carefully reading a first draft of
  \cref{appendix:accessible_sheaves} and for several helpful
  suggestions. The second author also would like to thank Kazuhiro
  Fujiwara and Nagoya University and Takeshi Saito and the University
  of Tokyo for their support and hospitality during his visit to
  Japan, where part of this paper was written. Finally, we thank an
  anonymous referee for helpful comments.
\end{acknowledgements}

\section{Dirac stacks}

We wish to define a Dirac stack to be a presheaf $X \colon \Aff^{\op} \to \spaces$ of anima on the category of affine Dirac schemes, which is well-behaved in the sense that
\begin{enumerate}
\item[(1)] $X$ to be determined by a small amount of data, and
\item[(2)] $X$ to be determined by local data,
\end{enumerate}
and both of these conditions require some care.

In the case of $(1)$, the category $\Aff$ of affine Dirac schemes is not small, so not every presheaf of anima on $\Aff$ can be written as a small colimit of representable presheaves. Our solution, which is similar to the one used by Clausen--Scholze~\cite{clausenscholzecondensed}, is to only consider the presheaves which can be expressed as small colimits of representable presheaves. We show in \cref{appendix:accessible_sheaves} that this condition is equivalent to the functor $X$ being accessible and that it leads to an excellent theory.

In the case of (2), this leads us to consider presheaves which satisfy the sheaf condition with respect to some Grothendieck topology. We would like to use the fpqc-topology, which is very strong but still subcanonical. In general, the category of fpqc-covering sieves of a fixed affine Dirac scheme does not have a small cofinal subcategory, and hence, it is not clear that a sheafication functor exists. We extend the work of Waterhouse \cite{waterhouse1975basically} to show our solution to~(1) also solves~(2) in that an accessible presheaf admits a sheafification which is again accessible.

\begin{remark}
We will show that there exists a fully faithful embedding
$$\xymatrix@C=10mm{
{ \Sch } \ar[r]^-{h'} &
{ \sheaves(\Aff) } \cr
}$$
of the category of Dirac schemes, which we defined in~\cite[Definition~2.29]{diracgeometry1}, into the $\infty$-category of Dirac stacks defined according to the above principles.
\end{remark}

\subsection{The flat topology}

The category $\Aff$ of affine Dirac schemes is not small, but it is coaccessible in the sense that its opposite category $\Aff^{\op} \simeq \CAlg(\Ab)$ is accessible. We show in \cref{prop:equivalent_conditions_for_accessability} that a presheaf
$$\xymatrix@C=8mm{
{ \Aff^{\op} } \ar[r]^-{\phantom{,}X\phantom{,}} &
{ \spaces } \cr
}$$
is an accessible functor if and only if it is a small colimit of representable presheaves. 

\begin{definition}
\label{def:dirac_prestacks}
The $\infty$-category of Dirac prestacks is the full subcategory
$$\presheaves(\Aff) \subset \Fun(\Aff^{\op},\spaces)$$
spanned by the accessible presheaves.
\end{definition}

We show in \cref{theorem:giraud_axioms_for_accessible_presheaves} that $\presheaves(\Aff)$ is complete and cocomplete and that it satisfies the Giraud axioms, with the exception that it is not presentable. We also show that the Yoneda embedding factors through a fully faithful functor
$$\xymatrix{
{ \Aff } \ar[r]^-{\phantom{,}h\phantom{,}} &
{ \mathcal{P}(\Aff) } \cr
}$$
and that this functor exhibits $\mathcal{P}(\Aff)$ as the cocompletion of $\Aff$.

We proceed to define the flat topology on $\Aff$ and to show that it gives rise to a well-behaved subcategory $\sheaves(\Aff) \subset \presheaves(\Aff)$ of accessible flat sheaves. We recall the general notions of sieves and Grothendieck topologies in \cref{appendix:effective_epimorphisms_of_sheaves}.

\begin{definition}
\label{def:fpqccovering}
A sieve $j \colon U \to h(S)$ on an affine Dirac scheme $S$ is a covering sieve for the flat topology if there exists a finite family $(f_i \colon T_i \to S)_{i \in I}$ of maps of affine Dirac schemes such that the induced map
$$\xymatrix{
{ T \simeq \coprod_{i \in I} T_i } \ar[r]^-{f} &
{ S } \cr
}$$
is faithfully flat and such that each $h(f_i) \colon h(T_i) \to h(S)$ factors through $j$.
\end{definition}

\begin{proposition}
The collection of covering sieves for the flat topology forms a finitary Grothendieck topology on $\Aff$.
\end{proposition}

\begin{proof}
The assumptions of~\cite[Proposition~A.3.2.1]{luriesag} are satisfied, since coproducts in $\Aff$ are universal and since flat maps in $\Aff$ are stable under composition and under base-change along arbitrary maps.
\end{proof}

\begin{proposition}
\label{proposition:explicit_conditions_to_be_an_fpqc_sheaf}
A presheaf $X \colon \Aff^{\op} \to \spaces$ is a sheaf for the flat topology if and only if it has the following properties:
\begin{enumerate}
\item[{\rm (i)}]For every finite family $(S_i)_{i \in I}$ of objects in $\Aff$, the induced map
$$\xymatrix{
{ X(\coprod_{i \in I}S_i) } \ar[r] &
{ \prod_{i \in I} X(S_i) } \cr
}$$
is an equivalence of anima.
\item[{\rm (ii)}]For every faithfully flat map $f \colon T \to S$ in $\Aff$, the induced map
$$\xymatrix{
{ X(S) } \ar[r] &
{ \varprojlim_{\Delta} X(T^{\times_S[-]}) } \cr
}$$
is an equivalence of anima.
\end{enumerate}
\end{proposition}

\begin{proof}
This is~\cite[Proposition~A.3.3.1]{luriesag}.
\end{proof}

\begin{definition}
\label{def:diracstack}
The $\infty$-category of Dirac stacks is the full subcategory
$$\sheaves(\Aff) \subset \presheaves(\Aff)$$
spanned by the sheaves for the flat topology. 
\end{definition}

Let us emphasize that for a presheaf $X \colon \Aff^{\op} \to \spaces$ to be a Dirac stack, it must be both accessible and a sheaf for the flat topology. The former is the case if and only if the equivalent conditions of \cref{prop:equivalent_conditions_for_accessability} are satisfied, and the latter is the case if and only if the conditions of \cref{proposition:explicit_conditions_to_be_an_fpqc_sheaf} are satisfied.

\begin{remark}
\label{rem:diracstack}
We will only consider sheaves on $\Aff$ for the flat topology. However, it is also possible to consider sheaves for the \'{e}tale topology, whose covering sieves are defined by substituting \'{e}tale for flat in \cref{def:fpqccovering}. Since \'{e}tale maps in $\Aff$ are flat by \cite[Theorem 1.9]{diracgeometry1}, the flat topology is finer than the \'{e}tale topology.
\end{remark}

A priori, it is not clear that the inclusion $\iota \colon \presheaves(\Aff) \to \sheaves(\Aff)$ admits a left adjoint ``sheafification'' functor. However, we now show that, by an argument of Waterhouse~\cite[Theorem~5.1]{waterhouse1975basically}, such a functor does indeed exist. It will be convenient to work with Dirac rings instead of affine Dirac schemes, and we will use the following notion throughout the proof.

\begin{definition}\label{definition:sub-extension}
A subextension of a map of Dirac rings $f \colon A \to B$ is a diagram of rings of the form
$$\xymatrix{
{ C } \ar[r]^-{g} \ar[d]^-{i} &
{ D } \ar[d]^-{j} \cr
{ A } \ar[r]^-{f} &
{ B } \cr
}$$
with $i$ and $j$ inclusions of sub-Dirac rings. The subextensions of $f \colon A \to B$ form a partially ordered set $\SubExt(f)$.
\end{definition}

\begin{remark}\label{remark:sub-extension}
The partially ordered set $\SubExt(f)$ of $f \colon A \to B$ is equivalent to the full subcategory of $\Fun(\Delta^1,\CAlg(\Ab))_{/f}$ spanned by the monomorphisms.
\end{remark}

\begin{lemma}
\label{lem:approximation_of_kappa_small_subextensions_by_ff_ones}
Let $\kappa$ be an uncountable cardinal, let $f \colon A \to B$ be a faithfully flat map of Dirac rings, and let $g' \colon C' \to D'$ be a subextension of $f \colon A \to B$ such that $C'$ and $D'$ are $\kappa$-small. There exists a subextension $g \colon C \to D$ of $f \colon A \to B$ such that $g' \colon C' \to D'$  is a subextension of $g \colon C \to D$, such that $C$ and $D$ are $\kappa$-small, and such that $g \colon C \to D$ is faithfully flat.
\end{lemma}

\begin{proof}
In the case of ordinary rings, this is proved in~\cite[Lemma~3.1]{waterhouse1975basically}. We recall the argument for the convenience of the reader and to verify that that it works without much change in the Dirac setting.

We first recall the necessary commutative algebra. By the equational criterion for flatness \cite[Theorem 3.8]{diracgeometry1}, a map of Dirac rings $f \colon A \to B$ is flat if and only if for every system of linear equations
\begin{equation}
\label{equation:equation_for_flatness}
    \textstyle{ \sum_{1 \leq k \leq n} y_kc_{ki} = 0 \hspace{12mm}
      (1 \leq i \leq m), }
\end{equation}
where the $c_{ki}$ are homogeneous elements of $A$, every solution $(y_k)_{1 \leq k \leq n}$, consisting of homogeneous elements of $B$, is a linear combination
$$\textstyle{ y_k = \sum_{1 \leq j \leq p} b_jz_{jk} \hspace{12mm} (1 \leq k \leq n), }$$
where the $b_j$ are homogeneous elements of $B$, of solutions $(z_{jk})_{1 \leq k \leq n}$ to (\ref{equation:equation_for_flatness}) consisting of homogeneous elements $z_{jk}$ of $A$.

By the equational criterion for faithful flatness, which we proved in \cite[Addendum~3.9]{diracgeometry1}, a flat map of Dirac rings $f \colon A \to B$ is faithfully flat if and only if it is a monomorphism and for every solution $(y_k)_{1 \leq k \leq n}$ consisting of homogeneous elements of $B$ to a system of linear equations
\begin{equation}
\label{equation:equation_for_faithful_flatness}
    \textstyle{ \sum_{1 \leq k \leq n} y_kc_{ki} = d_i \hspace{12mm}
      (1 \leq i \leq m), }
\end{equation}
where the $c_{ki}$ and $d_i$ are homogeneous elements of $A$, there exists a solution $(x_k)_{1 \leq k \leq n}$ consisting of homogeneous elements of $A$.

Now, in order to prove the lemma, we let $f \colon A \to B$ and $g' \colon C' \to D'$ be as in the statement. We note that, as a consequence of \cite[Addendum 3.23]{diracgeometry1}, the map $f \colon A \to B$ necessarily is a monomorphism, and therefore, we may assume that it is the inclusion of a sub-Dirac ring.

We will define subextensions $g_n \colon C_n \to D_n$ of $f \colon A \to B$ recursively for $n \geq 0$. We define $g_0 = g'$, and assuming that $g_n$ has been defined for $n < r$, we define $g_r$ as follows: We consider all systems of equations of the form~(\ref{equation:equation_for_flatness})~and~(\ref{equation:equation_for_faithful_flatness}) with $c_{ki}$ and $d_i$ homogeneous elements in $C_{r-1}$. Since $f \colon A \to B$ is flat, every solution $(y_k)$, consisting of homogeneous elements of $B$, to such a system of equations of the first kind can be written as a linear combination $y_k = \sum b_jz_{jk}$ with the $b_j$ and $z_{jk}$ homogeneous elements of $B$ and $A$, respectively. Similarly, since $f \colon A \to B$ is faithfully flat, for every solution $(y_k)$, consisting of homogeneous elements of $B$, to a system of equations of the second kind, there exists a solution $(x_k)$ to the same system of equations consisting of homogeneous elements of $A$. We now define $D_r \subset B$ to be the sub-Dirac ring obtained by adjoining to $D_{r-1}$ all the homogeneous elements $b_j$, $z_{jk}$, and $x_k$ of $B$ obtained in this way and define $g_r \colon C_r \to D_r$ to be the inclusion of $C_r = A \cap D_r$. This completes the definition of the subextensions $g_n \colon C_n \to D_n$ for $n \geq 0$.

We claim that their union $g \colon C \to D$ is the desired subextension of $f \colon A \to B$. Indeed, it is faithfully flat as both the above criteria hold for it by construction, and since $\kappa$ is uncountable, the Dirac rings $C$ and $D$ are $\kappa$-small.
\end{proof}

\begin{corollary}
\label{corollary:cofinality_of_ff_subexts_in_kappa_small_ones}
Let $f \colon A \to B$ be a faithfully flat map of Dirac rings and $\kappa$ an uncountable regular cardinal. Let $\SubExt_{\kappa}(f)$ be the poset of those subextensions $g \colon C \to D$ such that both $C$ and $D$ are $\kappa$-small, and let $\SubExt^{\ff}_{\kappa}(f)$ be the subposet of those $g \colon C \to D$ that, in addition, are faithfully flat. In this situation, the inclusion
$$\xymatrix{
{ \SubExt^{\ff}_{\kappa}(f) } \ar[r] &
{ \SubExt_{\kappa}(f) } \cr
}$$
is cofinal. In particular, the partially ordered set $\SubExt^{\ff}_{\kappa}(f)$ is $\kappa$-filtered.
\end{corollary}

\begin{proof}
The first part is a restatement of \cref{lem:approximation_of_kappa_small_subextensions_by_ff_ones}. The second part is a consequence of the first part, since $\SubExt_{\kappa}(f)$ is $\kappa$-filtered as it is a poset which admits $\kappa$-small upper bounds. 
\end{proof}

\begin{lemma}
\label{lemma:kappa_small_subextensions_cofinal_in_kappa_small_affines_over_tensor_product}
Let $f \colon A \to B$ be a faithfully flat map of Dirac rings, let $\kappa$ be an uncountable regular cardinal, and let $[n] \in \Delta_{+}$. The map
$$\xymatrix{
{ \SubExt_{\kappa}(f) } \ar[r] &
{ \CAlg(\Aff_{\kappa})_{/B^{\otimes_{A}[n]}} } \cr
}$$
from the poset of $\kappa$-small subextensions of $f$ to the category of $\kappa$-small Dirac rings over $B^{\otimes_{A} [n]}$ that to $g \colon C \to D$ assigns the induced map
$$\xymatrix{
{ D^{\otimes_{C} [n]} } \ar[r] &
{ B^{\otimes_{A} [n]} } \cr
}$$
is a $\varinjlim$-equivalence as a functor of $\infty$-categories. 
\end{lemma}

\begin{proof}
Using Joyal's version of Quillen's Theorem A, \cite[Theorem~2.19]{clausenjansen}, we must show that for any $k \colon R \to B^{\otimes_{A} [n]}$ with $R$ a $\kappa$-small Dirac ring, the slice category
$$\SubExt_{\kappa}(f)_{k/}$$
is weakly contractible. Explicitly, the objects of this category are pairs
$$(g \colon C \to D, h \colon R \to D^{\otimes_{C} [n]})$$
of a $\kappa$-small subextension of $f \colon A \to B$ and a map $h$ such that the composite
$$R \to D^{\otimes_{C} [n]} \to B^{\otimes_{A} [n]}$$
coincides with $k$. Maps are inclusions of subextensions that commute with the given maps from $R$. In particular, the category in question is a poset. We will show that it is a $\kappa$-filtered, which implies that it is weakly contractible. 

Since we are considering a poset, it suffices to show that every subset
$$S = \{(C_i \xrightarrow{\;g_i\;} D_i, R
  \xrightarrow{\;h_i\;} D_i^{\otimes _{C_i}[n]}) \mid i \in I\}$$
indexed by a $\kappa$-small index set $I$ has an upper bound. By first choosing an upper bound $g \colon C \to D$ of the subset
$$S' = \{C_i \xrightarrow{\;g_i\;} D_i \mid i \in I\} \subset
  \SubExt_{\lambda}(f),$$
we can assume that $C_i = C$ and $D_i = D$ for all $i \in I$.

With this accomplished, we consider the canonical comparison map
$$\xymatrix{
{ \varinjlim \prod_{i \in I} \Map(R,D^{\otimes_C[n]}) } \ar[r] &
{ \prod_{i \in I} \varinjlim \Map(R,D^{\otimes_C[n]}), } \cr
}$$
where the colimits are indexed by the partially ordered set $\SubExt_{\kappa}(f)$. Since this partially ordered set is $\kappa$-filtered, and since $I$ is $\kappa$-small, the map in question is an equivalence by~\cite[Proposition~5.3.3.3]{luriehtt}. Similarly, the map
$$\xymatrix{
{ \prod_{i \in I} \varinjlim \Map(R,D^{\otimes_C[n]}) } \ar[r] &
{ \prod_{i \in I} \Map(R,\varinjlim D^{\otimes_C[n]}) } \cr
}$$
is an equivalence, because $R$ is $\kappa$-small, and therefore, a $\kappa$-compact object in the category of Dirac rings. Finally, the map
$$\xymatrix{
{ \prod_{i \in I} \Map(R,\varinjlim D^{\otimes_C[n]}) } \ar[r] &
{ \prod_{i \in I} \Map(R,B^{ A[n]}) } \cr
}$$
is an equivalence, since $A \simeq \varinjlim C$ and $B \simeq \varinjlim D$, and since relative tensor products preserve filtered colimits in all three variables.

Considering the composite of these maps, we conclude that there exists a $\kappa$-small subextension $g' \colon C' \to D'$ of $f \colon A \to B$ that contains $g \colon C \to D$ as a subextension and a map of Dirac rings $h \colon R \to D'{}^{\otimes_{C'}[n]}$ such that for all $i \in I$, the map
$$\xymatrix{
{ R } \ar[r]^-{h_i} &
{ D^{\otimes_C[n]} } \ar[r] &
{ D'{}^{\otimes_{C'}[n]} } \cr
}$$
is equal to the map $h$. Thus, the pair $(g' \colon C' \to D', h \colon R \to D'{}^{\otimes_{C'}[n]})$ is the desired upper bound of the subset $S$.
\end{proof}

For every infinite cardinal $\kappa$, the flat topology on $\Aff$ induces a topology on the full subcategory $\Aff_{\kappa} \subset \Aff$. A presheaf $X \colon \Aff_{\kappa}^{\op} \to \spaces$ is a sheaf for this topology if and only if the conditions~(i)--(ii) of \cref{proposition:explicit_conditions_to_be_an_fpqc_sheaf} hold for all finite families $(S_i)_{i \in I}$ of objects in $\Aff_{\kappa}$ and for all faithfully flat maps $f \colon T \to S$ in $\Aff_{\kappa}$, respectively.  

\begin{theorem}
\label{theorem:waterhouse_preservation_of_flat_sheaves_under_lke}
For every uncountable regular cardinal $\lambda$, the left Kan extension
$$\xymatrix{
{ \presheaves(\Aff_{\lambda}) } \ar[r]^-{i_{!}} &
{ \presheaves(\Aff) } \cr
}$$
preserves sheaves with respect to the flat topology.
\end{theorem}

\begin{proof}
We extend Waterhouse's argument in~\cite[Theorem~5.1]{waterhouse1975basically} to the case of sheaves of anima. So we wish to show that if $X \in \presheaves(\Aff_{\lambda})$ satisfies the hypotheses~(i)--(ii) of \cref{proposition:explicit_conditions_to_be_an_fpqc_sheaf}, then so does $i_!(X) \in \presheaves(\Aff)$. 

First, if $X$ satisfies~(i), then by \cite[Lemma 5.5.8.14]{luriehtt}, it can be written as a sifted colimit of representable presheaves. Moreover, since $i_!$ preserves colimits and takes representable presheaves to representable presheaves, we deduce that also $i_!(X)$ can be written as a sifted colimit of representable presheaves, and hence, preserves finite products by~\cite[Proposition~5.5.8.10]{luriehtt}.

Second, suppose that $X$ satisfies~(ii) and let $f \colon A \to B$ be a faithfully flat map of Dirac rings. We must show that the induced map of anima
$$\xymatrix{
{ i_!(X)(\Spec(A)) } \ar[r] &
{ \varprojlim_{\,[n] \in \Delta} i_!(X)(\Spec(B^{\otimes_A[n]})) } \cr
}$$
is an equivalence. Now, for any $[n] \in \Delta_+$, the left Kan extension is given by
$$\textstyle{ i_!(X)(\Spec(B^{\otimes_A[n]})) \simeq \varinjlim 
  X(\Spec(R)), }$$
where the colimit is indexed by $R \to B^{\otimes_A[n]}$ in $((\Aff_{\lambda})^{\op})_{/B^{\otimes_A[n]}}$, and a combination of \cref{corollary:cofinality_of_ff_subexts_in_kappa_small_ones} and \cref{lemma:kappa_small_subextensions_cofinal_in_kappa_small_affines_over_tensor_product} shows that the functor 
$$\xymatrix{
{ \SubExt_{\lambda}^{\ff}(f) } \ar[r] &
{ ((\Aff_{\lambda})^{\op})_{/B^{\otimes_A[n]}} } \cr
}$$
that to the subextension $g \colon C \to D$ of $f \colon A \to B$ assigns $D^{\otimes_C[n]} \to B^{\otimes_A[n]}$ is a $\smash{ \varinjlim }$-equivalence. Hence, we can identify the map in question with the map
$$\xymatrix{
{ \varinjlim_{\,C \to D} X(\Spec(C)) } \ar[r] &
{ \varprojlim_{[n] \in \Delta} \varinjlim_{\,C \to D}  X(\Spec(D^{\otimes_C[n]})) } \cr
}$$
with the colimits indexed by $\SubExt_{\lambda}^{\ff}(f)$. Since the latter is $\lambda$-filtered, and since $\lambda$ is uncountable, we can commute the limit past the colimit, and therefore, it will suffice to show that the map
$$\xymatrix{
{ \varinjlim_{\,C \to D} X(\Spec(C)) } \ar[r] &
{ \varinjlim_{\,C \to D} \varprojlim_{[n] \in \Delta} X(\Spec(D^{\otimes_C[n]})) } \cr
}$$
is an equivalence. However, since every $g \colon C \to D$ in $\SubExt_{\lambda}^{\ff}(f)$ is faithfully flat, and since $X$ satisfies~(ii), this map is a colimit of equivalences, and hence, is itself an equivalence. This shows that $i_!(X)$ also satisfies~(ii) as desired.
\end{proof}

\begin{remark}
\label{rem: waterhouse_preservation_of_flat_sheaves_under_lke}
Let $\sheaves(\Aff_{\kappa}) \subset \presheaves(\Aff_{\kappa})$
be the full subcategory spanned by the sheaves for the
flat topology. \cref{theorem:waterhouse_preservation_of_flat_sheaves_under_lke}
shows that the canonical map
$$\xymatrix{
  { \varinjlim_{\lambda} \sheaves(\Aff_{\lambda}) } \ar[r] &
  { \sheaves(\Aff) } \cr
}$$
from the colimit indexed by the large partially ordered set of uncountable regular cardinals $\lambda$ in the universe of discourse of the large $\infty$-categories $\sheaves(\Aff_{\lambda})$ to the large $\infty$-category $\sheaves(\Aff)$ of accessible stacks on $\Aff$ is an equivalence.
\end{remark}

\begin{corollary}
\label{cor:sheafication_functor_exists_on_accessible_presheaves}
The inclusion of the full subcategory spanned by the sheaves for the flat topology admits a left adjoint sheafification functor,
$$\xymatrix@C=10mm{
{ \presheaves(\Aff) } \ar@<.7ex>[r]^-{\phantom{,}L\phantom{,}} &
{ \sheaves(\Aff), } \ar@<.7ex>[l]^-{\phantom{h}\iota\phantom{h}} \cr
}$$
and the sheafification functor $L$ preserves finite limits.
\end{corollary}

\begin{proof}
Given $X \in \presheaves(\Aff)$, by \cref{prop:equivalent_conditions_for_accessability} we can write $X \simeq i_!(X_0)$ for some cardinal $\lambda$ and $X_0 \in \presheaves(\Aff_{\lambda})$. Since $\Aff_{\lambda}$ is essentially small, the inclusion
\[
\iota_{\lambda} \colon \sheaves(\Aff_{\lambda}) \to \presheaves(\Aff_{\lambda})
\]
admits a left adjoint sheafification functor $L_{\lambda} \colon \presheaves(\Aff_{\lambda}) \to \sheaves(\Aff_{\lambda})$, because its domain and target are both presentable. Now, we claim that the map
$$\xymatrix@C=10mm{
{ X \simeq i_!(X_0) } \ar[r]^-{i_!\eta} &
{ i_!\iota_{\lambda}L_{\lambda}(X_0) } \cr
}$$
is a sheafification. Indeed, by \cref{theorem:waterhouse_preservation_of_flat_sheaves_under_lke}, the target is the underlying presheaf of a sheaf, and for every $Y \in \sheaves(\Aff)$, we have
$$\begin{aligned}
{} & \Map(i_!\iota_{\lambda}L_{\lambda}(X_0),\iota(Y)) \simeq \Map(\iota_{\lambda}L_{\lambda}(X_0),i^*\iota(Y)) \cr 
{} & \simeq \Map(\iota_{\lambda}L_{\lambda}(X_0),\iota_{\lambda}i^*(Y)) \simeq \Map(L_{\lambda}(X_0),i^*(Y)) \simeq \Map(X_0,\iota_{\lambda}i^*(Y)) \cr 
{} & \simeq \Map(X_0,i^*\iota(Y)) \simeq \Map(i_!(X_0),\iota(Y)) \simeq \Map(X,\iota(Y)) \cr
\end{aligned}$$
as desired. Hence, we conclude from  \cite[\href{https://kerodon.net/tag/02FV}{Tag 02FV}]{kerodon} that the desired left adjoint functor $L$ exists. It preserves finite limits, because each $L_{\lambda}$ does so.
\end{proof}

\begin{remark}
\label{rem:sheafication_as_localization_at_a_set_of_maps}
The sheafification functor $L \colon \presheaves(\Aff) \to \sheaves(\Aff)$ is a localization with respect to the large set $W$ consisting of the following colimit-interchange maps:
\begin{enumerate}
\item[{\rm (i)}]For every finite family $(S_i)_{i \in I}$ of objects in $\Aff$, the map
$$\xymatrix{
{ \coprod_{i \in I} h(S_i) } \ar[r] &
{ h(\coprod_{i \in I}S_i), } \cr
}$$
where $h \colon  \Aff \hookrightarrow \presheaves(\Aff)$ is the Yoneda embedding.
\item[{\rm (ii)}]For every faithfully flat map $f \colon T \to S$ in $\Aff$, the map
$$\xymatrix{
{ \varinjlim_{\Delta^{\op}} h(T^{\times_S[-]}) } \ar[r] &
{ h(S). } \cr
}$$
\end{enumerate}
Indeed, this follows from \cref{proposition:explicit_conditions_to_be_an_fpqc_sheaf}.
\end{remark}
  
We now show that the $\infty$-category of Dirac stacks satisfies the following variant on the $\infty$-categorical Giraud axioms.

\begin{theorem}
\label{theorem:giraud_axioms_for_dirac_stacks}
The $\infty$-category $\sheaves(\Aff)$ has the following properties:
\begin{enumerate}[leftmargin=8mm]
\item[{\rm (i)}]It is complete and cocomplete, and it is generated under small colimits by the essential image of the Yoneda embedding $h \colon \Aff \to \sheaves(\Aff)$, which consists of $\omega_1$-compact objects and is coaccessible.
\item[{\rm (ii)}]Colimits in $\sheaves(\Aff)$ are universal.
\item[{\rm (iii)}]Coproducts in $\sheaves(\Aff)$ are disjoint.
\item[{\rm (iv)}]Every groupoid in $\sheaves(\Aff)$ is       effective.
\end{enumerate}
\end{theorem}

\begin{proof} \cref{cor:sheafication_functor_exists_on_accessible_presheaves} shows that $\sheaves(\Aff)$ is a left exact localization of $\presheaves(\Aff)$, so completeness, cocompletness, the fact that $\sheaves(\Aff)$ is generated under small colimits by the Yoneda image, and parts (ii)--(iv) follow from \cref{theorem:giraud_axioms_for_accessible_presheaves}. It remains to prove that if $S$ is an affine Dirac scheme, then $h(S)$ is $\omega_1$-compact. Now, for every Dirac stack $X$, we have
$$\Map(h(S),X) \simeq X(S),$$
so the claim that $h(S)$ is $\omega_1$-compact is equivalent to the statement that $\omega_1$-filtered colimits in $\sheaves(\Aff)$ are calculated pointwise. By \cref{proposition:accessible_presheaves_are_complete_cocomplete}, all small colimits in $\presheaves(\Aff)$ are calculated pointwise, so it suffices to show that $\sheaves(\Aff) \subset \presheaves(\Aff)$ is closed under $\omega_1$-filtered colimits. But by \cref{proposition:explicit_conditions_to_be_an_fpqc_sheaf}, the sheaf condition is expressed in terms of $\omega_1$-small limits, so the pointwise colimit of an $\omega_1$-filtered diagram of sheaves is again a sheaf, because $\omega_1$-small limits and $\omega_1$-filtered colimits of anima commute by~\cite[Proposition~5.3.3.3]{luriehtt}, as $\omega_1$ is a regular cardinal.
\end{proof}

\begin{remark}
In the language of \cref{question:properties_of_macro_categories}, part~(i) of \cref{theorem:giraud_axioms_for_dirac_stacks} would be phrased as saying that $\sheaves(\Aff)$ is an $\omega_1$-macropresentable $\infty$-category.
\end{remark}

We next deduce from~\cite[Theorem~6.1.3.9]{luriehtt} the following fundamental descent theorem. We will use later to prove that $\infty$-categories of Dirac formal groups over Dirac stacks descend along effective epimorphisms.

\begin{theorem}
\label{theorem:slice_infty_categories_of_stacks_satisfy_descent}
The slice $\infty$-category functor
$$\xymatrix@C=20mm{
{ \sheaves(\Aff)^{\op} } \ar[r]^-{\sheaves(\Aff)_{/-}} &
{ \verylargecat_{\infty} } \cr
}$$
preserves small limits.
\end{theorem}

\begin{proof}
Suppose that $X \simeq \varinjlim_{k \in K}X_k$ is a colimit of a small diagram in $\sheaves(\Aff)$. We wish to prove that the canonical map
$$\xymatrix{
{ \sheaves(\Aff)_{/X} } \ar[r] &
{ \varprojlim_{k \in K^{\op}} \sheaves(\Aff)_{/X_k} } \cr
}$$
is an equivalence. As in \cref{rem: waterhouse_preservation_of_flat_sheaves_under_lke}, we may write
$$\textstyle{ \sheaves(\Aff) \simeq \varinjlim_{\lambda \in L} \sheaves(\Aff_{\lambda}) }$$
as the filtered colimit indexed by the partially ordered set $L$ consisting of the small uncountable regular cardinals $\lambda$ for which the diagram in equation factors through $\sheaves(\Aff_{\lambda}) \subset \sheaves(\Aff)$. Now, for every $\lambda \in L$, $\sheaves(\Aff_{\lambda})$ is an $\infty$-topos, so
$$\xymatrix{
{ \sheaves(\Aff_{\lambda})_{/X} } \ar[r] &
{ \varprojlim_{k \in K^{\op}} \sheaves(\Aff_{\lambda})_{/X_k} } \cr
}$$
is an equivalence by~\cite[Theorem~6.1.3.9]{luriehtt}. Therefore, also the induced map of colimits
$$\xymatrix{
{ \sheaves(\Aff)_{/X} \simeq \varinjlim_{\lambda \in L} \sheaves(\Aff_{\lambda})_{/X} } \ar[r] &
{ \varinjlim_{\lambda \in L} \varprojlim_{k \in K^{\op}} \sheaves(\Aff_{\lambda})_{/X_k} } \cr
}$$
is an equivalence, so it remains to argue that the canonical map
$$\xymatrix{
{ \varinjlim_{\lambda \in L}\varprojlim_{k \in K^{\op}} \sheaves(\Aff_{\lambda})_{/X_k} } \ar[r] &
{ \varprojlim_{k \in K^{\op}}\varinjlim_{\lambda \in L} \sheaves(\Aff_{\lambda})_{/X_k} } \cr
}$$
is an equivalence. But there exists a regular cardinal $\kappa$ such that $K$ is $\kappa$-small and $L$ is $\kappa$-filtered, so this follows from \cite[Proposition~5.3.3.3]{luriehtt} by the same argument as in the proof of~\cite[Lemma~4.5.3.1]{luriesag}.
\end{proof}

\subsection{Dirac schemes as stacks}

We show that that the category of Dirac schemes embeds fully faithfully into the $\infty$-category of Dirac stacks.

\begin{definition}
\label{definition:restricted_yoneda_embedding}
The restricted Yoneda embedding is the composition
$$\xymatrix{
{ \Sch } \ar[r] &
{ \Fun(\Sch^{\op},\spaces) } \ar[r] &
{ \Fun(\Aff^{\op},\spaces) } \cr
}$$
of the Yoneda embedding and the restriction along the canonical inclusion.
\end{definition}

If $X$ is a Dirac scheme, then we write $h(X)$ for its image by the restricted Yoneda embedding. We proceed to show that $h(X)$ is a Dirac stack.

\begin{lemma}
\label{lemma:presheaf_associated_to_affine_scheme_is_a_dirac_stack}
If $X$ is affine, then the presheaf $h(X)$ is a Dirac stack. 
\end{lemma}

\begin{proof}
The presheaf $h(X)$ is representable, and hence accessible, as the $\infty$-category $\Aff$ is coaccessible. The fact that $h(X)$ is a sheaf for flat topology is a restatement of faithfully flat descent, which we proved in \cite[Addendum 3.23]{diracgeometry1}.
\end{proof}

\begin{proposition}
\label{proposition:presheaf_represented_by_dirac_scheme_is_a_flat_sheaf}
If $X$ is a Dirac scheme, then the presheaf $h(X)$ is a sheaf for the flat topology.  
\end{proposition}

\begin{proof}
We must verify that $h(X)$ satisfies conditions~(i) and~(ii) of  \cref{proposition:explicit_conditions_to_be_an_fpqc_sheaf}. In the case of~(i), there is nothing to prove, and in the case of~(ii), we must show that if $h \colon T \to S$ is a faithfully flat map of affine Dirac schemes, then
$$\xymatrix{
{ h(X)(S) } \ar[r] &
{ \varprojlim_{\Delta} h(X)(T^{\times_S[-]}) } \cr
}$$
is an equivalence of anima. Now, since the presheaf $h(X)$ takes values in the full subcategory $\spaces_{\leq 0} \subset \spaces$ of $0$-truncated anima, the restriction
$$\xymatrix{
{ \varprojlim_{\Delta} h(X)(T^{\times_S[-]}) } \ar[r] &
{ \varprojlim_{\Delta_{\leq 1}} h(X)(T^{\times_S[-]}) } \cr
}$$
is an equivalence by~\cite[Proposition~A.1]{diracgeometry1}. So it suffices to show that
$$\xymatrix@C=10mm{
{ \Map(S, X) } \ar[r]^-{h^*} &
{ \Map(T, X) } \ar@<.7ex>[r]^-{d_0^*} \ar@<-.7ex>[r]_-{d_1^*} &
{ \Map(T \times_{S} T, X) } \cr
}$$
is a limit diagram of sets.

We first verify that $h^*$ is injective. So we let $f_1, f_2 \colon S \to X$ be maps of Dirac schemes with $f_1h = f_2h \colon T \to X$. Since $h$ induces a surjective map of underlying topological spaces by \cite[Proposition 3.15]{diracgeometry1}, we deduce that $f_1$ and $f_2$ have the same underlying map $p \colon |S| \to |X|$ of topological spaces. So we have $f_i = (p,\phi_i)$, and wish to show that the maps $\phi_1, \phi_2 \colon p^*\mathcal{O}_{X} \to \mathcal{O}_S$ of sheaves of Dirac rings coincide. This can be checked stalkwise. So we write $h = (q,\psi)$, let $s \in |X|$ with  $x = p(x)$, and choose $t \in |T|$ such that $s = q(t)$, where we use that $q \colon |T| \to |S|$ is surjective. Now, by the assumption that $f_1h = f_2h$, the two composites
$$\xymatrix{
{ \mathcal{O}_{X,x} } \ar[r]^-{\phi_{i,x}} &
{ \mathcal{O}_{S,s} } \ar[r]^-{\psi_s} &
{ \mathcal{O}_{T,t} } \cr
}$$
coincide. But the latter map is a flat local homomorphism of local Dirac rings, and therefore, it is injective by \cite[Proposition 3.16]{diracgeometry1}. This shows that $\phi_{1,x} = \phi_{2,x}$, and hence, that $f_1 = f_2$, as desired.

Next, we show that $g \colon T \to X$ with $gd_0 = gd_1 \colon T \times_ST \to X$ factors through a map $f \colon S \to X$. Such a factorization is necessarily unique, by what was proved above. We choose a covering $(W_i)_{i \in I}$ of $X$ by affine open subschemes, and consider the covering $(V_i)_{i \in I}$ of $T$ by the open subschemes $V_i \simeq T \times_XW_i \to T$. Since
$$|V_i \times_ST| = |T \times_SV_i| \subset |T \times_ST|,$$
we conclude from \cite[Proposition 3.24]{diracgeometry1} that there exists a unique covering $(U_i)_{i \in I}$ of $S$ by open subschemes such that $V_i \simeq T \times_SU_i \to T$ for all $i \in I$. We observe that any $f \colon S \to X$ with the desired property necessarily will satisfy $f(U_i) \subset W_i$. Thus, it suffices to show that for every $i \in I$, there exists a map of Dirac schemes $f_i \colon U_i \to W_i$ with the property that the composite map
$$\xymatrix@C=10mm{
{ V_i } \ar[r]^-{h|_{V_i}} &
{ U_i } \ar[r]^-{f_i} &
{ W_i } \ar[r] &
{ X } \cr
}$$
is equal to $g|_{V_i} \colon V_i \to X$.

Now, for every $i \in I$, we choose a covering $(U_{i,j})_{j \in J_i}$ of $U_i$ by affine open subschemes, and note that also $(V_{i,j})_{j \in J_i}$ with $V_{i,j} \simeq T \times_SU_{i,j} \to T$ is a covering of $V_i$ by affine open subschemes. Moreover, the map $h|_{V_{i,j}} \colon V_{i,j} \to U_{i,j}$ is a faithfully flat map of affine Dirac schemes, and we have a map
$$\xymatrix@C=10mm{
{ V_{i,j} } \ar[r]^-{g|_{V_{i,j}}} &
{ W_i } \cr
}$$
with the property its restriction along the two projections
$$\xymatrix{
{ V_{i,j} \times_{U_{i,j}}V_{i,j} } \ar@<.7ex>[r] \ar@<-.7ex>[r] &
{ V_{i,j} } \cr
}$$
coincide. As everything is affine, we conclude from \cref{lemma:presheaf_associated_to_affine_scheme_is_a_dirac_stack} that $g|_{V_{i,j}}$ admits a unique factorization through a map $f_{i,j} \colon U_{i,j} \to W_i$. By uniqueness, the maps $f_{i,j} \colon U_{i,j} \to W_i$ glue to give the desired $f_i \colon U_i \to W_i$. This completes the proof.
\end{proof}

\begin{lemma}
\label{lemma:scheme_a_colimit_of_its_basis_stable_under_intersection}
Let $X$ be a Dirac scheme, and let $\mathfrak{B}$ be a basis for the topology on its underlying space $|X|$ such that if $U,V \in \mathfrak{B}$, then $U \cap V \in \mathfrak{B}$. If $h(U)$ is accessible for all $U \in \mathfrak{B}$, then so is $h(X)$.
\end{lemma}

\begin{proof}
We know from \cref{proposition:presheaf_represented_by_dirac_scheme_is_a_flat_sheaf} that $h(X) \in \Fun(\Aff^{\op},\spaces)$ is a sheaf for the flat topology, and we wish to prove that it belongs to the full subcategory
$$\textstyle{ \varinjlim_{\kappa}\sheaves(\Aff_{\kappa}) \simeq \sheaves(\Aff) \subset \Fun(\Aff^{\op},\spaces). }$$
Since the inclusions $i_{\kappa} \colon \Aff_{\kappa} \to \Aff$ exhibit $\Aff$ as a colimit of the $\Aff_{\kappa}$, it suffices to show that there exists a $\lambda$ such that for all $\kappa \geq \lambda$, the diagram
$$\xymatrix@C=6mm{
{ \Aff_{\lambda} } \ar[rr] \ar[dr]_-(.35){i_{\lambda}^*h(X)} &&
{ \Aff_{\kappa} \ar[dl]^-(.35){i_{\kappa}^*h(X)}} \cr
{} &
{ \spaces } &
{} \cr
}$$
exhibits $i_{\kappa}^*h(X)$ as a left Kan extension of $i_{\lambda}^*h(X)$. Now, each $h(U)$ is accessible and $\mathfrak{B}$ is small, so we can find a $\lambda$ such that for all $\kappa \geq \lambda$, the diagram
$$\xymatrix@C=6mm{
{ \Aff_{\lambda} } \ar[rr] \ar[dr]_-(.35){i_{\lambda}^*h(U)} &&
{ \Aff_{\kappa} \ar[dl]^-(.35){i_{\kappa}^*h(U)}} \cr
{} &
{ \spaces } &
{} \cr
}$$
exhibits $i_{\kappa}^*h(U)$ as a left Kan extension of $i_{\lambda}^*h(U)$. Since left Kan extension preserves sheaves by \cref{theorem:waterhouse_preservation_of_flat_sheaves_under_lke}, we conclude that it will suffice to show that the map
$$\xymatrix{
{ \varinjlim_{U \in \mathfrak{B}}i_{\lambda}^*h(U) } \ar[r] &
{ i_{\lambda}^*h(X) } \cr
}$$
in $\sheaves(\Aff_{\lambda})$ is an equivalence.

To prove that this is so, let us write $h_{\lambda} \simeq i_{\lambda}^*h$. Since, locally on $S$, every map $\eta \colon S \to X$ from an affine Dirac scheme factors through some $U \in \mathfrak{B}$, the map
$$\xymatrix{
{ \coprod_{U \in \mathfrak{B}} h_{\lambda}(U) } \ar[r] &
{ h_{\lambda}(X) } \cr
}$$
in $\sheaves(\Aff_{\lambda})$ is an effective epimorphism. Therefore, it will suffice to show that the base-change
$$\xymatrix{
{ (\varinjlim_{U \in \mathfrak{B}}h_{\lambda}(U)) \times_{h_{\lambda}(X)}h_{\lambda}(U_0) } \ar[r] &
{ h_{\lambda}(U_0) } \cr
}$$
of the map in question along the map induced by the open immersion of any $U_0 \in \mathfrak{B}$ into $X$ is an equivalence. Since colimits in $\sheaves(\Aff_{\lambda})$ are universal, and since the Yoneda embedding preserves limits, the canonical map
$$\xymatrix{
{ \varinjlim_{U \in \mathfrak{B}}h_{\lambda}(U \times_XU_0) } \ar[r] &
{ (\varinjlim_{U \in \mathfrak{B}}h_{\lambda}(U)) \times_{h_{\lambda}(X)}h_{\lambda}(U_0) } \cr
}$$
in $\sheaves(\Aff_{\lambda})$ is an equivalence. Hence, we wish to show that the map of colimits
$$\xymatrix{
{ \varinjlim_{U \in \mathfrak{B}} h_{\lambda}(U \times_XU_0) } \ar[r] &
{ \varinjlim_{U \in \mathfrak{B}_{U_0}} h_{\lambda}(U) \simeq h_{\lambda}(U_0) } \cr
}$$
induced by the functor $- \times_XU_0 \colon \mathfrak{B} \to \mathfrak{B}_{/U_0}$, which exists by our assumption that the basis $\mathfrak{B}$ is closed under finite intersections, is an equivalence. In fact, we claim that this functor is a $\varinjlim$-equivalence. Indeed, by~\cite[Theorem~4.1.3.1]{luriehtt}, it suffices to show that
$$\mathfrak{B} \times_{\mathfrak{B}_{/U_0}}(\mathfrak{B}_{/U_0})_{V/}$$
is weakly contractible for every $V \in \mathfrak{B}_{/U_0}$. But we can identify this category with the poset $\mathfrak{B}_{V/}$, which has $V$ as a minimal element, so its classifying anima is contractible, as required.
\end{proof}

\begin{proposition}
\label{proposition:sheaf_represented_by_a_scheme_is_accessible_and_hence_a_stack}
If $X$ is a Dirac scheme, then the presheaf $h(X)$ is accessible.
\end{proposition}

\begin{proof}
The statement holds for affine Dirac schemes by \cref{lemma:presheaf_associated_to_affine_scheme_is_a_dirac_stack}, and we now use \cref{lemma:scheme_a_colimit_of_its_basis_stable_under_intersection} to bootstrap our way to general Dirac schemes.

First, if $X$ admits an open immersion $j \colon X \to S$ into an affine Dirac scheme, then we take $\mathfrak{B}$ to be the poset consisting of the open subschemes $U \subset X$ with the property that $j(U) = S_f \subset S$ is a distinguished open subscheme. This $\mathfrak{B}$ is a basis for the topology on $|X|$, it is stable under finite intersection, and $h(U)$ is accessible for every $U \in \mathfrak{B}$, because $S_f$ is affine. So $h(X)$ is accessible by  \cref{lemma:scheme_a_colimit_of_its_basis_stable_under_intersection}.

Finally, if $X$ is any Dirac scheme, then we take $\mathfrak{B}$ to be the poset of open subschemes $U \subset X$ with the property that $U$ admits an open immersion into some affine Dirac scheme. Again, this $\mathfrak{B}$ is a basis for the topology on $|X|$, it is stable under finite intersection, and by the case considered above $h(U)$ is accessible for all $U \in \mathfrak{B}$. So \cref{lemma:scheme_a_colimit_of_its_basis_stable_under_intersection} shows that $h(X)$ is accessible.
\end{proof}

\begin{theorem}
\label{theorem:presheaf_represented_by_a_scheme_is_a_dirac_stack_and_restricted_yoneda_is_fully_faithful}
The restricted Yoneda embedding defines a fully faithful functor
$$\xymatrix{
{ \Sch } \ar[r]^-{\phantom{,}h\phantom{,}} &
{ \sheaves(\Aff) } \cr
}$$
from the category of Dirac schemes to the $\infty$-category of Dirac stacks.
\end{theorem}

\begin{proof}
It follows from \cref{proposition:presheaf_represented_by_dirac_scheme_is_a_flat_sheaf} and \cref{proposition:sheaf_represented_by_a_scheme_is_accessible_and_hence_a_stack} that the restricted Yoneda embedding from \cref{definition:restricted_yoneda_embedding} takes values in $\sheaves(\Aff) \subset \Fun(\Aff^{\op},\spaces)$. So we get the desired functor and must prove that it is fully faithful. We let $X$ and $Y$ be Dirac schemes and wish to prove that the canonical map
$$\xymatrix{
{ \Map(Y,X) } \ar[r] &
{ \Map(h(Y),h(X)) } \cr
}$$
is an equivalence of anima. (In fact, the domain and target of this map are both $0$-truncated anima.) If $Y$ is affine, then this follows from the Yoneda lemma.

In general, if $\mathfrak{B}$ is a basis for the topology on the underlying space $|Y|$, which is stable under finite intersection, then the canonical map of Dirac schemes
$$\xymatrix{
{ \varinjlim_{U \in \mathfrak{B}} U } \ar[r] &
{ Y } \cr
}$$
is an equivalence by \cite[Theorem~2.32]{diracgeometry1}, and the canonical map of Dirac stacks
$$\xymatrix{
{ \varinjlim_{U \in \mathfrak{B}} h(U) } \ar[r] &
{ h(Y) } \cr
}$$
is an equivalence by the proof of \cref{lemma:scheme_a_colimit_of_its_basis_stable_under_intersection}. Thus, we can bootstrap our way from affine Dirac schemes to general Dirac schemes in the same way as we did in the proof of \cref{proposition:sheaf_represented_by_a_scheme_is_accessible_and_hence_a_stack}.
\end{proof}

A small diagram $X \colon K \to \Sch$ may not admit a colimit, and even if it does, then the colimit is generally not meaningful. Instead, it is the colimit of the composite diagram $h(X) \colon K \to \sheaves(\Aff)$, which always exists, that is meaningful. We show, however, that if $K$ is static, then the two colimits do agree.

\begin{addendum}
\label{addendum:schemes_are_stacks}
The functor $h \colon \Sch \to \sheaves(\Aff)$ preserves coproducts.
\end{addendum}

\begin{proof}
The functor in question preserves finite coproducts, so it suffices to show that for a small family of Dirac schemes $(X_k)_{k \in K}$, the canonical map
$$\xymatrix{
{ \varinjlim_{I \subset K}h(\coprod_{k \in I} X_k) } \ar[r] &
{ h(\varinjlim_{I \subset K} \coprod_{k \in I} X_k) \simeq h(\coprod_{k \in K}X_k), } \cr
}$$
where the colimits range over the filtered category of finite subsets $I \subset K$, is an equivalence. We claim that the left-hand colimit is calculated pointwise. Indeed, this claim is equivalent to the statement that the presheaf given by pointwise colimit satisfies the sheaf condition. But the flat topology on $\Aff$ if finitary and the diagram in question takes values in $\spaces_{\leq 0}$-valued presheaves, so by~\cite[Proposition~A.1]{diracgeometry1}, the sheaf condition amounts to the requirement that certain finite diagrams in $\spaces$ are limit diagrams. Hence, the claim follows from the fact that finite limits and filtered colimits of anima commute. Finally, given $T \in \Aff$, the canonical map
$$\xymatrix{
{ \varinjlim_{I \subset K} \Map(T,\coprod_{k \in I}X_k) } \ar[r] &
{ \Map(T,\coprod_{k \in K}X_k) } \cr
}$$
is an equivalence, since $T$ is quasi-compact.
\end{proof}

\subsection{Geometric maps}
\label{subsection:geometric_maps}

The $\infty$-category of all Dirac stacks is quite large, and in practice, one is often interested in the smaller subcategory of geometric Dirac stacks, which, informally, is the closure of the full subcategory of Dirac schemes under quotients by flat groupoids. The precise definition, however, is more delicate and follows Lurie~\cite[Chapter~26]{luriesagtemp} and T\"{o}en--Vezzosi~\cite{toenvezzosi}. Geometric Dirac stacks are formally similar to Artin stacks, the main differences being:
\begin{enumerate}
\item We consider sheaves of anima, as opposed to sheaves of groupoids.
\item We consider sheaves on affine Dirac schemes, as opposed to affine schemes.
\item We allow quotients by flat groupoids, as opposed to only smooth grouoids.
\end{enumerate}
More generally, we define the relative notion of a geometric map between Dirac stacks. We follow Lurie~\cite[Chapter~26]{luriesagtemp}, who sets up the analogous theory in derived algebraic geometry. The definition is recursive.

\begin{definition}
\label{definition:0_geometric_and_0_submersive_map}
Let $f \colon Y \to X$ be a map of Dirac stacks.
\begin{enumerate}
\item[(1)]The map $f \colon Y \to X$ is $0$-geometric if for every map $\eta \colon S \to X$ from a Dirac scheme, the base-change $f_S \colon Y_S \to S$ of $f$ along $\eta$ is equivalent to a map of Dirac schemes.
\item[(2)]The map $f \colon Y \to X$ is $0$-submersive if for every map $\eta \colon S \to X$ from a Dirac scheme, the base-change $f_S \colon Y_S \to S$ of $f$ along $\eta$ is an effective epimorphism and equivalent to a flat map of Dirac schemes.
\end{enumerate}
\end{definition}

\begin{remark}[Affine source suffices]
\label{remark:affine_source_suffices}
In Definition~\ref{definition:0_geometric_and_0_submersive_map}, it suffices to consider maps $\eta \colon S \to X$ with $S$ an affine Dirac scheme. However, the source of the base-change $f_S \colon Y_S \to S$ will typically not be an affine Dirac scheme, even if $S$ is affine.
\end{remark}

\begin{remark}[Effective epimorphisms]
\label{remark:effective_epis_of_fpqc_sheaves}
It follows from~\cref{proposition:characterization_of_effective_epimorphisms_of_sheaves} that a map of Dirac stacks $p \colon Y \to X$ is an effective epimorphism if and only if for every map $\eta \colon S \to X$ from an affine Dirac scheme, there exists a diagram
$$\xymatrix@C=10mm{
{ T } \ar[r]^{\tilde{\eta}} \ar[d]^-{q} &
{ Y } \ar[d]^-{p} \cr
{ S } \ar[r]^-{\eta} &
{ X } \cr
}$$
where $q$ is a faithfully flat map of affine Dirac schemes. So $p \colon Y \to X$ is an effective epimorphism if and only if, locally for the flat topology on the source, every map $\eta \colon S \to X$ from an affine Dirac scheme admits a factorization through $f$.
\end{remark}

\begin{remark}[Submersive maps of Dirac schemes]
\label{remark:submersive_maps_of_Dirac_schemes}
Let $f \colon Y \to X$ be a map of Dirac schemes. The map $f$ is automatically $0$-geometric, because the embedding of Dirac schemes in Dirac stacks preserves fiber products, and \cref{remark:effective_epis_of_fpqc_sheaves} shows that $f$ is $0$-submersive if and only if it is an $\fpqc$-covering in the sense that $f$ is flat and for every affine open $U \subset X$, there exists a quasi-compact open $V \subset f^{-1}(U)$ such that $f|_V \colon V \to U$ is surjective. This is the analogue for Dirac schemes of the definition of an $\fpqc$-covering given in~\cite[\href{https://stacks.math.columbia.edu/tag/022B}{Tag 022B}]{stacks-project}.
\end{remark}

\begin{definition}
\label{def:n_geometric_and_n_submersive_maps}
Let $f \colon Y \to X$ be a map of Dirac stacks and $n \geq 1$ an integer.
\begin{enumerate}
\item[(1)]The map $f \colon Y \to X$ is $n$-geometric if for every base-change $f_S \colon Y_S \to S$ along a map $\eta \colon S \to X$ from a Dirac scheme, there exists an $(n-1)$-submersive map $p \colon T \to Y_S$ from a Dirac scheme.
\item[(2)]The map $f \colon Y \to X$ is $n$-submersive if for every base-change $f_S \colon Y_S \to S$ along a map $\eta \colon S \to X$ from a Dirac scheme, there exists an $(n-1)$-submersive map $p \colon T \to Y_S$ from a Dirac scheme such that the composite map
$$\xymatrix{
{ T } \ar[r]^-{p} &
{ Y_S } \ar[r]^-{f_S} &
{ S } \cr
}$$
is $0$-submersive. 
\end{enumerate}
The map $f \colon Y \to X$ is geometric if it is $n$-geometric for some $n \geq 0$, and it is submersive if it is $n$-submersive for some $n \geq 0$. 
\end{definition}

\begin{remark}[Stability under base-change]
\label{remark:base-change}
In \cref{def:n_geometric_and_n_submersive_maps}, it suffices to consider maps $\eta \colon S \to X$ with $S$ an affine Dirac scheme. The properties of being $n$-geometric and being $n$-submersive are preserved under base-change along any map.
\end{remark}

\begin{warning}[The value of $n$ is irrelevant] 
\label{warning:only_geometricity_relevant}
Contrary to what \cref{def:n_geometric_and_n_submersive_maps} might suggest, it is the properties of a map of Dirac stacks of being geometric or submersive that are important, whereas the particular $n$ for which the map in question is $n$-geometric or $n$-submersive is of little relevance. However, the recursive nature of the definitions is useful for making inductive proofs.
\end{warning}

\begin{lemma}
\label{lem:n_geometric_means_nplusone_geometric}
Let $f \colon Y \to X$ be a map of Dirac stacks, and let $n \geq 0$ be an integer. If $f$ is $n$-geometric, then $f$ is also $(n+1)$-geometric. If $f$ is $n$-submersive, then $f$ is also $(n+1)$-submersive.
\end{lemma}

\begin{proof}
We argue by induction on $n \geq 0$. If $f$ is $0$-geometric, then, by definition, the base-change $f_S \colon Y_S \to S$ of $f$ along any map $\eta \colon S \to X$ from a Dirac scheme is equivalent to a map of Dirac schemes. Since the identity map $\id \colon Y_S \to Y_S$ is a $0$-submersive map from a Dirac scheme, so we conclude that $f$ is $1$-geometric.

So we let $n \geq 1$ and assume inductively that the statement has been proved for $k < n$. If $f$ is $n$-geometric, then for every base-change $f_S \colon Y_S \to S$ of $f$ along a map $\eta \colon S \to X$ from a Dirac scheme $S$, there exits an $(n-1)$-submersive map $p \colon T \to Y_S$ from a Dirac scheme. By the inductive hypothesis, the map $p$ is also $n$-submersive, which shows that $f$ is $(n+1)$-geometric. The same argument shows that if $f$ is $n$-submersive, then it is $(n+1)$-submersive.
\end{proof}

\begin{proposition}
\label{prop:submersiveeffectiveepimorphism}
If a map of Dirac stacks $f \colon Y \to X$ is $n$-submersive for some $n \geq 0$, then it is an effective epimorphism. 
\end{proposition}

\begin{proof}Since $\sheaves(\Aff)$ is generated under small colimits by the image of the Yoneda embedding, and since the property of being an effective epimorphism is local on the base, we may assume that $X$ is affine. We proceed by induction on $n \geq 0$, the case $n = 0$ being trivial, because a $0$-submersive map is an effective epimorphism by definition. So we let $n \geq 1$ and inductively assume the statement for $k < n$. Since $f \colon Y \to X$ is $n$-submersive and $X$ affine, there exists a factorization
$$\xymatrix{
{ T } \ar[r]^-{p} &
{ Y } \ar[r]^-{f} &
{ X } \cr
}$$
with $p$ an $(n-1)$-submersive map from a Dirac scheme such that $fp$ is $0$-submersive. Now, by the inductive hypothesis, both $p$ and $fp$ are effective epimorphisms, so we conclude from~\cite[Corollary~6.2.3.12]{luriehtt} that also $f$ is an effective epimorphism. 
\end{proof}

\begin{proposition}
\label{prop:n_geometric_maps_stability_under_composition}
Let $g \colon Z \to Y$ and $f \colon Y \to X$ be  maps of Dirac stacks, and let $n \geq 0$. If $f$ and $g$ are both $n$-geometric or both $n$-submersive, then so is $fg$. 
\end{proposition}

\begin{proof}
We argue by induction on $n \geq 0$. The geometric part of case $n = 0$ follows from definitions, and the submersive part follows from the fact that both effective epimorphisms and flat maps of Dirac schemes are closed under composition. So we let $n \geq 1$ and assume inductively that the statements have been proved for $k < n$. To prove the induction step, we consider the diagram
$$\xymatrix@C=10mm{
{ V } \ar[r]^-{q} &
{ Z_U } \ar[r]^-{p'} \ar[d]^-{g_U} &
{ Z_S } \ar[r]^-{\eta''} \ar[d]^-{g_S} &
{ Z } \ar[d]^-{g} \cr
{} &
{ U } \ar[r]^-{p} &
{ Y_S } \ar[r]^-{\eta'} \ar[d]^-{f_S} &
{ Y } \ar[d]^-{f} \cr
{} &&
{ S } \ar[r]^-{\eta} &
{ X } \cr
}$$
with $S$ a Dirac scheme and the squares cartesian. Suppose first that $f$ and $g$ are $n$-geometric. We first choose an $(n-1)$-submersive map $p \colon U \to Y_S$ from a Dirac scheme. Since $g$ is $n$-geometric, so is $g_U$, and hence, so we may choose an $(n-1)$-submersive map $q \colon V \to Z_U$ from a Dirac scheme. By the inductive hypothesis, the composition of the $(n-1)$-submersive maps $p'$ and $q$ is $(n-1)$-submersive. This shows that $fg$ is $n$-geometric. If $f$ and $g$ are $n$-submersive, then we may choose the maps $p$ and $q$ such that $f_Sp$ and $g_Uq$ are $0$-submersive. But then also 
$$(fg)_S(p'q) \simeq f_Sg_Sp'q \simeq (f_Sp)(g_Uq)$$
is $0$-submersive, so we conclude that $fg$ is $n$-submersive. 
\end{proof}

\begin{lemma}
\label{lemma:0_submersion_locally_we_have_lifts_along_effective_epimorphisms}
If $f \colon Y \to X$ is an effective epimorphism of Dirac stacks, then for every map $\eta \colon V \to X$ from a Dirac scheme, there exists a diagram
$$\xymatrix@C=10mm{
{ W } \ar[r] \ar[d]^-{p} &
{ Y } \ar[d]^-{f} \cr
{ V } \ar[r]^-{\eta} &
{ X } \cr
}$$
with $p \colon W \to V$ a $0$-submersion from a Dirac scheme $W$.
\end{lemma}

\begin{proof}
Let $(S_i)_{i \in I}$ be a family with $S_i \subset V$ affine open and $\bigcup_{i \in I}|S_i| = |V|$. It follows from \cref{remark:submersive_maps_of_Dirac_schemes} and \cref{addendum:schemes_are_stacks} that the induced map
$$\xymatrix{
{ S \simeq \coprod_{i \in I} S_i } \ar[r]^-{g} &
{ V } \cr
}$$
is $0$-submersive. Applying \cref{remark:effective_epis_of_fpqc_sheaves} to each $S_i$ separately, we find $0$-submersive maps $q_i \colon T_i \to S_i$ with $T_i$ affine and factorizations of the composite maps
$$\xymatrix{
{ T_i } \ar[r]^-{q_i} &
{ S_i } \ar[r] &
{ S } \ar[r]^-{g} &
{ V } \ar[r]^-{\eta} &
{ X } \cr
}$$
through $f \colon Y \to X$. Thus, we obtain the desired diagram with $W \simeq \coprod_{i \in I}T_i$.
\end{proof}

\begin{proposition}
\label{proposition:two_out_of_three_for_geometric_morphisms}
Let $g \colon Z \to Y$ and $f \colon Y \to X$ be  morphisms of Dirac stacks and $n \geq 0$.
\begin{enumerate}
\item[{\rm (a)}]Let $g$ be $n$-submersive. If $fg$ is $(n+1)$-geometric, then so is $f$.
\item[{\rm (b)}]Let $f$ be $(n+2)$-geometric. If $fg$ is $(n+1)$-geometric, so is $g$. 
\end{enumerate}
\end{proposition}

\begin{proof}
For part~(a), we consider the diagram
\[
\xymatrix@C=10mm{
{ V } \ar[r]^-{q} &
{ Z_S } \ar[r]^-{\eta''} \ar[d]^-{g_S} &
{ Z } \ar[d]^-{g} \cr
{} &
{ Y_S } \ar[r]^-{\eta'} \ar[d]^-{f_S} &
{ Y } \ar[d]^-{f} \cr
{} &
{ S } \ar[r]^-{\eta} &
{ X } \cr
}
\]
with $S$ a Dirac scheme and the squares cartesian. Since $fg$ is $(n+1)$-geometric, we can find an $n$-submersive map $q$ from a Dirac scheme, and since $g$ is $n$-submersive, so is its base-change $g_S$ and the composite map $g_Sq$ by \cref{prop:n_geometric_maps_stability_under_composition}. Hence, $g_Sq$ is the required $n$-submersive map to $Y_S$. 

For part~(b), we consider the diagram
$$\xymatrix@C=10mm{
{ W } \ar[r]^-{h} \ar@/_1pc/[dr]_-(.35){k} &
{ V } \ar[r]^-{q} &
{ Z_S } \ar[r]^-{\eta''} \ar[d]^-{g_S} &
{ Z } \ar[d]^-{g} \cr
{} &
{ U } \ar[r]^-{p} &
{ Y_S } \ar[r]^-{\eta'} \ar[d]^-{f_S} &
{ Y } \ar[d]^-{f} \cr
{} &&
{ S } \ar[r]^-{\eta} &
{ X } \cr
}$$
with $S$ a Dirac scheme and with the squares cartesian. The maps $p$ and $q$ are an $(n+1)$-submersive map from a Dirac scheme $U$ and an $n$-submersive map from a Dirac scheme $V$, and they exist by the assumption that $f$ and $fg$ be respectively $(n+2)$-geometric and $(n+1)$-geometric. Now, \cref{prop:submersiveeffectiveepimorphism} shows that $p$ is an effective epimorphism, so by \cref{lemma:0_submersion_locally_we_have_lifts_along_effective_epimorphisms}, we can find a $0$-submersive map of Dirac schemes $h \colon W \to V$ together with the indicated factorization of $g_S q h$ through $p$. 

Since $h$ is $0$-submersive, it is $n$-submersive by \cref{lem:n_geometric_means_nplusone_geometric}, and since also $q$ is $n$-submersive, it follows from \cref{prop:n_geometric_maps_stability_under_composition} that $qh$ is $n$-submersive. Similarly, since $k$ is $0$-geometric, it is $(n+1)$-geometric by \cref{lem:n_geometric_means_nplusone_geometric}, and since also $p$ is $(n+1)$-geometric, it follows from another application of \cref{prop:n_geometric_maps_stability_under_composition} that $pk$ is $(n+1)$-geometric. So part~(a) shows that $g_S$ is $(n+1)$-geometric, as
desired.
\end{proof}

\begin{proposition}
\label{proposition:geometricity_local_on_base}
Let $f \colon Y \to X$ be a map of Dirac stacks, and let $f' \colon Y' \to X'$ be the  base-change of $f$ along an effective epimorphism $p \colon X' \to X$.
\begin{enumerate}
\item[{\rm (a)}]If $f'$ is $n$-geometric with $n \geq 1$, then so is $f$.
\item[{\rm (b)}]If $f'$ is $n$-submersive with $n \geq 1$, then so is $f$.
\end{enumerate}
\end{proposition}

\begin{proof}
Let $\eta \colon S \to X$ be a Dirac scheme, and let $f_S \colon Y_S \to S$ be the base-change of $f$ along $\eta$. We wish to find an $(n-1)$-submersive map $U \to Y_S$ from a Dirac scheme. Using \cref{lemma:0_submersion_locally_we_have_lifts_along_effective_epimorphisms}, we can find the left-hand square
$$\xymatrix@C=10mm{
{ T } \ar[r]^-{\tilde{\eta}} \ar[d]^-{q} &
{ X' } \ar[d]^-{p} &
{ Y' } \ar[l]_-{f'} \ar[d]^-{p'} \cr
{ S } \ar[r]^-{\eta} &
{ X } &
{ Y } \ar[l]_-{f} \cr
}$$
with $q$ a $0$-submersive map of Dirac schemes. Since $f'$ is $n$-geometric, we can find an $(n-1)$-submersive map $r \colon U \to Y_T'$ from a Dirac scheme. But then
$$\xymatrix{
{ U } \ar[r]^-{r} &
{ Y_T' \simeq Y_T } \ar[r]^-{q_Y} &
{ Y_S } \cr
}$$
is $(n-1)$-submersive. Indeed, the map $q_Y$ is the base-change of the $0$-submersive map $q$, and hence, $0$-submersive, so by \cref{lem:n_geometric_means_nplusone_geometric} it is also $(n-1)$-submersive, and finally, the composite map is $(n-1)$-submersive by \cref{prop:n_geometric_maps_stability_under_composition}. This completes the proof of~(a), and the proof of~(b) is analogous.
\end{proof}

\begin{warning}
The statement of \cref{proposition:geometricity_local_on_base} does not hold for $n = 0$, since the property of being a Dirac scheme is not local for the flat topology. In fact, there exists an algebraic space $f \colon Y \to X$ over a scheme $X$ and an \'{e}tale surjection $p \colon Y' \to Y$ such that $Y'$ is a scheme, but $Y$ is not; see \cite[\S4.4.2]{vistoli2004notes}.
\end{warning}

\begin{lemma}
\label{lemma:n_submersive_morphism_of_schemes_is_0_submersive}
If a map of Dirac schemes $f \colon T \to S$ is $n$-submersive for some $n \geq 0$, then it is $0$-submersive. 
\end{lemma}

\begin{proof}
We proceed by induction on $n \geq 0$, the case $n = 0$ being trivial. So we let $n \geq 1$ and assume that the statement has been proved for $k < n$. Since $f \colon T \to S$ is $n$-submersive, there exists, by definition, a scheme $V$ and an $(n-1)$-submersive map $p \colon V \to T$ such that the composite map
$$\xymatrix{
{ V } \ar[r]^-{p} &
{ T } \ar[r]^-{f} &
{ S } \cr
}$$
is $0$-submersive. By the inductive hypothesis, the map $p$ is $0$-submersive, and since both $fp$ and $p$ are flat, we deduce that $f$ is flat as well. Similarly, since both $fp$ and $p$ are effective epimorphisms, so is $f$. Hence, $f$ is a $0$-submersion.
\end{proof}

\begin{proposition}
\label{proposition:n_geometric_and_submersive_implies_n_submersive}
Let $f \colon Y \to X$ be an $n$-geometric map of Dirac stacks. If $f$ is $m$-submersive for some $m \geq 0$, then $f$ is $n$-submersive.
\end{proposition}

\begin{proof}
Using~\cref{lem:n_geometric_means_nplusone_geometric} and downward induction, it suffices to show that if $f$ is $n$-geometric and $(n+1)$-submersive, then it is $n$-submersive. Since both properties are detected after base-change to a Dirac scheme, we may assume that $X$ is a Dirac scheme. If $n = 0$, then there is nothing to prove, so we assume that $n \geq 1$. Since $f$ is $n$-geometric, and since $X$ is a Dirac scheme, there exists an $(n-1)$-submersive map $p \colon U \to Y$ from a Dirac scheme. We claim that the composite map
$$\xymatrix{
{ U } \ar[r]^-{p} &
{ Y } \ar[r]^-{f} &
{ X } \cr
}$$
is a $0$-submersive map of Dirac schemes, so that $f$ is also $n$-submersive. 

To prove the claim, we consider the diagram
$$\xymatrix@C=10mm{
{ W } \ar[r]^-{r} &
{ V_U } \ar[r]^-{p'} \ar[d]^-{q_U} &
{ V } \ar[d]^-{q} &
{} \cr
{} &
{ U } \ar[r]^-{p} &
{ Y } \ar[r]^{f} &
{ X } \cr
}$$
with $q$ an $n$-submersive map from a Dirac scheme and the square cartesian. The map $q$ exists, because $f$ is $(n+1)$-submersive. The base-change $q_U$ is also $n$-submersive, so we can find the $(n-1)$-submersive map $r$ from a Dirac scheme such that $q_Ur$ is $0$-submersive. The composite map $(fp)(q_Ur)$ is an $n$-submersive map between Dirac schemes, and hence, $0$-submersive by \cref{lemma:n_submersive_morphism_of_schemes_is_0_submersive}. Thus, $(fp)(q_Ur)$ and $q_Ur$ are both flat maps of Dirac schemes, so $fp$ is flat, and $(fp)(q_Ur)$ and $q_Ur$ are both effective epimorphisms, so $fp$ is an effective epimorphism. This proves the claim.
\end{proof}

We have now established the notions of maps of Dirac schemes being geometric and submersive and proved that both are stable under composition and base-change as well as local on the base for the flat topology, so we proceed to specialize and consider the absolute notion of a geometric Dirac stack.

\begin{definition}
\label{def:geometric_stack_absolute}
A Dirac stack $X$ is $n$-geometric if the unique map
$$\xymatrix{
{ X } \ar[r]^-{f} &
{ \Spec(\mathbb{Z}) } \cr
}$$
is $n$-geometric. It is geometric if it is $n$-geometric for some $n \geq 0$. 
\end{definition}

\begin{remark}
In our definition, a Dirac stack is $0$-geometric if and only if it is represented by a Dirac scheme. Different choices are possible, but as we explained in \cref{warning:only_geometricity_relevant}, it is the property of being geometric that is important and not the property of being $n$-geometric for a particular $n \geq 0$.
\end{remark}

We record some basic stability properties of geometric Dirac stacks.

\begin{proposition}
\label{prop:n_geometric_maps_stability_under_compositionabsolute}
Let $f \colon Y \to X$ be a map of Dirac stacks, and let $n \geq 0$.
\begin{enumerate}
\item[{\rm (a)}]Suppose that $f$ is $n$-geometric. If $X$ is $n$-geometric, then so is $Y$.
\item[{\rm (b)}]Suppose that $f$ is $n$-submersive. If $Y$ is $(n+1)$-geometric, then so is $X$.
\item[{\rm (c)}]Suppose that $X$ is $(n+2)$-geometric. If $Y$ is $(n+1)$-geometric, then so is $f$.
\end{enumerate}
In particular, all maps between geometric Dirac stacks are geometric. 
\end{proposition}

\begin{proof}
This is \cref{prop:n_geometric_maps_stability_under_composition} and \cref{proposition:two_out_of_three_for_geometric_morphisms}. 
\end{proof}

\begin{proposition}
\label{proposition:geometric_stacks_closed_under_finite_limits}
The full subcategory $\sheaves(\Aff)^{\geom} \subset \sheaves(\Aff)$ spanned by the geometric stacks is closed under finite limits. More precisely, the final Dirac stack is $0$-geometric, and given a cartesian square of Dirac stacks
$$\xymatrix{
{ Y' } \ar[r]^-{g'} \ar[d]^-{f'} &
{ Y } \ar[d]^-{f} \cr
{ X' } \ar[r]^-{g} &
{ X } \cr
}$$
with $X'$, $X$, and $Y$ $n$-geometric, also $Y'$ is $n$-geometric.
\end{proposition}

\begin{proof}Only the claim concerning fiber products needs proof. If $n = 0$, then $X$, $X'$, and $Y$ are Dirac schemes, and hence so is $Y'$. If $n \geq 1$, then by \cref{prop:n_geometric_maps_stability_under_compositionabsolute}~(c), the map $f$ is $n$-geometric, and therefore, so is $f'$. Thus, \cref{prop:n_geometric_maps_stability_under_compositionabsolute}~(a) implies that $Y'$ is $n$-geometric as stated.
\end{proof}

\begin{theorem}
\label{theorem:geometricity_in_terms_of_groupoids}
The full subcategory $\sheaves(\Aff)^{\geom} \subset \sheaves(\Aff)$ spanned by the geometric Dirac stacks has the following properties:
\begin{enumerate}
\item[{\rm (1)}]It contains all small coproducts of affine Dirac schemes.
\item[{\rm (2)}]It is closed under colimits of groupoids with submersive face maps.
\end{enumerate}
Moreover, it is minimal with these properties.
\end{theorem}

\begin{proof}
We first show that $\sheaves(\Aff)^{\geom}$ satisfies~(1) and~(2). For~(1), we recall from \cref{addendum:schemes_are_stacks} that a small coproduct of affine Dirac schemes is a Dirac scheme, hence, a geometric Dirac stack. For~(2), we let $S \colon \Delta^{\op} \to \sheaves(\Aff)$ be a groupoid with colimit $X$. By \cref{theorem:giraud_axioms_for_dirac_stacks}, $S$ is effective in the sense that the square
$$\xymatrix@C=10mm{
{ S_1 } \ar[r]^-{d_1} \ar[d]^-{d_0} &
{ S_0 } \ar[d]^-{f} \cr
{ S_0 } \ar[r]^-{f} &
{ X } \cr
}$$
is cartesian. In this situation, we wish to show that if $d_0$ is a submersive map between geometric Dirac stacks, then $X$ is geometric. If we choose an $n \geq 0$ such that both $S_0$ and $S_1$ are $n$-geometric, then \cref{prop:n_geometric_maps_stability_under_compositionabsolute}~(c) shows that $d_0$ is $n$-geometric, and thus \cref{proposition:n_geometric_and_submersive_implies_n_submersive} shows that $d_0$ is $n$-submersive. Since (the horizontal) $f$ is an effective epimorphism, we conclude from \cref{proposition:geometricity_local_on_base} and from the fact that $d_0$ is $n$-submersive that (the vertical) $f$ is $n$-submersive. But then \cref{prop:n_geometric_maps_stability_under_compositionabsolute}~(b) shows that $X$ is $(n+1)$-geometric as desired.

Finally, we let $\mathcal{C} \subset \sheaves(\Aff)^{\geom}$ be a full subcategory satisfying~(1) and~(2) and show, by induction on $n \geq 0$, that $\mathcal{C}$ contains all $n$-geometric Dirac stacks. For $n = 0$, we must show that $\mathcal{C}$ contains every Dirac scheme $X$. We first suppose that $X \simeq \coprod_{i \in I}X_i$ is a coproduct of Dirac schemes, each of which can be embedded as an open subscheme $X_i \subset S_i$ of an affine Dirac scheme. We can write each $X_i$ as the union of a family $(T_{i,j})_{j \in J_i}$ of distinguished open subschemes $T_{i,j} \subset S_i$, and therefore, it follows from \cref{addendum:schemes_are_stacks} that the canonical map
$$\xymatrix{
{ V_0 \simeq \coprod_{i \in I}\coprod_{j \in J_i} T_{i,j} } \ar[r]^-{p} &
{ X } \cr
}$$
is $0$-submersive. Now, since $p$ is an effective epimorphism, its \v{C}ech nerve
$$\xymatrix{
{ \Delta_+^{\op} } \ar[r]^-{V} &
{ \sheaves(\Aff) } \cr
}$$
is a colimit diagram, and since coproducts in $\sheaves(\Aff)$ are disjoint, we have
$$\textstyle{ V_k \simeq V_0^{\times_X[k]} \simeq \coprod_{i \in I} \coprod_{j_0,\dots,j_k \in J_i} T_{i,j_0} \times_{S_i} \dots \times_{S_i} T_{i,j_k} }$$
for $k \geq 0$. But each summand is an affine Dirac scheme, since distinguished open subschemes of an affine Dirac scheme are closed under finite intersections, so we conclude from properties~(1) and~(2) that $X$ is contained in $\mathcal{C}$.

Next, if $X$ is a general Dirac scheme, then we choose a family $(S_i)_{i \in I}$ of affine open subschemes such that $\bigcup_{i \in I}|S_i| = |X|$. In this case, the canonical map
$$\xymatrix{
{ V_0 \simeq \coprod_{i \in I} U_i } \ar[r]^-{p} &
{ X } \cr
}$$
is an effective epimorphism, and its \v{C}ech nerve $V \colon \Delta_+^{\op} \to \sheaves(\Aff)$ has the property that for all $k \geq 0$, the Dirac scheme $V_k$ is a coproduct of open subschemes of affine Dirac schemes. Hence, each $V_k$ is contained in $\mathcal{C}$ by what was proved in the previous paragraph, and therefore, so is $X$. This proves the case $n = 0$.

Finally, we let $n \geq 1$ and assume, inductively, that $\mathcal{C}$ contains all $(n-1)$-geometric Dirac stacks. Given an $n$-geometric Dirac stack $X$, we choose an $(n-1)$-submersive map $p \colon V_0 \to X$ from a Dirac scheme and consider its \v{C}ech nerve
$$\xymatrix{
{ \Delta_+^{\op} } \ar[r]^-{V} &
{ \sheaves(\Aff). } \cr
}$$
Since $p$ is $(n-1)$-submersive and $V|_{\Delta^{\op}}$ an effective groupoid, we conclude that all face maps in $V|_{\Delta^{\op}}$ are $(n-1)$-submersive. Hence, by \cref{prop:n_geometric_maps_stability_under_composition}, so is every iterated face map $V_k \to V_0$. But $V_0$ is $(n-1)$-geometric, being a Dirac scheme, so \cref{prop:n_geometric_maps_stability_under_compositionabsolute} shows that $V_k$ is $(n-1)$-geometric for all $k \geq 0$. By the inductive hypothesis, $V|_{\Delta^{\op}}$ is a groupoid in $\mathcal{C}$ with submersive face maps, so we conclude from property~(2) that $X$ is contained in $\mathcal{C}$, as desired.
\end{proof}

\begin{remark}[Flat maps]
\label{remark:flat_maps_between_dirac_stacks}
The notion of a flat map of Dirac schemes may be generalized to the notion of a flat map of Dirac stacks as follows. A map of Dirac stacks $f \colon Y \to X$ is flat if, in every diagram of the form
$$\xymatrix@C+=10mm{
{ T } \ar[r]^-{h} \ar[d]^-{g} &
{ Y } \ar[d]^-{f} \cr
{ S } \ar[r]^-{p} &
{ X } \cr
}$$
with $S$ and $T$ Dirac schemes and with $p \colon S \to X$ and $(g,h) \colon T \to Y_S$ submersive, the map of Dirac schemes $g \colon T \to S$ is flat. One can show that a map of Dirac stacks is submersive if and only if it is a flat effective epimorphism. In particular, if a flat map of Dirac stacks admits a section, then it is submersive. It follows that, in \cref{theorem:geometricity_in_terms_of_groupoids}, the property~(2) is equivalent to:
\begin{enumerate}
    \item[(2')]It is closed under colimits of groupoids with flat face maps.
\end{enumerate}
Thus, the theorem shows that the full subcategory spanned by the geometric Dirac stacks is precisely the closure under colimits of groupoids with flat face maps of the full subcategory spanned by the Dirac schemes.
\end{remark}

\section{Coherent cohomology}
\label{sec:quasicoherentsheaves}

In this section, we construct the functor that to a Dirac stack $X$ assigns the presentably symmetric monoidal stable $\infty$-category $\QCoh(X)$ of a quasi-coherent $\mathcal{O}_X$-modules that we outlined in the introduction.

\subsection{Quasi-coherent $\mathcal{O}_X$-modules}

We first define the functor that to an affine Dirac scheme $S$ assigns its presentably symmetric monoidal stable $\infty$-category of quasi-coherent $\mathcal{O}_S$-modules. The definition uses animation, which we first recall following~\cite[Section~5.5.8]{luriehtt} and~\cite{cesnaviciusscholze}.

Let $\mathcal{C}$ be a presentable $1$-category. An object $P \in \mathcal{C}$ is compact $1$-projective if the map $\Map(P,-) \colon \mathcal{C} \to \Set$ corepresented by $P$ preserves small sifted colimits, and we write $\pi_0 \colon \mathcal{C}^{\cp} \to \mathcal{C}$ for the inclusion of the full subcategory spanned by the compact $1$-projective objects. Now, if $\mathcal{C}$ is generated under small sifted colimits by $\pi_0 \colon \mathcal{C}^{\cp} \to \mathcal{C}$, then its animation is defined to be the initial map
$$\xymatrix@C=10mm{
{ \mathcal{C}^{\cp} } \ar[r]^-{h} &
{ \mathcal{C}^{\an} } \cr
}$$
to an $\infty$-category that admits small sifted colimits. So any map $f_0 \colon \mathcal{C}^{\cp} \to \mathcal{D}$ to an $\infty$-category $\mathcal{D}$ that admits small sifted colimits admits a unique factorization
$$\xymatrix@C=6mm{
{ \mathcal{C}^{\cp} } \ar[rr]^{h} \ar[dr]_-(.35){f_0} &&
{ \mathcal{C}^{\an} } \ar[dl]^-(.32){f} \cr
{} &
{ \mathcal{D} } &
{} \cr
}$$
such that $f$ preserves sifted colimits. The $\infty$-category $\mathcal{C}^{\an}$ is presentable and the map $h \colon \mathcal{C}^{\cp} \to \mathcal{C}^{\an}$ preserves finite coproducts. Moreover, if also $\mathcal{D}$ admits all small colimits and if $f_0$ preserves finite coproducts, then $f$ preserves all small colimits. In particular, the inclusion $\pi_0 \colon \mathcal{C}^{\cp} \to \mathcal{C}$ factors uniquely as
$$\xymatrix@C=6mm{
{ \mathcal{C}^{\cp} } \ar[rr]^{h} \ar[dr]_-(.35){\pi_0} &&
{ \mathcal{C}^{\an} } \ar[dl]^-(.32){\pi} \cr
{} &
{ \mathcal{C} } &
{} \cr
}$$
and the map $\bar{\pi} \colon \tau_{\leq 0}(\mathcal{C}^{\an}) \to \mathcal{C}$ induced by $\pi$ is an equivalence, so that the right adjoint of $\pi$ can be identified with the inclusion $s \colon \mathcal{C} \to \mathcal{C}^{\an}$ of the full subcategory spanned by the $0$-truncated objects.

\begin{lemma}
\label{lem:animation_of_r_modules_exists_and_is_prestable}
If $R$ is a Dirac ring, then the animation $\Mod_R(\Ab)^{\an}$ of the abelian category $\Mod_R(\Ab)$ exists and is a complete Grothendieck prestable $\infty$-category. 
\end{lemma}

\begin{proof}
The compact $1$-projective objects in the abelian category $\mathcal{C} \simeq \Mod_R(\Ab)$ are the graded $R$-modules which are finitely generated and projective, so $\mathcal{C}^{\cp} \subset \mathcal{C}$ generates $\mathcal{C}$ under sifted colimits. Therefore, the animation $i \colon \mathcal{C}^{\cp} \to \mathcal{C}$ exists and agrees, by~\cite[Proposition~5.5.8.15]{luriehtt}, with the Yoneda embedding $h \colon \mathcal{C}^{\cp} \to \presheaves_{\Sigma}(\mathcal{C}^{\cp})$ into the full subcategory $\presheaves_{\Sigma}(\mathcal{C}^{\cp}) \subset \presheaves(\mathcal{C}^{\cp})$ spanned by the functors $(\mathcal{C}^{\cp})^{\op} \to \spaces$ that preserve finite products. Moreover, since $\mathcal{C}^{\cp}$ is additive, $\presheaves_{\Sigma}(\mathcal{C}^{\cp})$ is complete Grothendieck prestable by~\cite[Proposition~C.1.5.7 and Remark~C.1.5.9]{luriesag}.
\end{proof}

\begin{definition}
\label{definition:connective_qcoh_sheaves_on_an_affine}
If $S \simeq \Spec(R)$ is an affine Dirac stack, then
$$\QCoh(S)_{\geq 0} \simeq \Mod_{R}(\Ab)^{\an}$$
is the animation of the abelian category of graded $R$-modules, and
$$\xymatrix{
{ \QCoh(S)_{\geq 0} } \ar[r]^-{i} &
{ \QCoh(S) \simeq \Sp(\Mod_R(\Ab)^{\an}) } \cr
}$$
is its stabilization.
\end{definition}

\begin{warning}
While the inclusion $i \colon \QCoh(S)_{\geq 0 } \to \QCoh(S)$ is a stabilization for affine Dirac schemes, this will not be true for general Dirac stacks.
\end{warning}

To assure the reader that the definition in terms of animation is very natural, we record two equivalent definitions. For (3) below, we recall from~\cite[\S2.1]{diracgeometry1} that the abelian category of graded abelian groups is equivalent (as a symmetric monoidal category) to the heart of the Beilinson $t$-structure on graded spectra, so that any Dirac ring determines a commutative algebra in graded spectra.

\begin{proposition}
\label{proposition:various_description_of_qcoh_on_an_affine}
If $S \simeq \Spec(R)$ is an affine Dirac scheme, then the following $\infty$-categories are canonically equivalent:
\begin{enumerate}
    \item[{\rm (1)}]The $\infty$-category $\QCoh(X) \simeq \Sp(\Mod_R(\Ab)^{\an})$.
    \item[{\rm (2)}]The derived $\infty$-category $\mathcal{D}(\Mod_R(\Ab))$ of the Grothendieck abelian category of $R$-modules in graded abelian groups.
    \item[{\rm (3)}]The $\infty$-category $\Mod_{R}(\Fun(\mathbb{Z},\Sp))$ of $R$-modules in $\mathbb{Z}$-graded spectra.
\end{enumerate}
\end{proposition}

\begin{proof}The $\infty$-categories~(1)--(3) are all presentable and stable and equipped with right complete $t$-structures. So it suffices to compare connective parts, and by the universal property of animation, it suffices to verify that~(2) and~(3) are generated under sifted colimits by the subcategories spanned by compact projective objects, which coincide with $\Mod_R(\Ab)^{\cp}$. In the case of~(2), this is~\cite[Proposition~1.3.3.14]{lurieha}, and in the case of~(3), the heart of the Beilinson $t$-structure is given by
$$\Mod_R(\Fun(\mathbb{Z},\Sp))^{\heartsuit} \simeq \Mod_R(\Fun(\mathbb{Z},\Sp)^{\heartsuit}) \simeq \Mod_R(\Ab),$$
and the connective part $\Mod_R(\Fun(\mathbb{Z},\Sp)_{\geq 0}$ is generated under sifted colimits by subcategory spanned by the compact projective objects in the heart, which precisely is $\Mod_R(\Ab)^{\cp}$.
\end{proof}

The symmetric monoidal structure on $\Mod_R(\Ab)$ restricts to one on the full subcategory spanned by the compact 1-projective objects, which, in turn, gives rise to symmetric monoidal structures on $\QCoh(S)_{\geq 0}$ and $\QCoh(S)$, as well as on the Yoneda embedding $h$ and on the canonical inclusion $i$.

\begin{lemma}
\label{lemma:qcoh_affine_symmetric_monoidal}
If $S \simeq \Spec(R)$ is an affine Dirac scheme, then the functors
$$\xymatrix@C=8mm{
{ \Mod_R(\Ab)^{\cp} } \ar[r]^-{h} &
{ \QCoh(S)_{\geq 0} } \ar[r]^-{i} &
{ \QCoh(S) } \cr
}$$
promote uniquely to symmetric monoidal functors such that the tensor products on the middle and right-hand terms preserve colimits in each variable.
\end{lemma}

\begin{proof}
This is \cite[Propositions~4.8.1.10 and 2.2.1.9]{lurieha}.
\end{proof}

The symmetric monoidal category of compact 1-projective $R$-modules in graded abelian groups is covariantly functorial in $R$ through extension of scalars, and it follows from~\cite[Proposition~4.8.1.10]{lurieha} that this gives rise to a functor
$$\xymatrix{
{ \Aff^{\op} } \ar[r] &
{ \CAlg(\LPr)^{\Delta^1} } \cr
}$$
that to $S$ assigns the symmetric monoidal functor
$$\xymatrix@C=10mm{
{ \QCoh(S)_{\geq 0}^{\otimes} } \ar[r]^{i^{\otimes}} &
{ \QCoh(S)^{\otimes} } \cr
}$$
provided by \cref{lemma:qcoh_affine_symmetric_monoidal}. Since we were unable to find this argument spelled out in detail in the literature, we take this opportunity to do so here.

\begin{construction}
\label{construction:qcoh_as_a_functor_valued_in_calg_prl}
Let $\verylargecat{\!}_{\infty}^{\Sigma}, \verylargecat{\!}_{\infty}^L \subset \verylargecat_{\infty}$ be the subcategories spanned by the (possibly large) $\infty$-categories that admit finite coproducts and small colimits, respectively, and the functors between them that preserve the colimits in question. We recall from~\cite[Corollary~4.8.1.4]{lurieha} that the canonical inclusion
$$\xymatrix{
{ \verylargecat{\!}_{\infty}^{L} } \ar[r] &
{ \verylargecat{\!}_{\infty}^{\Sigma} } \cr
}$$
promotes to a lax symmetric monoidal functor such that the induced map
$$\xymatrix{
{ \CAlg(\verylargecat{\!}_{\infty}^{L}) } \ar[r] &
{ \CAlg(\verylargecat{\!}_{\infty}^{\Sigma}) } \cr
}$$
is the canonical inclusion of the $\infty$-category of symmetric monoidal $\infty$-categories, whose underlying $\infty$-categories admit small colimits and whose tensor products preserve small colimits in each variable, in the $\infty$-category of symmetric monoidal $\infty$-categories, whose underlying $\infty$-categories admit finite coproducts and whose tensor products preserve finite coproducts in each variable. In this situation, it follows from \cite[Proposition~4.8.1.10]{lurieha} that the diagram
$$\xymatrix{
{ \CAlg(\verylargecat{\!}_{\infty}^{L}) } \ar[r] \ar[d] &
{ \CAlg(\verylargecat{\!}_{\infty}^{\Sigma}) } \ar[d] \cr
{ \verylargecat{\!}_{\infty}^{L} } \ar[r] &
{ \verylargecat{\!}_{\infty}^{\Sigma} } \cr
}$$
is left adjointable. If $\mathcal{C}$ is small, then the left adjoint of the bottom horizontal map takes $\mathcal{C}$ to the $\infty$-category $\mathcal{P}_{\Sigma}(\mathcal{C})$ of finite product-preserving presheaves of small anima, which by \cite[Theorem~5.5.1.1]{luriehtt} is presentable. Thus, the left adjoint of the top horizontal map restricts to a map
$$\xymatrix{
{ \CAlg(\Cat_{\infty}^{\Sigma}) } \ar[r] &
{ \CAlg(\LPr) } \cr
}$$
from the $\infty$-category of small symmetric monoidal $\infty$-categories, whose underlying $\infty$-categories admit finite coproducts and whose tensor products preserve finite coproducts in each variable, to the $\infty$-category of presentably symmetric monoidal $\infty$-categories. It follows that the composition
$$\xymatrix{
{ \Aff^{\op} } \ar[r] &
{ \CAlg(\Cat_{\infty}^{\Sigma}) } \ar[r] &
{ \CAlg(\LPr) } \cr
}$$
of the functor that to $S \simeq \Spec(R)$ assigns $\Mod_R(\Ab)^{\cp,\otimes}$ and the restricted left adjoint is a functor that to $S$ assigns $\QCoh(S)_{\geq 0}^{\otimes}$. Finally, the map
$$\xymatrix@C=10mm{
{ \QCoh(S)_{\geq 0}^{\otimes} } \ar[r]^-{i^{\otimes}} &
{ \QCoh(S)^{\otimes} } \cr
}$$
is given by tensor product with $\Sp \in \LPr$, because $\QCoh(S)_{\geq 0}$ is presentable, and therefore, is functorial.
\end{construction}

We recall from  \cref{remark:universal_property_of_accessible_presheaves_as_cocompletion} that the $\infty$-category $\presheaves(\Aff)$ of Dirac prestacks is the free cocompletion of the $\infty$-category $\Aff$ of affine Dirac schemes. Therefore, the following definition is meaningful.

\begin{definition}
\label{definition:connective_qcoh_sheaf_functor_on_prestacks}
The unique small limit-preserving extension
$$\xymatrix@C=4mm{
{} &
{ \Aff^{\op} } \ar[dl]_-{h} \ar[dr]^-{\QCoh^{\otimes}} &
{} \cr
{ \presheaves(\Aff)^{\op} } \ar[rr]^-{\QCoh^{\otimes}} &&
{ \CAlg(\LPr) } \cr
}$$
of the functor provided by \cref{construction:qcoh_as_a_functor_valued_in_calg_prl} 
is said to be the functor that to a Dirac prestack $X$ assigns its presentably symmetric monoidal $\infty$-category of quasi-coherent $\mathcal{O}_X$-modules.
\end{definition}

If $f \colon Y \to X$ is a map of Dirac prestacks, then we write
$$\xymatrix{
{ \QCoh(X) } \ar@<.7ex>[r]^-{f^*} &
{ \QCoh(Y) } \ar@<.7ex>[l]^-{f_*} \cr
}$$
for the adjoint functors provided by \cref{definition:connective_qcoh_sheaf_functor_on_prestacks}. We call $f^*$ and $f_*$ the inverse image and the direct image, respectively. In fact, \cref{definition:connective_qcoh_sheaf_functor_on_prestacks} promotes $f^*$ to a symmetric monoidal functor, and therefore, by \cite[Corollary~7.3.2.7]{lurieha}, it also promotes $f_*$ to a lax symmetric monoidal functor.

\begin{warning}
\label{warning:symmetric_monoidal_category}
We will often abuse notation and simply write $\QCoh(X)$ both for the presentably symmetric monoidal $\infty$-category $\QCoh(X)^{\otimes}$ and for its underlying presentable $\infty$-category $\QCoh(X)$.
\end{warning}

Let $X$ be a Dirac prestack, and let $p \colon X \to \Spec(\mathbb{Z})$ be the unique map. Since the inverse image $p^*$ is symmetric monoidal, it follows in particular that
$$\mathcal{O}_X \simeq p^*(\mathbb{Z})$$
is the tensor unit of $\QCoh(X)$. More generally, its spin-$s$ Serre twist defined by
$$\mathcal{O}_X(s) \simeq p^*(\mathbb{Z}(s))$$
is a tensor-invertible object of $\QCoh(X)$.

\begin{definition}
\label{definition:structure_sheaf_and_serre_twist}
Let $X$ be a Dirac prestack, let $p \colon X \to \Spec(\mathbb{Z})$ be the unique map, and let $s$ be a half-integer. The spin-$s$ Serre twist of a quasi-coherent $\mathcal{O}_X$-module $\mathcal{F}$ is the quasi-coherent $\mathcal{O}_X$-module
$$\mathcal{F}(s) \simeq \mathcal{F} \otimes \mathcal{O}_X(s).$$
\end{definition}

\begin{remark}
\label{remark:concrete_description_of_qcoh}
We may describe $\QCoh(X)$ for a Dirac prestack $X$ more concretely as follows. By \cref{prop:equivalent_conditions_for_accessability}, we may write any prestack $X \in \presheaves(\Aff)$ as the colimit of the diagram given by the composite map
$$\xymatrix{
{ (\Aff_{\kappa})_{/X} } \ar[r] &
{ \presheaves(\Aff_{\kappa})_{/X} } \ar[r] &
{ \presheaves(\Aff_{\kappa}) } \ar[r] &
{ \presheaves(\Aff), } \cr
}$$
for some small cardinal $\kappa$. Hence, by definition, we can identify $\QCoh(X)$ with the limit of the diagram given by the composite map
$$\xymatrix@C+=12mm{
{ ((\Aff_{\kappa})_{/X})^{\op} } \ar[r] &
{ (\Aff_{\kappa})^{\op} } \ar[r] &
{ \Aff^{\op} } \ar[r]^-{\QCoh} &
{ \CAlg(\LPr). } \cr
}$$
\end{remark}

\begin{remark}
Let $f \colon Y \to X$ be a map of Dirac prestack, such that the induced map of Dirac stacks $L(f) \colon L(Y) \to L(X)$ is an equivalence. We will show in \cref{thm:qcohstackification} below that the $f^* \colon \QCoh(X) \to \QCoh(Y)$ is an equivalence.
\end{remark}

\subsection{The canonical $t$-structure}

We showed in \cref{lem:animation_of_r_modules_exists_and_is_prestable} that if $S \simeq \Spec(R)$ is an affine Dirac scheme, then the $\infty$-category
$$\QCoh(S)_{\geq 0} \simeq \Mod_R(\Ab)^{\an}$$
is complete Grothendieck prestable. So by \cite[Proposition~C.1.2.9]{luriesag}, its stabilization $\QCoh(S)$ has a canonical $t$-structure, which is both left complete and right complete. Moreover, its connective part is the essential image of the fully faithful functor
$$\xymatrix{
{ \QCoh(S)_{\geq 0} } \ar[r]^-{i} &
{ \QCoh(S), } \cr
}$$
and its coconnective part is stable under filtered colimits.

\begin{remark}
Under the equivalence $\QCoh(S) \simeq \mathcal{D}(\Mod_{R}(\Ab))$ from \cref{proposition:various_description_of_qcoh_on_an_affine}, the canonical $t$-structure corresponds to the $t$-structure on the derived $\infty$-category with $\mathcal{F}$ connective if and only if $\pi_i(\mathcal{F}) \in \Mod_{R}(\Ab)$ is zero for all $i < 0$. 
\end{remark}

\begin{proposition}
\label{prop:tstructure}
Let $X \in \presheaves(\Aff)$ be a Dirac prestack.
\begin{enumerate}
\item[{\rm (1)}]There exists a unique $t$-structure on $\QCoh(X)$ such that $\mathcal{F} \in \QCoh(X)_{\geq 0}$ if and only if $\eta^*(\mathcal{F}) \in \QCoh(S)_{\geq 0}$ for every map $\eta \colon S \to X$ with $S$ affine.
\item[{\rm (2)}]If $\mathcal{G} \in \QCoh(X)$ and $\eta^*(\mathcal{G}) \in \QCoh(S)_{\leq 0}$ for every map $\eta \colon S \to X$ with $S$ affine, then $\mathcal{G} \in \QCoh(X)_{\leq 0}$.
\item[{\rm (3)}]For every half-integer $s$, $\mathcal{O}_{X}(s) \in \QCoh(X)^{\heartsuit}$.
\end{enumerate}
\end{proposition}

\begin{proof}
To prove~(1), we use \cref{prop:equivalent_conditions_for_accessability} to write $X$ as a colimit of a small diagram $S \colon K \to \presheaves(\Aff)$ of affine Dirac stacks. Hence,
$$\textstyle{ \QCoh(X) \simeq \varprojlim_K \QCoh(S), }$$
and we wish to show that
$$\textstyle{ \QCoh(X)_{\geq 0} \simeq \varprojlim_K \QCoh(S)_{\geq 0} }$$
is the connective part of a $t$-structure. But $\QCoh(X)_{\geq 0}$ and $\QCoh(X)$ are both presentable, and therefore, it suffices by \cite[Proposition~1.4.4.11]{lurieha} to show that
$$\QCoh(X)_{\geq 0} \subset \QCoh(X)$$
is closed under small colimits and extensions, which is clear. This proves~(1).

To prove~(2), given $\mathcal{G} \in \QCoh(X)$ be as in the statement, we must show that for every $\mathcal{F} \in \QCoh(X)_{\geq 0}$, the mapping anima $\Map(\mathcal{F},\mathcal{G})$ is $0$-truncated. We again write $X \simeq \varinjlim_KS$ as a colimit of a small diagram $S \colon K \to \Aff$ of affine Dirac stacks, and for $k \in K$, we write $\eta_k \colon S_k \to X$ be the canonical map. This exhibits the mapping anima in question as a limit
$$\textstyle{ \Map(\mathcal{F},\mathcal{G}) \simeq \varprojlim_{k \in K} \Map(\eta_k^*(\mathcal{F}),\eta_k^*(\mathcal{G})) }$$
of anima, each of which is $0$-truncated, by assumption. But then also the limit is $0$-truncated, since the inclusion $\spaces_{\leq 0} \to \spaces$ preserves limits. This proves~(2).

Finally, we observe that~(3) holds, because
$\eta^*(\mathcal{O}_X(s)) \simeq \mathcal{O}_S(s) \in \QCoh(S)^{\heartsuit}$ for every map $\eta \colon S \to X$ with $S$ affine.
\end{proof}

\begin{remark}\label{rem:tstructure}The opposite implication of~(2) in \cref{prop:tstructure} fails in general. Indeed, the inverse image functor $\eta^* \colon \QCoh(X) \to \QCoh(S)$ need not in general  preserve coconnectivity, even if both $S$ and $X$ are affine.
\end{remark}

We recall that an exact functor
$F \colon \mathcal{C} \to \mathcal{D}$ between stable
$\infty$-categories with $t$-structures is defined to be right
$t$-exact, if $F(\mathcal{C}_{\geq 0}) \subset \mathcal{D}_{\geq 0}$,
and to be left $t$-exact, if
$F(\mathcal{C}_{\leq 0}) \subset \mathcal{D}_{\leq 0}$. It is $t$-exact, if it
is both left $t$-exact and right $t$-exact.

\begin{addendum}\label{add:tstructure}
If $f \colon Y \to X$ is a map of Dirac stacks, then the extension of scalars functor and the restriction of scalars functor
$$\xymatrix{
{ \QCoh(X) } \ar@<.7ex>[r]^-{f^*} &
{ \QCoh(Y) } \ar@<.7ex>[l]^-{f_*} \cr
}$$
are right $t$-exact and left $t$-exact, respectively.
\end{addendum}

\begin{proof}
Since $f^*$ and $f_*$ are adjoint, the statement that the former is right $t$-exact is equivalent to the statement that the latter is left $t$-exact. To show that $f^*$ is right $t$-exact, we let $\mathcal{F} \in \QCoh(X)_{\geq 0}$ and wish to show that $f^*(\mathcal{F}) \in \QCoh(Y)_{\geq 0}$. But if $\eta \colon S \to Y$ is map from an affine Dirac scheme, then so is $f \circ \eta \colon S \to X$, so
$$\eta^*(f^*(\mathcal{F})) \simeq (f \circ \eta)^*(\mathcal{F}) \in \QCoh(S)_{\geq 0},$$
since $\mathcal{F} \in \QCoh(X)_{\geq 0}$. 
\end{proof}

\begin{remark}[Comparison with derived $\infty$-category]
\label{remark:qcoh_and_derived_infty_categories}
If $X$ is any Dirac prestack, then there is a canonical exact functor
$$\xymatrix{
{ \mathcal{D}^b(\QCoh(X)^{\heartsuit}) } \ar[r] &
{ \QCoh(X) } \cr
}$$
from the bounded derived $\infty$-category of the abelian category of quasi-coherent $\mathcal{O}_X$-modules. In this generality, the functor and the bounded derived $\infty$-category are constructed in~\cite[Corollary~7.59]{bunke2019controlled}. However, it may not be fully faithful, even for qcqs schemes, as Verdier's counterexample in~\cite[Expos\'{e}~II, Appendice~I]{SGA6} shows. In fact, in general, the $\infty$-category $\mathcal{D}^b(\QCoh(X)^{\heartsuit})$ may not even be locally small, unlike the presentable $\infty$-category $\QCoh(X)$.
\end{remark}

\begin{remark}
If $S \simeq \Spec(R)$ is affine, then the canonical $t$-structure on $\QCoh(S) \simeq \mathcal{D}(\Mod_R(\Ab))$ is compatible with filtered colimits by \cite[Proposition~1.3.5.21]{lurieha}. We will show that this is true, more generally, if $X$ is geometric in the sense of \cref{def:geometric_stack_absolute}. In particular, for $X$ geometric, the $\infty$-category $\QCoh(X)_{\geq 0}$ is Grothendieck prestable in the sense of \cite[Definition~C.1.4.2]{luriesag}.
\end{remark}

\subsection{Descent for quasi-coherent $\mathcal{O}_X$-modules}
\label{subsection:descent_for_qcoh}

We proceed to show that the functor $\QCoh$ of \cref{definition:connective_qcoh_sheaf_functor_on_prestacks} admits a factorization through the sheafification functor $L \colon \presheaves(\Aff) \to \sheaves(\Aff)$, and we begin with the following lemma.

\begin{lemma}
\label{lem:flat_cartesian_squares_of_affines_are_right_adjointable}
Given a cartesian square of affine Dirac schemes
$$\xymatrix@C=10mm{
{ T' } \ar[r]^-{g'} \ar[d]^-{f'} &
{ T } \ar[d]^-{f} \cr
{ S' } \ar[r]^{g} &
{ S } \cr
}$$
in which either $f$ or $g$ or both are flat, the diagram
$$\xymatrix@C=12mm{
{ \QCoh(T') } &
{ \QCoh(T) } \ar[l]_-{\;g'{}^*} \cr
{ \QCoh(S') } \ar[u]_-{f'{}^*} &
{ \QCoh(S) } \ar[u]_-{f^*} \ar[l]_-{\;g^*} \cr
}$$
is right adjointable in the sense that the canonical map
$$\xymatrix{
{ f^*g_* } \ar[r] &
{ g_*'g'{}^*f^*g_* \simeq g_*'f'{}^*g^*g_* } \ar[r] &
{ g_*'f'{}^* } \cr
}$$
is an equivalence.
\end{lemma}

\begin{proof}
We write $S \simeq \Spec(A)$, $S' \simeq \Spec(A')$, $T \simeq \Spec(B)$, and $T' \simeq \Spec(B')$. The right adjoints $g_*$ and $g_*'$ exist, and since we consider affine Dirac schemes only, they are given by restriction of scalars, and hence, they preserve small colimits. Therefore, it suffices to show that for the generator $A'$ of the presentable stable $\infty$-category $\QCoh(S')$, the map $f^*g_*(A') \to g_*'f'{}^*(A')$ is an equivalence. But this map is given by the canonical projection
$$\xymatrix{
{ A' \otimes_A^LB } \ar[r] &
{ A' \otimes_AB \simeq B', } \cr
}$$
which is an equivalence by the assumption that either $f$ or $g$ or both are flat.
\end{proof}

\begin{remark}
If we built animated Dirac geometry by replacing Dirac rings by animated Dirac rings, then the analogue of \cref{lem:flat_cartesian_squares_of_affines_are_right_adjointable} would hold without any flatness assumptions.
\end{remark}

The following generalization to the relatively affine case will be useful later. Here, given a morphism $f \colon Y \rightarrow X$ of Dirac stacks, we say it is flat affine if for any map $\eta: S \rightarrow X$ from an affine scheme, the pullback $Y \times_{X} S$ is an affine scheme and the induced morphism $Y \times_{X} S \rightarrow S$ is flat. 

\begin{proposition}
\label{prop:pushforward_along_affine_flat_satisfies_base_change}
If $f \colon Y \to X$ is a flat affine map of Dirac stacks,  then the following hold:
\begin{enumerate}
\item[{\rm (1)}]The functor $f_* \colon \QCoh(Y) \to \QCoh(X)$ is $t$-exact and admits a right adjoint.
\item[{\rm (2)}]If $g \colon X' \to X$ is any map of Dirac stacks and if $f' \colon Y' \to X'$ is the base-change of $f$ along $g$, then
$$\xymatrix@C=10mm{
{ \QCoh(X) } \ar[r]^-{f^*} \ar[d]^-(.45){g^*} &
{ \QCoh(Y) } \ar[d]^(.45){g'{}^*} \cr
{ \QCoh(X') } \ar[r]^-{f'{}^*} &
{ \QCoh(Y') } \cr
}$$
is right adjointable in the sense that the base-change map
$$\xymatrix{
{ g^*f_* } \ar[r] &
{ f_*'f'{}^*g^*f_* \simeq f_*'g'{}^*f^*f_* } \ar[r] &
{ f_*'g'{}^* } \cr
}$$
is an equivalence.
\end{enumerate}
\end{proposition}

\begin{proof}
We first prove~(1). Since $f \colon Y \to X$ is affine, we have an adjunction
$$\xymatrix{
{ \Aff_{/Y} } \ar@<.7ex>[r]^-{u} &
{ \Aff_{/X} } \ar@<.7ex>[l]^-{v} \cr
}$$
with $u$ given by composition with $f$ and $v$ given by base-change along $f$. A counit $p$ of the adjunction induces a natural transformation
$$\xymatrix{
{ \QCoh(S) } \ar[r]^-{p^*} &
{ \QCoh(S \times_X Y) } \cr
}$$
of functors from $(\Aff_{/X})^{\op}$ to $\Cat_{\infty}$ given by extension of scalars along the canonical projection $p \colon S \times_X Y \to  S$, and we have a commutative diagram
$$\xymatrix@C=-16mm{
{ \varprojlim_{\,(S,S \to X)} \QCoh(S) } \ar[rr]^-{f^*} \ar[dr]_-(.45){p^*} &&
{ \varprojlim_{(T,T \to Y)} \QCoh(T) } \ar[dl] \cr
{} &
{ \varprojlim_{(S,S \to X)} \QCoh(S \times_XY) } &
{} \cr
}$$
with the top horizontal map and the right-hand slanted map the maps of limits induced by $u^{\op}$ and $v^{\op}$, respectively. Moreover, the latter map is an equivalence. Indeed, this follows from~\cite[Example~2.21]{luriehtt}, since $v$ admits a left adjoint, and hence, is a $\smash{ \varinjlim }$-equivalence. Hence, to prove that $f_*$ admits a right adjoint, we may instead prove that $p^*$ admits a right adjoint $p_*$ which itself admits a further right adjoint.

We may view $p^*$ as a map in $\LPr$, so it is clear that it has a right adjoint $p_*$, and to understand $p_*$, we employ~\cite[Proposition~4.7.4.19]{lurieha}. The assumption~(*) in loc.~cit.~that for every map $h \colon T \to S $ in $\Aff_{/X}$, the square
$$\xymatrix@C=12mm{
{ \QCoh(S) } \ar[r]^-{(p_S)^*} \ar[d]^-{h^*} &
{ \QCoh(S \times_XY) } \ar[d]^-{h'{}^*} \cr
{ \QCoh(T) } \ar[r]^-{(p_T)^*} &
{ \QCoh(T \times_XY) } \cr
}$$
be right adjointable holds by \cref{lem:flat_cartesian_squares_of_affines_are_right_adjointable}, since $f \colon Y \to X$ is flat affine. So we conclude that $p_*$ is the map of limits induced by a natural transformation
$$\xymatrix@C=12mm{
{ \QCoh(S \times_XY) } \ar[r]^-{(p_S)_*} &
{ \QCoh(S) } \cr
}$$
of functors from $(\Aff_{/X})^{\op}$ to $\verylargecat_{\infty}$.

Now, for every $(S,S \to X)$, we may identify $(p_{S})_{*}$ with the restriction of scalars along a map of Dirac rings, and hence, it is right $t$-exact. This implies that the induced map of limits $p_*$ is right $t$-exact, which shows that $f_*$ is right $t$-exact. But by \cref{add:tstructure}, $f_*$ is also left $t$-exact, so we conclude that $f_*$ is $t$-exact.

Similarly, for every $(S,S \to X)$, the functor $(p_S)_*$ admits a right adjoint given by coextension of scalars. Therefore, we may view $(p_S)_*$ as a natural transformation of functors from $(\Aff_{/X})^{\op}$ to $\LPr$, and since the forgetful functor
$$\xymatrix{
{ \LPr} \ar[r] &
{ \verylargecat_{\infty} } \cr
}$$
preserves limits, we conclude that $p_*$ admits a right adjoint.\footnote{\,The right adjoint of $p_*$ is typically not the map of limits induced by a natural transformation.} So $f_*$ admits a right adjoint, which completes the proof of~(1).

It remains to prove~(2). First, if $X'$ is affine, then with notation as in the proof of~(1) above, we may equivalently prove that the diagram
$$\xymatrix{
{ \varprojlim_{(S,S\to X)} \QCoh(S) } \ar[r]^-{p^*} \ar[d]^-{\pr_g} &
{ \varprojlim_{(S,S \to X)} \QCoh(S \times_XY) } \ar[d]^-{\pr_g} \cr
{ \QCoh(X') } \ar[r]^-{(p_{X'})^*} &
{ \QCoh(X' \times_XY) } \cr
}$$
is right adjointable. But this follows from~\cite[Proposition~4.7.4.19]{lurieha}, since for every $\eta \colon S \to X$ and $h \colon T \to S$ with $S$ and $T$ affine, the diagram
$$\xymatrix@C=10mm{
{ \QCoh(S) } \ar[r]^-{(p_S)^*} \ar[d]^-{h^*} &
{ \QCoh(S \times_XY) } \ar[d]^-{h'{}^*} \cr
{ \QCoh(T) } \ar[r]^-{(p_T)^*} &
{ \QCoh(T\times_XY) } \cr
}$$
is right adjointable by \cref{lem:flat_cartesian_squares_of_affines_are_right_adjointable}.

In general, we use that a map in $\QCoh(X')$ is an equivalence if and only if its image by the inverse image $h^* \colon \QCoh(X') \to \QCoh(X'')$ is an equivalence for every map $h \colon X'' \to X'$ with $X''$ affine. So we consider the diagram
$$\xymatrix@C=10mm{
{ \QCoh(X) } \ar[r]^-{f^*} \ar[d]^-{g^*} &
{ \QCoh(Y) } \ar[d]^-{g'{}^*} \cr
{ \QCoh(X') } \ar[r]^-{f'{}^*} \ar[d]^-{h^*} &
{ \QCoh(Y') } \ar[d]^-{h'{}^*} \cr
{ \QCoh(X'') } \ar[r]^-{f''{}^*} &
{ \QCoh(Y'') } \cr
}$$
with $f'$ the base-change of $f$ along $g$ and $f''$ the base-change of $f'$ along $h$. Since $X''$ is affine, both the outer square and the bottom square are right adjointable by what was already proved. Thus, in the diagram of canonical maps
$$\xymatrix@C=0mm{
{ h^*g^*f_* } \ar[rr] \ar[dr] &&
{ h^*f_*'g'{}^* } \ar[dl] \cr
{} &
{ f_*''h'{}^*g'{}^*, } \cr
}$$
the two slanted maps are equivalences, and hence, so is the top horizontal map. Therefore, we conclude that the canonical map
$$\xymatrix{
{ f^*g_* } \ar[r] &
{ f_*'g'{}^* } \cr
}$$
is an equivalence, as we wanted to proved.
\end{proof}

\begin{remark}
A useful way to interpret the second part of \cref{prop:pushforward_along_affine_flat_satisfies_base_change} is that if $f$ is flat affine, then $f_*$ can be calculated locally on the target. 
\end{remark}

\begin{remark}
A consequence of \cref{prop:pushforward_along_affine_flat_satisfies_base_change} is that coherent cohomology can be calculated by means of ``relative injective'' resolutions in the situation
$$\xymatrix@C=5mm{
{ Y } \ar[rr]^-{f} \ar[dr]_-(.4){q} &&
{ X } \ar[dl]^-(.4){p} \cr
{} &
{ S } &
{} \cr
}$$
with $f$ affine flat and an effective epimorphism, with $Y$ affine, and with $S \simeq \Spec(\mathbb{Z})$. The direct image $q_*$ is $t$-exact, because $q$ is a map between affine Dirac schemes, and $f_*$ is $t$-exact, by \cref{prop:pushforward_along_affine_flat_satisfies_base_change}. Therefore, given a resolution
$$\xymatrix@C=6mm{
{ 0 } \ar[r] &
{ \mathcal{F} } \ar[r] &
{ \mathcal{G}^0 } \ar[r] &
{ \mathcal{G}^1 } \ar[r] &
{ \cdots } \cr
}$$
in the abelian category $\QCoh(X)^{\heartsuit}$ such that the terms $\mathcal{G}^n$ belong to the essential image of $f_* \colon \QCoh(Y)^{\heartsuit} \to \QCoh(X)^{\heartsuit}$, the coherent cohomology
$$p_*(\mathcal{F}) \in \QCoh(S)_{\leq 0} \simeq \mathcal{D}(\Ab)_{\leq 0}$$
is represented by the complex
$$\xymatrix@C=6mm{
{ p_*(\mathcal{G}^0) } \ar[r] &
{ p_*(\mathcal{G}^1) } \ar[r] &
{ \cdots } \cr
}$$
in the abelian category $\QCoh(S)^{\heartsuit} \simeq \Ab$. In the homotopy theory literature, the objects in the essential image of $f_* \colon \QCoh(Y)^{\heartsuit} \to \QCoh(X)^{\heartsuit}$ are referred to as relative injective objects.
\end{remark}

We now proceed to show that the $\infty$-category of quasi-coherent $\mathcal{O}_X$-modules descends along effective epimorphisms.

\begin{proposition}
\label{proposition:test_for_flatness_in_terms_of_one_submersion_descent_for_qcoh}
If $f \colon T \to S$ is a faithfully flat map of affine Dirac schemes, then the diagram $\QCoh(T^{\times_S[-]}) \colon \Delta_+ \to \LPr$ is a limit diagram.
\end{proposition}

\begin{proof}
We let $U \simeq T^{\times_S[-]} \colon \Delta_+^{\op} \to \Aff$ be the \v{C}ech nerve of $f$, and wish to show that $\QCoh(U) \colon \Delta_+ \to \LPr$ is a limit diagram. By \cite[Lemma~6.2.3.7]{luriehtt}, we may instead show that its restriction $\QCoh(U)|_{\Delta_{s+}} \colon \Delta_{s+} \to \LPr$ along the inclusion of the (nonfull) subcategory $\Delta_{s+} \subset \Delta_+$ spanned by the injective order-preserving maps is a limit diagram. The advantage of restricting to $\Delta_{s+} \subset \Delta_+$ is that, in the diagram $U|_{\Delta_{s+}} \colon \Delta_{s+}^{\op} \to \Aff$, every map is a base-change of $f$, and hence, flat. So by \cref{prop:pushforward_along_affine_flat_satisfies_base_change}, for every map $\theta \colon [n] \to [m]$ in $\Delta_{s+}$, we have a diagram
$$\xymatrix{
{ \QCoh(U_m)_{\geq 0} } \ar@<.7ex>[r]^-{i} \ar[d]^-{\theta^*} &
{ \QCoh(U_m) } \ar@<.7ex>[l]^-{r} \ar@<.7ex>[r]^-{p} \ar[d]^-{\theta^*} &
{ \QCoh(U_m)_{\leq -1} } \ar@<.7ex>[l]^-{s} \ar[d]^-{\theta^*} \cr
{ \QCoh(U_n)_{\geq 0} } \ar@<.7ex>[r]^-{i} &
{ \QCoh(U_n) } \ar@<.7ex>[l]^-{r} \ar@<.7ex>[r]^-{p} &
{ \QCoh(U_n)_{\leq -1} } \ar@<.7ex>[l]^-{s} \cr
}$$
in which the rows are semi-orthogonal decompositions. It follows that the retractions $r$ define a map of diagrams from $\QCoh(U)|_{\Delta_{s+}}$ to the limit of the diagrams
$$\xymatrix{
{ \cdots } \ar[r]^-{\Omega} &
{ \QCoh(U)_{\geq 0}|_{\Delta_{s+}} } \ar[r]^-{\Omega} &
{ \QCoh(U)_{\geq 0}|_{\Delta_{s+}}, } \cr
}$$
and this map is an equivalence, by the definition of $\QCoh(U)$ as the stabilization of $\QCoh(U)_{\geq 0}$. Therefore, since limits commute with limits, it will suffice to show that the diagram $\QCoh(U)_{\geq 0}|_{\Delta_{s+}} \colon \Delta_{s+} \to \LPr$, or equivalently, the diagram
$$\xymatrix@C=18mm{
{ \Delta_+ } \ar[r]^-{\QCoh(U)_{\geq 0}} &
{ \verylargecat_{\infty} } \cr
}$$
is a limit diagram. To this end, we use that by (the dual of) Barr--Beck--Lurie theorem~\cite[Theorem~4.7.5.3]{lurieha}, it suffices to verify the following:
\begin{enumerate}
\item[(1)]The functor $f^* \colon \QCoh(U_{-1})_{\geq 0} \to \QCoh(U_{0})_{\geq 0}$ is comonadic.
\item[(2)]For every map $\theta \colon [m] \to n$ in $\Delta_+$, the diagram 
$$\xymatrix@C=10mm{
{ \QCoh(U_m)_{\geq 0} } \ar[r]^-{d_0^*} \ar[d]^-{\theta^*} &
{ \QCoh(U_{m+1})_{\geq 0} } \ar[d]^-{([0]*\theta)^*} \cr
{ \QCoh(U_n)_{\geq 0} } \ar[r]^-{d_0^*} &
{ \QCoh(U_{n+1})_{\geq 0} } \cr
}$$
is right adjointable.
\end{enumerate}
But~(1) holds by~\cite[Proposition~D.6.4.6]{luriesag}, since the map $f \colon T \to S$ is faithfully flat, and~(2) holds by  \cref{lem:flat_cartesian_squares_of_affines_are_right_adjointable}, because the map $d_0 \colon U_m \to U_{m+1}$ is flat.
\end{proof}

\begin{theorem}
\label{thm:qcohstackification}
Up to contractible choice, there exists a unique factorization 
$$\xymatrix@C=-2mm{
{ \presheaves(\Aff)^{\op} } \ar[rr]^-{L} \ar[dr]_-(.3){\QCoh} &&
{ \sheaves(\Aff)^{\op} } \ar[dl]^-(.3){\QCoh} \cr
{} &
{ \CAlg(\LPr) } &
{} \cr
}$$
through the sheafification functor $L$, and moreover, the functor
$$\xymatrix@C=12mm{
{ \sheaves(\Aff)^{\op} } \ar[r]^-{\QCoh} &
{ \CAlg(\LPr) } \cr
}$$
takes small colimits of Dirac stacks to limits of symmetric monoidal $\infty$-categories. 
\end{theorem}

\begin{proof}We must show that the functor $\QCoh \colon \presheaves(\Aff)^{\op} \to \CAlg(\LPr)$ takes the colimit-interchange maps in \cref{rem:sheafication_as_localization_at_a_set_of_maps} to equivalences of $\infty$-categories, and since the forgetful functor $\CAlg(\LPr) \to \LPr$ preserves and reflects limits, it suffices to show that the functor $\QCoh \colon \presheaves(\Aff)^{\op} \to \LPr$ does so. Thus, we must show that for every faithfully flat map $f \colon T \to S$ of affine Dirac schemes, the map
$$\xymatrix{
{ \QCoh(S) } \ar[r] &
{ \varprojlim_{\,\Delta} \QCoh(T^{\times_S[-]}) } \cr
}$$
is an equivalence, and that for every pair $(S_1,S_2)$ of affine Dirac schemes, the map
$$\xymatrix{
{ \QCoh(S_1 \sqcup S_2) } \ar[r] &
{ \QCoh(S_1) \times \QCoh(S_2) } \cr
}$$
is an equivalence.
The former statement follows from \cref{proposition:test_for_flatness_in_terms_of_one_submersion_descent_for_qcoh}, whereas the latter statement follows from the idempotent decomposition of graded modules over a product of two graded rings.
\end{proof}

\begin{remark}
We note that \cref{thm:qcohstackification} is equivalent to the statement that for every Dirac prestack $X$, the map induced by unit of the sheafification adjunction
$$\xymatrix{
{ \QCoh((\iota \circ L)(X)) } \ar[r]^-{\eta^*} &
{ \QCoh(X) } \cr
}$$
is an equivalence. In homotopy theory, this fact was first observed by Hopkins~\cite{hopkins1999complex}.
\end{remark}

\begin{corollary}
\label{cor:effectiveepimorphismqcoh}
If $p \colon Y \to X$ is an effective epimorphism in $\sheaves(\Aff)$, then the following hold:
\begin{enumerate}
    \item[{\rm(1)}]The diagram $\QCoh(Y^{\times_X[-]}) \colon \Delta_+ \to \CAlg(\LPr)$ is a limit diagram.
    \item[{\rm(2)}]The inverse image $p^* \colon \QCoh(X)
    \to \QCoh(Y)$ is conservative.
\end{enumerate}
\end{corollary}

\begin{proof}That $p \colon Y \to X$ is an effective epimorphism means, by definition, that the diagram $Y^{\times_X[-]} \colon \Delta_+^{\op} \to \sheaves(\Aff)$ is a colimit diagram. Thus, the assertion~(1) follows from \cref{thm:qcohstackification}. The assertion~(2) follows from~(1) and from the folklore fact, which we prove in \cite[Lemma B.1]{diracgeometry1}, that if $\mathcal{C} \colon \Delta \to \Cat_{\infty}$ is a cosimplicial $\infty$-category, then the canonical map $\varprojlim_{\,\Delta} \mathcal{C} \to \mathcal{C}^0$ is conservative.
\end{proof}

\subsection{Quasi-coherent $\mathcal{O}_X$-modules on geometric stacks} 
An advantage of geometric Dirac stacks is that the $t$-structure on the stable $\infty$-category of quasi-coherent $\mathcal{O}_X$-modules is particularly well-behaved, as we now show.

\begin{theorem}
\label{theorem:qcoh_on_geometric_dirac_stack_is_grothendieck_prestable}
If $X$ is a geometric Dirac stack, then the following hold:
\begin{enumerate}
\item[{\rm (1)}]The canonical $t$-structure on the presentable stable $\infty$-category $\QCoh(X)$ is left complete, right complete, and compatible with filtered colimits.
\item[{\rm (2)}]The presentable $\infty$-category $\QCoh(X)_{\geq 0}$ is complete Grothendieck prestable, and the inclusion $i \colon \QCoh(X)_{\geq 0} \to \QCoh(X)$ is a stabilization.
\item[{\rm (3)}]For $\mathcal{F} \in \QCoh(X)$, the following are equivalent:
\begin{enumerate}
\item[{\rm (a)}]\ $\mathcal{F} \in \QCoh(X)_{\leq 0}$.
\item[{\rm (b)}]There exists a submersion $p \colon V \to X$, where $V$ is a small coproduct of affine Dirac schemes, such that $p^*(\mathcal{F}) \in \QCoh(V)_{\leq 0}$.
\item[{\rm (c)}]For every submersion $p \colon V \to X$, where $V$ is a small coproduct of affine Dirac schemes, we have $p^*(\mathcal{F}) \in \QCoh(V)_{\leq 0}$.
\end{enumerate}
\end{enumerate}
\end{theorem}

\begin{proof}
We note~(1) and~(2) are equivalent by~\cite[Corollary~C.3.1.4]{luriesag}.
To prove the theorem, let us say that a geometric stack $X$ is \emph{good} if the statements~(1)--(3) hold for $X$. So we wish to prove that all geometric stacks are good, and by the minimality part of \cref{theorem:geometricity_in_terms_of_groupoids}, it suffices to verify the following:
\begin{enumerate}
\item[(i)]If $X \simeq \coprod_{\alpha} U_{\alpha}$ is a small coproduct of affine Dirac schemes, then $X$ is good.
\item[(ii)]If $X \colon \Delta^{\op} \to \sheaves(\Aff)$ is a groupoid in good geometric stacks with submersive face maps, then the colimit $\varinjlim_{[n] \in \Delta^{\op}} X_n$ is good. 
\end{enumerate}
To verify~(i), we let $X \simeq \coprod_{\alpha} U_{\alpha}$ be a small coproduct of affine Dirac schemes and first show $X$ satisfies~(1) in the statement. By the definition of $\QCoh$, we have
$$\textstyle{ \QCoh(X) \simeq \prod_{\alpha} \QCoh(U_{\alpha}), }$$
and by definition of the canonical $t$-structure, we have
$$\textstyle{ \QCoh(X)_{\geq 0} \simeq \prod_{\alpha} \QCoh(U_{\alpha})_{\geq 0}. }$$
Moreover, since a product of anima is contractible if and only if all of its factors are so, the coconnective part of the canonical $t$-structure is given by 
$$\textstyle{ \QCoh(X)_{\leq 0} \simeq \prod_{\alpha} \QCoh(U_{\alpha})_{\leq 0}. }$$
By \cref{proposition:various_description_of_qcoh_on_an_affine}, for $U_{\alpha} \simeq \Spec(R_{\alpha})$, we have an equivalence 
$$\QCoh(U_{\alpha}) \simeq \mathcal{D}(\Mod_{R_{\alpha}}(\Ab))$$
with the derived $\infty$-category of an abelian category with enough projectives, and hence, by~\cite[Propositions~1.3.3.16, 1.3.5.21, and 1.3.5.24]{lurieha}, the canonical $t$-structure on $\QCoh(U_{\alpha})$ is left and right complete and compatible with filtered colimits. But the properties of being left and right complete are expressed in terms of limits, and hence, are stable under products of $\infty$-categories, as is compatibility with filtered colimits by \cite[Proposition~C.3.2.4]{luriesag}. This shows that~(1) holds for $X$.

To prove that~(3) holds for $X$, let $\mathcal{F} \in \QCoh(X)$. It is clear that~(c) implies~(b), so we proceed to prove that~(a) implies~(c) and that~(b) implies~(a). We first assume that $\mathcal{F} \in \QCoh(X)_{\leq 0}$ and let $p \colon V \to X$ be a submersion with $V \simeq \coprod_{\beta} V_{\beta}$ a small coproduct of affine Dirac schemes. We wish to show that $p^*(\mathcal{F}) \in \QCoh(V)_{\leq 0}$, or equivalently, that $p_{\beta}^*(\mathcal{F}) \in \QCoh(V_{\beta})_{\leq 0}$ for all $\beta$, where $j_{\beta} \colon V_{\beta} \to V$ the canonical open immersion and $p_{\beta} = pj_{\beta}$. Since $V_{\beta}$ is quasi-compact, the map $p_{\beta} \colon V_{\beta} \to X$ factors as $p_{\beta} = i_{\beta}q_{\beta}$ with $i_{\beta} \colon X_{\beta} \to X$ the canonical open immersion of some finite sub-coproduct $X_{\beta} = \coprod_{\alpha \in I_{\beta}} U_{\alpha}$. So we have a diagram of Dirac schemes
$$\xymatrix@C=10mm{
{ V_{\beta} } \ar[r]^-{j_{\beta}} \ar[d]^-{q_{\beta}} &
{ V } \ar[d]^-{p} \cr
{ X_{\beta} } \ar[r]^-{i_{\beta}} &
{ X } \cr
}$$
with $p$ a submersion and with $j_{\beta}$ and $i_{\beta}$ open immersions. Now, we conclude from  \cref{lemma:n_submersive_morphism_of_schemes_is_0_submersive} that $p$ is flat, and hence, so is $pj_{\beta} \simeq i_{\beta}q_{\beta}$. Since $i_{\beta}$ is an open immersion, it follows that $q_{\beta}$ is flat. Since $q_{\beta}$ is a flat map between affine Dirac schemes, we conclude that $q_{\beta}^*$ is $t$-exact. But $i_{\beta}^*$ is also $t$-exact, because $i_{\beta}$ is the inclusion of a summand, and hence, so is $p_{\beta}^* \simeq q_{\beta}^*i_{\beta}^*$. We conclude that $p_{\beta}^*(\mathcal{F}) \in \QCoh(V_{\beta})_{\leq 0}$, as we wanted to prove. This shows that~(a) implies~(c).

Finally, to prove that~(b) implies~(a), we let $\mathcal{F} \in \QCoh(X)$ and let $p \colon V \to X$ be a submersion with $\smash{ V \simeq \coprod_{\beta} V_{\beta} }$ a small coproduct of affine Dirac schemes such that $p^*(\mathcal{F}) \in \QCoh(V)_{\leq 0}$. Now, for every $\alpha$, we conclude from \cref{lemma:n_submersive_morphism_of_schemes_is_0_submersive} and \cref{remark:effective_epis_of_fpqc_sheaves} that we have a diagram
$$\xymatrix@C=10mm{
{ W_{\alpha} } \ar[r]^-{j_{\alpha}} \ar[d]^-{q_{\alpha}} &
{ V } \ar[d]^-{p} \cr
{ U_{\alpha} } \ar[r]^-{i_{\alpha}} &
{ X } \cr
}$$
with $i_{\alpha}$ and $j_{\alpha}$ open immersions and with $q_{\beta}$ a faithfully flat map between affine Dirac schemes. Hence, in the diagram
$$\xymatrix@C=12mm{
{ \QCoh(X) } \ar[r]^-{(i_{\alpha}^*)} \ar[d]^-{p^*} &
{ \prod_{\alpha} \QCoh(U_{\alpha}) } \ar[d]^-{\prod_{\alpha}q_{\alpha}^*} \cr
{ \QCoh(V) } \ar[r]^-{(j_{\alpha}^*)} &
{ \prod_{\alpha} \QCoh(W_{\alpha}), } \cr
}$$
the top horizontal map is an equivalence, and the bottom horizontal map preserves (connectivity and) coconnectivity, whereas the right-hand vertical map preserves and reflects (connectivity and) coconnectivity. Since $p^*(\mathcal{F}) \in \QCoh(V)_{\leq 0}$, we find that $\mathcal{F} \in \QCoh(X)_{\leq 0}$, as desired. This shows that~(b) implies~(a), which completes the proof of~(i).

It remains to prove~(ii). We claim that if $f \colon Y \to X$ is any submersion of good geometric Dirac stacks and if $\mathcal{F} \in \QCoh(X)$, then $\mathcal{F} \in \QCoh(X)_{\leq 0}$ if and only $f^*(\mathcal{F}) \in \QCoh(Y)_{\leq 0}$. Indeed, if we choose a submersion $p \colon V \to Y$ with $V$ a small coproduct of affine Dirac stacks, then also $q \simeq f p \colon V \to X$ is a submersion, so~(3) shows that $\mathcal{F} \in \QCoh(X)_{\leq 0}$ if and only if $q^*(\mathcal{F}) \simeq p^*(f^*(\mathcal{F})) \in \QCoh(V)_{\leq 0}$ if and only if $f^*(\mathcal{F}) \in \QCoh(Y)_{\leq 0}$, as claimed.

Now, to prove~(ii), we let $X \colon \Delta^{\op} \to \sheaves(\Aff)$ be groupoid in good geometric stacks, whose face maps are submersions, and show that
$$\textstyle{ |X| \simeq \varinjlim_{[n] \in \Delta_{\phantom{s}}^{\op}}X_n \simeq \varinjlim_{[n] \in \Delta_s^{\op}} X_n }$$
is good. Here $\Delta_s \subset \Delta$ is the subcategory spanned by the injective order-preserving maps, and the right-hand equivalence follows from~\cite[Lemma~6.5.3.7]{luriehtt}. The claim that we proved above shows that the diagram $\QCoh(X)_{\geq 0} \colon \Delta_s \to \LPr$ admits a factorization through the canonical inclusion
$$\xymatrix{
{ \Groth_{\infty}^{\lex} } \ar[r] &
{ \LPr } \cr
}$$
of the (non-full) subcategory spanned by the Grothendieck prestable $\infty$-categories and the functors between them that preserve small colimits and finite limits. But Lurie proves in~\cite[Proposition~C.3.2.4]{luriesag} that this subcategory is closed under small limits, so we conclude that the $\infty$-category
$$\textstyle{ \QCoh(|X|)_{\geq 0} \simeq \varprojlim_{[n] \in \Delta_s} \QCoh(X_n)_{\geq 0} }$$
is Grothendieck prestable. It is also complete, since completeness is a property expressed in terms of limits and truncations, both of which are preserved by limits along left-exact functors, and it follows from~\cite[Corollary~C.3.2.5]{luriesag} that
$$\xymatrix{
{ \QCoh(|X|)_{\geq 0} } \ar[r]^-{i} &
{ \QCoh(|X|) } \cr
}$$
is a stabilization. This shows that~(2), or equivalently~(1), holds for $|X|$.

Finally, to see that also~(3) holds for $|X|$, let $p \colon V \to |X|$ be a submersion with $V$ a small coproduct of affine Dirac schemes. Since $f \colon X_0 \to |X|$ is an effective epimorphism, we may use \cref{lemma:0_submersion_locally_we_have_lifts_along_effective_epimorphisms} to obtain a diagram
$$\xymatrix@C=10mm{
{ W } \ar[r]^-{q} \ar[d]^-{g} &
{ X_0 } \ar[d]^-{f} \cr
{ V } \ar[r]^-{p} &
{ |X| } \cr
}$$
where both $g$ and $q$ are submersions from a small coproduct of affine Dirac schemes. Since $X_{0}$ and $V$ are good, we conclude that in the diagram
$$\xymatrix@C=10mm{
{ \QCoh(|X|) } \ar[r]^-{p^*} \ar[d]^-{f^*} &
{ \QCoh(V) } \ar[d]^-{g^*} \cr
{ \QCoh(X_0) } \ar[r]^-{q^*} &
{ \QCoh(W), } \cr
}$$
the maps $g^*$ and $q^*$ preserve and reflect coconnectivity, and so does $f^*$, being the stabilization of a cocontinuous left-exact functor between prestable $\infty$-categories. Therefore, also $p^*$ preserves and reflect coconnectivity, which shows that~(3) holds for $|X|$. This completes the proof of~(ii), and hence, of the theorem.
\end{proof}

\begin{remark}
\label{remark:flat_morphisms_of_dirac_stacks_preserve_coconnectivity}
As outlined in \cref{remark:flat_maps_between_dirac_stacks}, it is possible to generalize the notion of a flat map from Dirac schemes to geometric Dirac stacks. Moreover, one may show that if $f \colon Y \to X$ is a flat map between geometric Dirac stacks, then the inverse image functor $f^* \colon \QCoh(X) \to \QCoh(Y)$ is $t$-exact.
\end{remark}

\subsection{Affine maps and vector bundles}
\label{subsection:affine_maps_and_vector_bundles}

We discuss the relationship between quasi-coherent $\mathcal{O}_X$-modules and affine maps of stacks, which we now define.

\begin{definition}
\label{definition:affine_map}
A map of Dirac stacks $f \colon Y \to X$ is affine if the domain of its base-change $f_S \colon Y_S \to S$ along every map $\eta \colon S \to X$ from an affine Dirac scheme is an affine Dirac scheme.
\end{definition}

We write $\sheaves(\Aff)_{/X}^{\aff}$ for the full subcategory $\sheaves(\Aff)_{/X}$ spanned by the affine maps $f \colon Y \to X$.

\begin{example}
\label{example:affine_map}
A map of Dirac stacks $f \colon Y \to S$ to an affine Dirac scheme is affine if and only if $Y$ is an affine Dirac scheme. Hence, the Dirac spectrum functor
$$\xymatrix@C=12mm{
{ \CAlg(\QCoh(S)^{\heartsuit})^{\op} } \ar[r]^-{\Spec} &
{ \sheaves(\Aff)_{/S}^{\aff} } \cr
}$$
is an equivalence of $\infty$-categories.
\end{example}

It is clear that the property of being affine is preserved under base-change, and we now prove that it is local on the base.

\begin{proposition}
\label{proposition:being_affine_is_local_on_base}
If $p \colon X' \to X$ is an effective epimorphism of Dirac stacks, then
$$\xymatrix{
{ \sheaves(\Aff)_{/X}^{\aff} } \ar[r] \ar[d] &
{ \sheaves(\Aff)_{/X'}^{\aff} } \ar[d] \cr
{ \sheaves(\Aff)_{/X} } \ar[r] &
{ \sheaves(\Aff)_{/X'} } \cr
}$$
is a cartesian diagram of $\infty$-categories.
\end{proposition}

\begin{proof}We first suppose that $p \colon X' \to X$ is a faithfully flat map $q \colon T \to S$ of affine Dirac schemes and consider the expanded diagram
$$\xymatrix{
{ \sheaves(\Aff)_{/S}^{\aff} } \ar[r] \ar[d] &
{ \varprojlim_{[n] \in \Delta} \sheaves(\Aff)_{/T^{\times_S[n]}}^{\aff} } \ar[r] \ar[d] &
{ \sheaves(\Aff)_{/T}^{\aff}\phantom{.} } \ar[d] \cr
{ \sheaves(\Aff)_{/S} } \ar[r] &
{ \varprojlim_{[n] \in \Delta} \sheaves(\Aff)_{/T^{\times_S[n]}} } \ar[r] &
{ \sheaves(\Aff)_{/T}. } \cr
}$$
The right-hand square is cartesian, since the property of being affine is preserved under base-change, and the left horizontal maps are both equivalences. Indeed, for the bottom map, this follows from \cref{theorem:slice_infty_categories_of_stacks_satisfy_descent}, and for the top map, it follows from faithfully flat descent for modules~\cite[Proposition~1.4]{diracgeometry1} and \cref{example:affine_map}.

Next, in the general case, we must show that if $f \colon Y \to X$ is a map of Dirac stacks with the property that its base-change $f' \colon Y' \to X'$ along $p$ is affine, then $f$ is affine. So we let $\eta \colon S \to X$ be a map from an affine Dirac stack and must show that the base-change $f_S \colon Y_S \to X$ of $f$ along $\eta$ is affine. Since $p$ is an effective epimorphism, we conclude from \cref{proposition:characterization_of_effective_epimorphisms_of_sheaves} that there exists a diagram
$$\xymatrix{
{ T } \ar[r]^-{\xi} \ar[d]^-{q} &
{ X' } \ar[d]^-{p} \cr
{ S } \ar[r]^-{\eta} &
{ X } \cr
}$$
with $q$ a faithfully flat map of affine Dirac schemes. By what was already proved, it suffices to show that the base-change $f_T \colon Y_T \to T$ of $f$ along $\eta q$ is affine. But $f_T$ is also the base-change of $f'$ along $\xi$, and since $f'$ is affine, so is $f_T$.
\end{proof}

\begin{example}
\label{example:being_affine_is_local_on_the_base}
A geometric Dirac stack $X$ admits an effective epimorphism
$$\xymatrix{
{ S \simeq \coprod_{i \in I} S_i } \ar[r]^-{p} &
{ X } \cr
}$$
from a coproduct in $\sheaves(\Aff)$ of a family of affine Dirac schemes. Since coproducts in $\sheaves(\Aff)$ are universal, \cref{proposition:being_affine_is_local_on_base} shows that $f \colon Y \to X$ is affine if and only if its base-change $f_i' \colon Y_{S_i} \to S_i$ along $p_i \colon S_i \to X$ is affine for all $i \in I$.
\end{example}

Let $X$ be a Dirac stack, and let $\kappa$ be a regular cardinal such that $X$ is a left Kan extension along $\Aff_{\kappa} \to \Aff$. \cref{theorem:almost_topos} shows that the canonical map
$$\xymatrix{
{ \sheaves(\Aff)_{/X} } \ar[r] &
{ \varprojlim_{\eta} \sheaves(\Aff)_{/S} } \cr
}$$
to the limit indexed by $\eta \colon S \to X$ in $((\Aff_{\kappa})_{/X})^{\op}$ is an equivalence, and it follows directly from \cref{definition:affine_map} that this equivalence restricts to an equivalence of the respective full subcategories spanned by the affine maps.

\begin{definition}
\label{definition:relative_dirac_spectrum}
Let $X$ be a Dirac stack. The Dirac spectrum relative to $X$ is the left-hand vertical map in the essentially unique diagram
$$\xymatrix{
{ \CAlg(\QCoh(X)^{\heartsuit})^{\op} } \ar[r] \ar[d]^-{\,\Spec} &
{ \varprojlim_{\eta} \CAlg(\QCoh(S)^{\heartsuit})^{\op} } \ar[d]^-{\,\Spec} \cr
{ \sheaves(\Aff)_{/X}^{\aff} } \ar[r] &
{ \varprojlim_{\eta} \sheaves(\Aff)_{/S}^{\aff}, } \cr
}$$
whose right-hand vertical map is provided by the Yoneda embedding.
\end{definition}

\begin{remark}
\label{remark:relative_dirac_spectrum}
In \cref{definition:relative_dirac_spectrum}, the right-hand vertical and bottom maps are equivalence. The top horizontal map is a right adjoint, and therefore, so is the Dirac spectrum relative to $X$. Its left adjoint takes an affine map $f \colon Y \to X$ to the quasi-coherent $\mathcal{O}_X$-algebra $f_*^{\heartsuit}\mathcal{O}_Y \simeq f_*\mathcal{O}_Y$.
\end{remark}

\begin{warning}
\label{warning:relative_dirac_spectrum_not_always_an_equivalence}
In general, the Dirac spectrum relative to $X$, or equivalently, the top horizontal map in \cref{definition:relative_dirac_spectrum} is not an equivalence. If it were, then
$$X \mapsto \CAlg(\QCoh(X)^{\heartsuit})$$
would take colimits of Dirac stacks to limits of categories. To see that this is not the case, we let $k$ be a commutative ring, and let $X$ be the Dirac stack obtained by gluing two copies of the affine line $\mathbb{A}_k^1$ along the origin. Let $f \colon X \to X'$ be the unique map to the corresponding pushout in the category of Dirac schemes. It follows from~\cite[Theorem~16.2.0.2]{luriesag} that the map
$$\xymatrix{
{ \QCoh(\Spec(k[x,y]/(xy))_{\geq 0} \simeq \QCoh(X')_{\geq 0} } \ar[r]^-{f^*} &
{ \QCoh(X)_{\geq 0} } \cr
}$$
is an equivalence, and therefore, so is the induced map of hearts. Thus, it suffices to show that the diagram of categories
$$\xymatrix{
{ \CAlg(\Mod_{k[x,y]/(xy)}(\Ab)) } \ar[r]^-{g'{}^*} \ar[d]^-{f'{}^*} &
{ \CAlg(\Mod_{k[x]}(\Ab)) } \ar[d]^-{f^*} \cr
{ \CAlg(\Mod_{k[y]}(\Ab)) } \ar[r]^-{g^*} &
{ \CAlg(\Mod_k(\Ab)) } \cr
}$$
is not cartesian. If it were, then for every $k[x,y]/(xy)$-algebra $A$, the diagram
$$\xymatrix{
{ A } \ar[r] &
{ f_*'f'{}^*A \times g_*'g'{}^*A } \ar@<.7ex>[r] \ar@<-.7ex>[r] &
{ f_*'g_*g^*f'{}^* A \simeq g_*'f_*f^*g'{}^*A } \cr
}$$
would be a limit diagram. But for $A = k[x,y]/(xy,x-y) = k[t]/t^2$, this becomes
$$\xymatrix{
{ k[t]/t^2 } \ar[r] &
{ k \times k } \ar@<.7ex>[r] \ar@<-.7ex>[r] &
{ k, } \cr
}$$
which is not a limit diagram.
\end{warning}

\begin{theorem}
\label{theorem:relative_dirac_spectrum_equivalence_for_geometric_dirac_stacks}
If $X$ is a geometric Dirac stack, then
$$\xymatrix@C=12mm{
{ \CAlg(\QCoh(X)^{\heartsuit})^{\op} } \ar[r]^-{\Spec} &
{ \sheaves(\Aff)_{/X}^{\aff} } \cr
}$$
is an equivalence.
\end{theorem}

\begin{proof}
Let $\mathcal{C} \subseteq \sheaves(\Aff)^{\geom}$ be the full subcategory spanned by the geometric Dirac stacks $X$ for which the Dirac spectrum relative to $X$ is an equivalence. To show that $\mathcal{C} = \sheaves(\Aff)^{\geom}$, it suffices by \cref{theorem:geometricity_in_terms_of_groupoids} to show that
\begin{enumerate}
\item[{\rm (1)}]$\mathcal{C}$ contains all small coproducts of affine Dirac stacks, and that
\item[{\rm (2)}]$\mathcal{C}$ closed under colimits of groupoids with submersive face maps.    
\end{enumerate}
For (1), we observe that the source and target of the Dirac spectrum relative to $X$ both take coproducts of Dirac stacks to products of $\infty$-categories. Since the Dirac spectrum relative to an affine Dirac stack is an equivalence by \cref{example:affine_map}, we conclude that~(1) holds. For (2), we let $X \simeq \varinjlim_{\Delta^{\op}}S$ be a colimit of a groupoid in $\mathcal{C}$ with submersive face maps and consider the diagram
$$\xymatrix@C=10mm{
{ \CAlg(\QCoh(X)^{\heartsuit})^{\op} } \ar[r]^-{\Spec} \ar[d] &
{ \sheaves(\Aff)_{/X}^{\aff} } \ar[d] \cr
{ \varprojlim_{[n] \in \Delta} \CAlg(\QCoh(S_n)^{\heartsuit})^{\op} } \ar[r]^-{\Spec} &
{ \varprojlim_{[n] \in \Delta} \sheaves(\Aff)_{/S_n}^{\aff}. } \cr
}$$
Now, the left-hand vertical map is an equivalence by \cref{theorem:qcoh_on_geometric_dirac_stack_is_grothendieck_prestable}, 
the right-hand vertical map is an equivalence by \cref{theorem:slice_infty_categories_of_stacks_satisfy_descent} and \cref{proposition:being_affine_is_local_on_base}, and the lower horizontal map is an equivalence by assumption. Hence, also the top horizontal map is an equivalence, which shows that $X$ is contained in $\mathcal{C}$. So~(2) holds.
\end{proof}

Finally, we define the $\infty$-category of vector bundles on a Dirac stack. To this end, we consider the composition
$$\xymatrix@C=10mm{
{ (\QCoh(X)^{\heartsuit})^{\op} } \ar[r]^-{\Sym} &
{ \CAlg(\QCoh(X)^{\heartsuit})^{\op} } \ar[r]^-{\Spec} &
{ \sheaves(\Aff)_{/X}^{\aff} } \cr
}$$
of the symmetric algebra and the relative Dirac spectrum functors. We write
$$\xymatrix{
{ \mathbb{A}_X(\mathcal{F}) \simeq \Spec(\Sym_{\mathcal{O}_X}(\mathcal{F})) } \ar[r]^-{p} &
{ X } \cr
}$$
for the affine map given by the image of the quasi-coherent $\mathcal{O}_X$-module $\mathcal{F}$ under this composite functor and call it the affine space associated with $\mathcal{F}$.

The symmetric algebra has a natural grading, which gives rise to an action of the multiplicative group $\mathbb{G}_m$ on the affine space functor. It follows that the affine space functor promotes to a functor
$$\xymatrix@C=10mm{
{ (\QCoh(X)^{\heartsuit})^{\op} } \ar[r]^-{V_X'} &
{ \LMod_{\mathbb{G}_m}(\sheaves(\Aff)_{/X}^{\aff}) } \cr
}$$
that takes values in the $\infty$-category of affine Dirac stacks over $X$ with a left action by the underlying $\mathbf{E}_1$-group of $\mathbb{G}_m$. The functor $V_X'$ admits a left adjoint that to an affine map with $\mathbb{G}_m$-action $p \colon V \to X$ assigns the degree $1$ part of the quasi-coherent $\mathcal{O}_X$-algebra $p_*\mathcal{O}_V$. It is fully faithful if $X$ is geometric, but not in general.

\begin{definition}
\label{definition:perfect_of_Tor_amplitude_0}
Let $X$ be a Dirac stack. A quasi-coherent $\mathcal{O}_X$-module $\mathcal{E}$ is locally free of finite rank if for every map $\eta \colon S \to X$ with $S$ affine, $\eta^*(\mathcal{E}) \in \QCoh(S)$ is compact and in the heart of the $t$-structure.
\end{definition}

If a quasi-coherent $\mathcal{O}_X$-module $\mathcal{E}$ is locally free of finite rank, then we also say that $\mathcal{E}$ is perfect with $\Tor$-amplitude in $[0,0]$, and we write $$\QCoh(X)^{\loc} \subset \QCoh(X)$$
for the full subcategory spanned by the quasi-coherent $\mathcal{O}_X$-modules that are locally free of finite rank. It is contained in the heart of the $t$-structure.

\begin{example}
If $S \simeq \Spec(R)$ is affine, then $\QCoh(S)^{\loc}$ is canonically equivalent to the category of finitely generated projective graded $R$-modules.
\end{example}

Let $X$ be a Dirac stack. If we write $X \simeq \varinjlim_{\,\alpha} S_{\alpha}$ as a small colimit of affine Dirac schemes, then, by \cref{definition:perfect_of_Tor_amplitude_0}, the canonical map
$$\xymatrix@C=10mm{
{ \QCoh(X)^{\loc} } \ar[r] &
{ \varprojlim_{\,\alpha} \QCoh(S_{\alpha})^{\loc} } \cr
}$$
is an equivalence. It follows that the difficulties appearing in comparing affine maps and quasi-coherent $\mathcal{O}_X$-modules outlined in \cref{warning:relative_dirac_spectrum_not_always_an_equivalence} do not occur in the locally free of finite rank case, so that the functor $V_X'$ restricts to a fully faithful functor $$\xymatrix@C=10mm{
{ (\QCoh(X)^{\loc})^{\op} } \ar[r]^-{V_X} &
{ \LMod_{\mathbb{G}_m}(\sheaves(\Aff)_{/X}^{\aff}). } \cr
}$$

\begin{definition}
\label{definition:vector_bundle}
The $\infty$-category of vector bundles on $X$ is the full subcategory
$$\Vect(X) \subset \LMod_{\mathbb{G}_m}(\sheaves(\Aff)_{/X}^{\aff})$$
given by the essential image of the functor $V_X$.
\end{definition}

\begin{remark}
\label{remark:vector_bundle}
Since $\QCoh(X)^{\loc}$ and $\Vect(X)$ are canonically anti-equivalent, they are both equivalent to $1$-categories. The inverse of $V_X$ assigns to a vector bundle $(p \colon V \to X, \mu \colon \mathbb{G}_m \times V \to V)$ the degree $1$ part of the graded $\mathcal{O}_X$-algebra $p_*^{\heartsuit}\mathcal{O}_V$.
\end{remark}

\begin{remark}
\label{remark:vector_bundle_covariant}
The full subcategory $\QCoh(X)^{\loc} \subset \QCoh(X)^{\heartsuit}$ agrees with the full subcategory of dualizable $\mathcal{O}_X$-modules. In particular, it is canonically selfdual, so it is possible to identify $\QCoh(X)^{\loc}$ and $\Vect(X)$ via the covariant equivalence
$$\xymatrix{
{ \QCoh(X)^{\loc} \simeq (\QCoh(X)^{\loc})^{\op} } \ar[r]^-{V_X} &
{ \Vect(X). } \cr
}$$
However, we will follow Grothendieck~\cite[Definition 1.7.8]{egaII} and use the contravariant equivalence $V_X$ between $\QCoh(X)^{\loc}$ and $\Vect(X)$, which, among other things, leads to the correct geometric interpretation of the Lie algebra of a formal group.
\end{remark}

\section{Formal stacks and formal groups}

We extend the definition of formal groups to the setting of Dirac stacks. To do so, we follow Lurie~\cite{lurieell1,lurieell2} and define the $\infty$-category $\FGroup(X)$ of formal groups over a Dirac stack $X$ to the $\infty$-category of abelian group objects in the $\infty$-category $\Hyp(X)$ of formal hyperplanes over $X$, which we define first. We show that the $\infty$-category $\FGroup(X)$ descends along effective epimorphisms $f \colon X' \to X$.

\subsection{Formal completion and quasi-coherent $\mathcal{O}_X$-modules}

We define the formal completion of a Dirac scheme along a closed immersion and prove that, in the case of a regular closed immersion, the $\infty$-category of quasi-coherent modules on the formal completion participates in a ``local cohomology'' stable recollement.

\begin{definition}
\label{definition:formal_completion_along_closed_immersion}
Let $g \colon Z \to X$ be a closed immersion of Dirac schemes defined by a quasi-coherent ideal $\mathcal{I} \subset \mathcal{O}_X$. For each $m \geq 0$, the $m$-th infinitesimal neighbourhood is the Dirac scheme given by 
$$\xymatrix@C=10mm{
{ Z^{(m)} \simeq \Spec(\mathcal{O}_X/\mathcal{I}^{m+1}) } \ar[r]^-{g^{(m)}} &
{ X. } \cr
}$$
The formal completion of $X$ along $g$ is the map of Dirac stacks
$$\xymatrix{
{ Y \simeq \varinjlim_{\,m}Z^{(m)} } \ar[r]^-{j} &
{ X } \cr
}$$
from the colimit in $\sheaves(\Aff)$ of the infinitesimal neighborhoods.
\end{definition}

\begin{remark}
\label{remark:filtered_colimits_in_sheaves_levelwise_over_affines}
The colimit in \cref{definition:formal_completion_along_closed_immersion} agrees with the colimit calculated in the $\infty$-category $\presheaves(\Aff)$, or equivalently, the latter colimit is automatically a sheaf for the flat topology on $\Aff$. Indeed, since the Dirac stacks $Z^{(m)}$ are $N$-truncated for a fixed $N$, namely, $N = 0$, and since the flat topology on $\Aff$ is finitary, the sheaf condition amounts to the statement that a finite diagram is a limit diagram. But this condition is preserved under filtered colimits by \cite[Proposition~5.3.3.3]{luriehtt}.
\end{remark}

In Grothendieck's philosophy, the formal completion of $X$ along $g \colon Z \to X$ plays the role of the non-existing tubular neighborhood of $Z$ in $X$. So we should expect that $j \colon Y \to X$ behave as an open and affine immersion with $i \colon U \simeq X \smallsetminus Z \to X$ as its complement. The following result bears out this expectation.

\begin{theorem}\label{thm:recollement}
Let $X$ be a Dirac scheme and let $\eta \colon Z \to X$ be a regular closed immersion. Let $j \colon Y \to X$ be the formal completion of $X$ along $\eta \colon Z \to X$, and let $i \colon U \simeq X \smallsetminus Z \to X$ be the inclusion of the open complement of $S$. In this situation, there is a stable recollement
\vspace{-1mm}
$$\begin{xy}
(0,0)*+{ \QCoh(U) };
(30.5,0)*+{ \QCoh(X) };
(61,0)*+{ \QCoh(Y) };
(7,0)*+{}="1";
(23,0)*+{}="2";
(38,0)*+{}="3";
(54,0)*+{}="4";
(0,-5)*+{};
{ \ar@<.5ex>@/^.9pc/_-{i^!} "1";"2";};
{ \ar^-{i_* \simeq\, i_!} "2";"1";};
{ \ar@<-1ex>@/_.9pc/_-{i^*} "1";"2";};  
{ \ar@<.5ex>@/^.9pc/_-{j_*} "3";"4";};
{ \ar^-{j^! \simeq\, j^*} "4";"3";};
{ \ar@<-1ex>@/_.9pc/_-{j_!} "3";"4";};  
\end{xy}$$
and, in addition, the functor $j_*$ is $t$-exact.
\end{theorem}

\begin{proof}
We will prove the statement in detail in the case, where $X$ is affine and $\eta \colon Z \to X$ is a closed immersion defined by a regular sequence $(f_1,\dots,f_d)$ of homogeneous global sections. To reduce the general case to this special case, we will admit the following base-change result, the proof of which is similar to that of \cite[Corollary~3.4.2.2]{luriesag}: If $f \colon U \to X$ is a qcqs and flat map of Dirac schemes, and if $f' \colon U' \to X'$ is the base-change of $f$ along any map $g \colon X' \to X$ of Dirac schemes, then the square
$$\xymatrix{
{ \QCoh(X) } \ar[r]^-{f^*} \ar[d]^-{g^*} &
{ \QCoh(U) } \ar[d]^-{g'{}*} \cr
{ \QCoh(X') } \ar[r]^-{f'{}^*} &
{ \QCoh(U') } \cr
}$$
is right adjointable in the sense that the base-change map
$$\xymatrix{
{ g^*f_* } \ar[r] &
{ f_*'f'{}^*g^*f_* \simeq f_*'g'{}^*f^*f_* } \ar[r] &
{ f_*'g'{}^* } \cr
}$$
is an equivalence. The open immersion $i \colon U \to X$ is flat, and we claim that it is qcqs. Indeed, the property of being qcqs is local on the target, so we assume that $X$ is affine and that $\eta \colon Z \to X$ is defined by a regular sequence of global sections $(f_1,\dots,f_d)$. But then $i \colon U = X_{f_1} \cup \dots \cup X_{f_d} \to X$ is the open immersion of a finite union of distinguished open subsets, and hence, is quasi-compact, and it is also quasi-separated, because the diagonal $\Delta \colon U \to U \times_XU$ is an isomorphism.

We now suppose that $\mathfrak{B}$ is a basis for the topology of $|X|$, which is closed under finite intersections, and that the statement has been proved for all $V \in \mathfrak{B}$. The proof of \cref{lemma:scheme_a_colimit_of_its_basis_stable_under_intersection} shows that
$$\xymatrix{
{ \QCoh(X) } \ar[r] &
{ \varprojlim_{V \in \mathfrak{B}} \QCoh(V) } \cr
}$$
is an equivalence. A stable recollement is determined uniquely by a functor $i_* \simeq i_!$ with the property that it is fully faithful and admits both a left adjoint $i^*$ and a right adjoint $i^!$. The base-change property above implies, by \cite[Proposition 4.7.4.19]{lurieha}, that the $\mathfrak{B}$-indexed family of functors
$$\xymatrix@C=13mm{
{ \QCoh(U \cap V) } \ar[r]^-{(i|_{U \cap V})_*} &
{ \QCoh(V) } \cr
}$$
promotes to a map of $\mathfrak{B}^{\triangleleft}$-indexed limit diagrams in $\LPr$. In particular, the map of limits
$$\xymatrix{
{ \QCoh(U) } \ar[r]^-{i_*} &
{ \QCoh(X) } \cr
}$$
is a map in $\LPr$, and hence, admits a right adjoint $i^!$, as desired. It is also fully faithful, since this property is preserved under limits. Moreover, the fiber of $i^!$ is identified with $j_* \colon \QCoh(Y) \to \mathrm{QCoh}(X)$, since limits commute, and $j_*$ is $t$-exact, since it is left $t$-exact and since the property of being right $t$-exact can be tested after base-change to $V \in \mathfrak{B}$.

Now, for a general $X$, we take $\mathfrak{B}$ be the poset of open
subschemes $V \subset X$ that admit an open immersion in an affine
Dirac scheme. And if $X$ admits an open immersion $j \colon X \to S$
in an affine Dirac scheme, then we take $\mathfrak{B}$ to be the poset
of open subschemes $V \subset X$ with the property that $j(V) \subset
S$ is a distinguished open subscheme. This reduces us to the special
case, where $X \simeq \Spec(A)$ is affine, and where the graded ideal
$I \subset A$ that defines $\eta \colon Z \to X$ is generated by a
regular sequence  $(f_1,\dots,f_d)$ of homogeneous elements. In this
case, we find as in~\cite[Chapter~7]{luriesag} that there are
semi-orthogonal decompositions 
$$\xymatrix{
{ \Mod_A(\mathcal{D}(\Ab))^{\Nil(I)} } \ar@<.7ex>[r]^-{j_!} &
{ \Mod_A(\mathcal{D}(\Ab)) } \ar@<.7ex>[l]^-{j^!} \ar@<.7ex>[r]^-{i^*} &
{ \Mod_A(\mathcal{D}(\Ab))^{\Loc(I)} } \ar@<.7ex>[l]^-{i_*} \cr
}$$
\vspace{-3mm}
$$\xymatrix{
{ \Mod_A(\mathcal{D}(\Ab))^{\Cpl(I)} } \ar@<-.7ex>[r]_-{j_*} &
{ \Mod_A(\mathcal{D}(\Ab)) } \ar@<-.7ex>[l]_-{j^*} \ar@<-.7ex>[r]_-{i^!} &
{ \Mod_A(\mathcal{D}(\Ab))^{\Loc(I)} } \ar@<-.7ex>[l]_-{i_!} \cr
}$$
with $i_* \simeq i_!$ and that the composite adjoint functors
$$\xymatrix{
{ \Mod_A(\mathcal{D}(\Ab))^{\Nil(I)} } \ar@<.7ex>[r]^-{j_!} &
{ \Mod_A(\mathcal{D}(\Ab)) } \ar@<.7ex>[l]^-{j^!} \ar@<.7ex>[r]^-{j^*} &
{ \Mod_A(\mathcal{D}(\Ab))^{\Cpl(I)} } \ar@<.7ex>[l]^-{j_*} \cr
}$$
are equivalences. Moreover, one proves as in~\cite[Proposition~7.2.3.1]{luriesag} that
$$\xymatrix{
{ \QCoh(U) } \ar[r]^-{i_*} &
{ \QCoh(X) \simeq \Mod_A(\mathcal{D}(\Ab)) } \cr
}$$
is fully faithful and that its essential image is the full subcategory spanned by the $I$-local $A$-modules. So we must prove that the functor
$$\xymatrix{
{ \QCoh(Y) \simeq \varprojlim_m\Mod_{A/I^{m+1}}(\mathcal{D}(\Ab)) } \ar[r]^-{j_*} &
{ \QCoh(X) \simeq \Mod_A(\mathcal{D}(\Ab)) } \cr
}$$
is fully faithful and that its essential image is the full subcategory spanned by the $I$-complete $A$-modules. The functor $j_*$ factors as
$$\xymatrix{
{ \QCoh(Y) } \ar[r]^-{G} &
{ \Mod_A(\mathcal{D}(\Ab))^{\Cpl(I)} } \ar[r] &
{ \QCoh(X) \simeq \Mod_A(\mathcal{D}(\Ab)), } \cr
}$$
and since the full subcategory spanned by the $I$-complete $A$-modules is closed under small limits, the functor $G$ admits a left adjoint functor
$$\xymatrix@C+=10mm{
{ \QCoh(Y) \simeq \varprojlim_m\Mod_{A/I^{m+1}}(\mathcal{D}(\Ab)) } \ar@<-.7ex>[r]_-{G} &
{ \Mod_A(\mathcal{D}(\Ab))^{\Cpl(I)}. } \ar@<-.7ex>[l]_-{F} \cr
}$$
It follows from the definition of $I$-completeness that the unit map
$$\xymatrix{
{ M } \ar[r] &
{ (GF)(M) \simeq \varprojlim_m M \otimes_AA/I^{m+1} } \cr
}$$
is an equivalence. By the triangle identities, this, in turn, implies that the natural transformation $GFG \to G$ induced by the counit is an equivalence, so it remains to prove that $G$ is conservative. The functor $G$ takes a compatible family $(M_m)_{m \geq 1}$ with $M_m \in \Mod_{A/I^{m+1}}(\mathcal{D}(\Ab))$ to the limit $M \simeq \varprojlim_mM_m \in \Mod_A(\mathcal{D}(\Ab))^{\Cpl(I)}$, and we must show that if $M \simeq 0$, then $M_m \simeq 0$ for all $m \geq 1$.

To do so, we rewrite the limit in question. We write $P$ for the set of positive integers considered as a category with a single map $i \to j$ if $i \geq j$ and recall from~\cite[Proposition~5.3.1.20]{luriehtt} that the diagonal $\Delta \colon P \to P^d$ is a $\varprojlim$-equivalence. Given $e = (e_1,\dots,e_d) \in P^d$, we let $I_e = (f_1^{e_1},\dots,f_d^{e_d}) \subset A$ and observe that $I_{\Delta(m)} \subset I^m \subset I_{\Delta([m/d])}$. It follows that the canonical maps define equivalences
$$\begin{aligned}
\QCoh(Y) {} & \simeq \textstyle{ \varprojlim_{m \in P} \Mod_{A/I^m}(\mathcal{D}(\Ab)) } \cr
{} & \simeq \textstyle{ \varprojlim_{m \in P} \Mod_{A/I_{\Delta(m)}}(\mathcal{D}(\Ab)) } \cr
{} & \simeq \textstyle{ \varprojlim_{e \in P^d} \Mod_{A/I_e}(\mathcal{D}(\Ab)). } \cr
\end{aligned}$$
Hence, given a compatible family $(M_e)_{e \in P^d}$ with $M_e \in \Mod_{A/I_e}(\mathcal{D}(\Ab))$, we must prove that if $M \simeq \varprojlim_eM_e \simeq 0$, then $M_e \simeq 0$ for all $e \in P^d$, and to do so, we proceed by induction on $d \geq 1$. We first consider the case $d = 1$, where we abbreviate $f = f_1$ and $e = e_1$. Given any $N \in \Mod_A(\mathcal{D}(\Ab))$, we write $N/f^e$ for the cofiber of $f^e \cdot \id \colon N \to N$. Since $f \in A$ is a non-zero-divisor, the canonical map $A/f^e \to A/(f^e)$ is an equivalence. Hence, the cofiber $A/f^e$ has a ring structure. We observe that
$$\begin{aligned}
M_e/f {} & \simeq M_e \otimes_{A/(f^e)}((A/f^e)/f) \simeq M_e \otimes_{A/(f^e)}((A/f)/f^e) \cr
{} & \simeq M_e \otimes_{A/(f^e)}(A/f \oplus A/f[1]) \simeq M_1 \oplus M_1[1] \cr
\end{aligned}$$
and that, under this equivalence, the map $M_e \to M_{e-1}$ induces the identity map on the first summand and the zero map on the second summand. It follows that
$$\textstyle{ M/f \simeq \varprojlim_eM_e/f \simeq M_1. }$$
So if $M \simeq 0$, then $M_1 \simeq 0$, which implies that $M_e \simeq 0$ for all $e \in P$ by d\'{e}vissage. This proves the case $d = 1$. So we assume that the statement holds for $d = r$ and proceed to prove that it holds for $d = r+1$. We assume that $M \simeq \varprojlim_eM_e \simeq 0$ and must show that $M_e \simeq 0$ for all $e \in P^{r+1}$. We first write
$$\textstyle{ M \simeq \varprojlim_eM_e \simeq \varprojlim_{e_{r+1}}\varprojlim_{(e_1,\dots,e_r)} M_e \simeq \varprojlim_{e_{r+1}} N_{e_{r+1}}. }$$
By the inductive hypothesis, it will suffice to show that $N_{e_{r+1}} \simeq 0$ for all $e_{r+1} \in P$. But the $A$-module structure on $N_{e_{r+1}}$ extends to a structure of module over the completed Dirac ring
$$\textstyle{ \widehat{A} \simeq \varprojlim_{(e_1,\dots,e_r)}A/(f_1^{e_1},\dots,f_r^{e_r}), }$$
and the image of $f_{r+1}$ in $\smash{ \widehat{A} }$ is a non-zero-divisor. Indeed, it follows from~\cite[Theorem~16.1]{matsumura} that the image of $f_{r+1}$ in $A/(f_1^{e_1},\dots,f_r^{e_r})$ is a non-zero-divisor for all $(e_1,\dots,e_r) \in P^r$. Therefore, the same argument as in the case $d = 1$ shows that $N_{e_{r+1}} \simeq 0$ for all $e_{r+1} \in P$, as desired. This establishes the two semi-orthogonal decompositions.

Finally, we know that $j_*$ is left $t$-exact and must prove that it is also right $t$-exact, or equivalently, that it restricts to a functor
$$\xymatrix{
{ \varprojlim_m \Mod_{A/I^{m+1}}(\mathcal{D}(A))_{\geq 0} } \ar[r]^-{j_*} &
{ \Mod_A(\mathcal{D}(\Ab))_{\geq 0}. } \cr
}$$
But this follows from the Milnor sequence and from the fact that the maps
$$\xymatrix{
{ M_m } \ar[r] &
{ M_m \otimes_{A/I^{m+1}}A/I^m \simeq M_{m-1} } \cr
}$$
induce surjections on $\pi_0$.
\end{proof}

\subsection{Formal hyperplanes}

We define a map of Dirac stacks $q \colon Y \to X$ to be a formal hyperplane if, locally on $X$, it is isomorphic to the formal affine space associated with a locally free $\mathcal{O}_X$-module $\mathcal{E}$. In older literature, including \cite{goerss2008quasi,messing1972crystals}, formal hyperplanes are called formal Lie varieties.

Let $S$ be a Dirac scheme, and let $\mathcal{E}$ be an $\mathcal{O}_S$-module locally free of finite rank. The rank is locally constant, and we do not require that it be constant. We define the associated affine space over $S$ to be the Dirac $S$-scheme
$$\xymatrix@C=10mm{
{ \mathbb{A}_S(\mathcal{E}) = \Spec(\Sym_{\mathcal{O}_S}(\mathcal{E})) } \ar[r]^-{p} &
{ S, } \cr
}$$
and we define the zero section
$$\xymatrix@C=10mm{
{ S \simeq \mathbb{A}_S(0) } \ar[r]^-{\eta} &
{ \mathbb{A}_S(\mathcal{E}) }
}$$
to be the map induced by the unique map of $\mathcal{O}_S$-modules $\mathcal{E} \to 0$. Finally, we define the formal affine space over $S$ associated with $\mathcal{E}$ to be the formal completion
$$\xymatrix@C=10mm{
{ \affinespace_S(\mathcal{E}) = \Spf(\Sym_{\mathcal{O}_S}(\mathcal{E})) } \ar[r]^-{q} &
{ S } \cr
}$$
of the affine space $p \colon \mathbb{A}_S(\mathcal{E}) \to S$ along the zero section.

\begin{definition}
\label{def:formal_hyperplane_over_a_general_stack}
Let $X$ be a Dirac stack. A map of Dirac stacks $q \colon Y \to X$ is a formal hyperplane over $X$ if for every map $\eta \colon S \to X$ from an affine Dirac stack $S$, its base-change along $\eta$ is equivalent to a formal affine space.
\end{definition}

We write $\Hyp(X) \subset \sheaves(\Aff)_{/X}$ for the full subcategory spanned by the formal hyperplanes over $X$. It is a priori an $\infty$-category, but we now show that it is in fact equivalent to a $1$-category.

\begin{lemma}
\label{lemma:formal_hyperplanes_form_a_1-category}
If $q \colon Y \to X$ is a formal hyperplane, then for any map of Dirac stacks $\eta \colon Z \to X$, the mapping anima $\Map(Z,Y)$ calculated in $\sheaves(\Aff)_{/X}$ is $0$-truncated. In particular, the $\infty$-category $\Hyp(X)$ is equivalent to a $1$-category.
\end{lemma}

\begin{proof}
The inclusion $\spaces_{\leq 0} \to \spaces$ of $0$-truncated anima preserves limits, and it also preserves filtered colimits. In particular, by writing $Z$ as a colimit of affine Dirac schemes, we may assume that $Z \simeq S$ is affine. Since $q$ is a formal hyperplane, its base-change $q_S \colon Y_S \to S$ along $\eta \colon S \to X$ is equivalent to the formal affine space associated with some locally free $\mathcal{O}_S$-module of finite rank $\mathcal{E}$. Hence, we have
$$\textstyle{
\Map_X(S,Y) \simeq \Map_S(S,Y_S) \simeq \Map_S(S,\varinjlim_{\,m}S^{(m)}) \simeq \varinjlim_{\,m} \Map_S(S,S^{(m)}),
}$$
where the latter equivalence follows from \cref{remark:filtered_colimits_in_sheaves_levelwise_over_affines}. The right-hand side is a filtered colimit of $0$-truncated anima, and hence, is itself $0$-truncated.
\end{proof}

It follows immediately from the definition that if $f \colon X' \to X$ is any map of Dirac stacks, then the base-change along $f$ restricts to a functor
$$\xymatrix@C=10mm{
{ \Hyp(X) } \ar[r]^-{f^*} &
{ \Hyp(X'). } \cr
}$$
There functors are part of the functor
$$\xymatrix@C+=15mm{
{ \sheaves(\Aff)^{\op} } \ar[r]^-{\Hyp(-)} &
{ \verylargecat_{\infty} } \cr
}$$
which classifies the cartesian fibration
$$\xymatrix@C=10mm{
{ \Hyp } \ar[r]^-{p} &
{ \sheaves(\Aff) } \cr
}$$
given by the restriction of the target map $\Fun(\Delta^{1}, \sheaves(\Aff)) \rightarrow \sheaves(\Aff)$ to the full subcategory spanned by the formal hyperplanes; see~\cite[Notation~6.1.3.4]{luriehtt}.

\begin{warning}
\label{warning:formal_hyperplanes_dont_satisfy_descent}
The $\infty$-category $\Hyp(X)$ of formal hyperplanes over $X$ does generally not descend along effective epimorphisms, even if $X$ is a Dirac scheme. We refer to \cite[Warning~1.1.22]{lurieell2} for an instructive counterexample.
\end{warning}

We proceed to show that for pointed formal hyperplanes, the above descent issue goes away. As far as formal groups are concerned, it makes no difference if we work with formal hyperplanes or pointed formal hyperplanes. By definition, the $\infty$-category of pointed Dirac stacks over $X$ is the slice category 
$$(\sheaves(\Aff)_{/X})_{*} = (\sheaves(\Aff)_{/X})_{X \backslash} \simeq (\sheaves(\Aff)_{X\backslash})_{/X}.$$
Its objects are diagrams $\sigma \colon \Delta^{2} \to \sheaves(\Aff)$ of the form
$$\xymatrix@C=6mm{
{} &
{ Y } \ar[dr]^-{q} &
{} \cr
{ X } \ar[rr]^-{\id_X} \ar[ur]^-{\eta} &
{} &
{ X } \cr
}$$
and we view $q$ as a structure map and $\eta$ as a section thereof.

\begin{definition}
\label{definition:pointed_formal_hyperplane}Let $X$ be a Dirac stack. A pointed Dirac stack over $X$ is a pointed formal hyperplane if its base-change along every map $\eta \colon S \to X$ from an affine Dirac scheme is equivalent to
$$\begin{xy}
(12,14)*+{ \affinespace_S(\mathcal{E}) }="1";
(0,0)*+{ S }="2";
(24,0)*+{ S }="3";
{ \ar^-{\eta} "1";"2";};
{ \ar^-{q} "3";"1";};
{ \ar^-{\id_S} "3";"2";};
\end{xy}$$
for some locally free $\mathcal{O}_S$-module of finite rank $\mathcal{E}$. Here $\eta$ is the zero section.
\end{definition}

We will show that a pointed formal hyperplane admits a canonical exhaustive filtration. A map of Dirac stacks $f \colon Z \to X$ is $(-1)$-truncated if and only if the map of anima $f_S \colon Z(S) \to X(S)$ is $(-1)$-truncated for all affine Dirac schemes $S$. By~\cite[Proposition~6.2.3.17]{luriehtt}, the property of being $(-1)$-truncated is preserved by base-change along arbitrary maps and descends along effective epimorphisms. Also, if $f \colon Z \to X$ is $(-1)$-truncated and $\eta \colon Y \to X$ any map, then the anima of lifts
$$\begin{xy}
(12,14)*+{ Z }="1";
(0,0)*+{ Y }="2";
(24,0)*+{ X }="3";
{ \ar@{-->} "1";"2";};
{ \ar^-{f} "3";"1";};
{ \ar^-{\eta} "3";"2";};
\end{xy}$$
is $(-1)$-truncated so that for $\eta$ to factor through $f$ is a property. We now extend the definition of infinitesimal neighborhoods to Dirac stacks as follows.

\begin{definition}
\label{definition:infinitesimal_neighborhood}
Let $g \colon Z \to X$ be a $(-1)$-truncated map of Dirac stacks, and let $m \geq 0$ be an integer. The $m$th infinitesimal neighbourhood of $Z$ in $X$ is the essentially unique $(-1)$-truncated map of Dirac stacks $$\xymatrix@C=10mm{
{ Z_g^{(m)} } \ar[r]^-{g^{(m)}} &
{ X } \cr
}$$
such that a map $\eta \colon S \to X$ from an affine Dirac scheme factors through $g^{(m)}$ if and only if there exists an effective epimorphism $p \colon T' \to S$ from an affine Dirac scheme and a closed immersion $i \colon T \to T'$ defined by a quasi-coherent ideal, whose $(m+1)$th power is zero, such that $\eta|_T \colon T \to X$ factors through $g$.
\end{definition}

We will often abbreviate and write $g^{(m)} \colon Z^{(m)} \to X$ for the $m$th infinitesimal neighborhood of $Z$ in $X$. \cref{definition:infinitesimal_neighborhood} implicitly contains the assertion that $Z^{(m)}$ is a Dirac stack. We now verify that this is indeed the case. 

\begin{proposition}
\label{proposition:infinitesimal_neighbourhoods_are_dirac_stacks}
Let $g \colon Z \to X$ be a $(-1)$-truncated map of Dirac stacks, and let $g^{(m)} \colon Z^{(m)} \to X$ be its $m$th infinitesimal neighborhood with $m \geq 0$. In this situation, the presheaf $Z^{(m)} \colon \Aff^{\op} \to \spaces$ is a Dirac stack.
\end{proposition}

\begin{proof}
We must show that the presheaf $Z^{(m)}$ is accessible and that it is a sheaf for the flat topology. To this end, we define $g^{(m)}{}' \colon Z^{(m)}{}' \to X$ to be the essentially unique $(-1)$-truncated map of presheaves such that a map $\eta \colon S \to X$ from an affine Dirac scheme factors through $g^{(m)}{}'$ if and only if there exists a closed immersion $i \colon S' \to S$ defined by a quasi-coherent ideal, whose $(m+1)$th power is zero, such that $\eta |_{S'} \colon S' \to X$ factors through $g$. By definition, the inclusion $Z^{(m)}{}' \to Z^{(m)}$ is a sheafification for the flat topology, so in order to show that $Z^{(m)}$ is a Dirac stack, it will suffice to show that $Z^{(m)}{}'$ is accessible. 

There exists a regular cardinal $\kappa$ such that $X$ and $Z$ are $\kappa$-accessible in the sense that they agree with the left Kan extensions of their restrictions along the inclusion $i \colon \Aff_{\kappa} \to \Aff$ of the full subcategory of $\kappa$-small affine Dirac schemes. We claim that also $Z^{(m)}{}'$ is $\kappa$-accessible.

To prove the claim, we consider the diagram
$$\xymatrix{
{ \varinjlim Z^{(m)}{}'(T) } \ar[r] \ar[d] &
{ Z^{(m)}{}'(S) } \ar[d] \cr
{ \varinjlim X(T) } \ar[r] &
{ X(S) } \cr
}$$
with the colimits indexed by $p \colon S \to T$ in $(\Aff_{\kappa})_{S/}$. The bottom horizontal map is an equivalence, by the assumption that $X$ is $\kappa$-accessible, and we wish to show that so is the top horizontal map. The right-hand vertical map is $(-1)$-truncated by definition, and the left-hand vertical map is $(-1)$-truncated, since it is a filtered colimit of $(-1)$-truncated maps. Indeed, by \cref{lemma:filteredness_of_overcategories_of_accessible_cat}, the category $(\Aff_{\kappa})_{S/}$ is $\kappa$-filtered, and hence, filtered. We conclude that also the top horizontal map is $(-1)$-truncated, so it remains to show that it is essentially surjective.

Let us switch to ring notation. Given a point $x$ of $Z^{(m)}{}'(\Spec(A))$, we must show that there exits a map of Dirac rings $B \to A$ such that $B$ is $\kappa$-small and such that $x$ is in the essential image of $Z^{(m)}{}'(\Spec(B)) \to Z^{(m)}{}'(\Spec(A))$. By definition, we can identify $x$ with a point of $X(\Spec(A))$ with the property that there exists an ideal $I \subset A$ with $I^{m+1} = 0$ such that the image of $x$ by
$$\xymatrix{
{ X(\Spec(A)) } \ar[r] &
{ X(\Spec(A/I)) } \cr
}$$
belongs to $Z(\Spec(A/I))$. If we write $A$ as the $\kappa$-filtered colimit $A \simeq \varinjlim B$ of its $\kappa$-small sub-Dirac rings $B \subset A$, then also $A/I \simeq \varinjlim B/(B \cap I)$. Moreover, since $X$ and $Z$ are assumed to be $\kappa$-accessible, they preserve $\kappa$-filtered colimits, so we conclude that there exists a $\kappa$-small sub-Dirac ring $A_0 \subset A$ such that $x$ lifts to a point $x_0$ of $X(\Spec(A_0))$ with the property that its image by
$$\xymatrix{
{ X(\Spec(A_0)) } \ar[r] &
{ X(\Spec(A_0/(A_0 \cap I))) } \cr
}$$
belongs to $Z(\Spec(A_0/(A_0 \cap I)))$. But we have $(A_0 \cap I)^{m+1} = 0$ as an ideal of $A_0$, because $I^{m+1} = 0$ and because $A_0 \subset A$ is a sub-Dirac ring, so we conclude that the point $x_0$ belongs to $Z^{(m)}{}'(\Spec(A_0))$ as desired.
\end{proof}

\begin{remark}[Base-change for infinitesimal neighbourhoods]
\label{remark:infinitesimal_neighborhood_and_base_change}
Let $g \colon Z \to X$ be a $(-1)$-truncated map of Dirac stacks, let $m \geq 0$ be an integer, and let $f \colon X' \to X$ be any map of Dirac stacks. In this situation, the base-change $g' \colon Z' \to X'$ of $g$ along $f$ is $(-1)$-truncated, and $g'{}^{(m)} \colon Z'{}^{(m)} \to X'$ is the base-change of $g^{(m)} \colon Z^{(m)} \to X$ along $f$.
\end{remark}

We also verify that the notion of infinitesimal neighborhoods for Dirac stacks of \cref{definition:infinitesimal_neighborhood} extends the notion of infinitesimal neighborhoods of a closed embedding Dirac schemes introduced in \cref{definition:formal_completion_along_closed_immersion}. 

\begin{lemma}
\label{lemma:infinitesimal_neighborhood_of_closed_immersion_of_schemes}
Let $g \colon Z \to X$ be a closed immersion of Dirac schemes. The induced map $h(g) \colon h(Z) \to h(X)$ of Dirac stacks is $(-1)$-truncated and the canonical map
$$\xymatrix{
{ h(Z_g^{(m)}) } \ar[r] &
{ h(Z)_{h(g)}^{(m)} } \cr
}$$
is an equivalence for every $m \geq 0$.
\end{lemma}

\begin{proof}It is clear that the map in the statement exists and is $(-1)$-truncated, so only the essential surjectivity of this map is at issue. So we let $\eta \colon S \to X$ be a map from an affine Dirac scheme such that $h(\eta)$ factors through $h(g)^{(m)}$ and  wish to show that $\eta$ factors $g^{(m)}$.

Since $h(\eta)$ factors through $h(g)^{(m)}$, there exists, by definition, a faithfully flat map $p \colon T' \to S$ from an affine Dirac scheme and a closed immersion $i \colon T \to T'$ defined by a quasi-coherent ideal of $\mathcal{O}_{T'}$, whose $(m+1)$th power is zero, such that $\eta|_T \colon T \to X$ factors through $g$. We now consider the diagram
$$\xymatrix{
{ T } \ar[r]^{i} &
{ T' } \ar[r]^-{p} &
{ S } \ar[r]^-{\eta} &
{ X } \cr
{ Z_T^{(m)} } \ar[r] \ar[u] &
{ Z_{T'}^{(m)} } \ar[r] \ar[u] &
{ Z_S^{(m)} } \ar[r] \ar[u] &
{ Z^{(m)} } \ar[u]_-{g^{(m)}} \cr
{ Z_T } \ar[r] \ar[u] &
{ Z_{T'} } \ar[r] \ar[u] &
{ Z_S } \ar[r] \ar[u] &
{ Z } \ar[u]_-{h_m} \cr
}$$
with all squares cartesian. Since $g^{(m)}$ is a monomorphism, the statement that $\eta$ factors through $g^{(m)}$, which we wish to prove, is equivalent to the statement that the base-change of $g^{(m)}$ along $\eta$ is an isomorphism. And since $p$ is faithfully flat, the common statement, in turn, is equivalent to the statement that the base-change of $g^{(m)}$ along $\eta\circ p$ is an isomorphism. By the same reasoning, the assumption that $\eta|_T$ factors through $g \simeq g^{(m)}h_m$ implies that the base-changes of both $g^{(m)}$ and $h_m$ along $\eta \circ p \circ i$ are isomorphisms.

We now consider the two left-hand columns in diagram above together with the horizontal maps between them. The maps in this part of the diagram all induce homeomorphisms of underlying spaces, and moreover, the maps in the left-hand column are isomorphisms. Hence, if $I,J \subset \mathcal{O}_{T'}$ are the quasi-coherent ideals that define $g_{T'}$ and $i$, respectively, then $I \subset J$. But then $I^{m+1} \subset J^{m+1} = 0$, which shows that the base-change of $g^{(m)}$ along $\eta \circ p$ is an isomorphism, as desired.
\end{proof}

\begin{definition}
\label{definition:formal_completion_along_monomorphism}
Let $g \colon Z \to X$ be a $(-1)$-truncated map of Dirac stacks. The formal completion of $X$ along $g$ is the induced $(-1)$-truncated map
$$\xymatrix{
{ Y_g \simeq \varinjlim_{\,m}Z_g^{(m)} } \ar[r]^-{j} &
{ X } \cr
}$$
from the colimit of the infinitesimal neighborhoods of $Z$ in $X$. The Dirac stack $X$ is formally complete along $g$ if the map $j$ is an equivalence.
\end{definition}

We next show that the formal completion along a $(-1)$-truncated map of Dirac stacks is an idempotent operation in the following precise sense.

\begin{proposition}
\label{proposition:formal_completion_is_formally_complete}
Let $g \colon Z \to X$ be a $(-1)$-truncated map of Dirac stacks, let $j \colon Y \to X$ be the formal completion of $X$ along $g$, and let $h \colon Z \to Y$ be the unique map such that $g \simeq j \circ h$. In this situation, the map $j$ induces an equivalence
$$\xymatrix{
{ Z_h^{(m)} } \ar[r] &
{ Z_g^{(m)} } \cr
}$$
for all $m \geq 0$. In particular, the Dirac stack $Y$ is formally complete along $h$.
\end{proposition}

\begin{proof}It is clear that the map in the statement exists and is $(-1)$-truncated. So we let $\eta \colon S \to X$ be a map from an affine Dirac scheme and assume that $\eta$ factors through $g^{(m)}$. In particular, it factors through $j$, so that $\eta \simeq j \circ \eta'$ for a unique map $\eta' \colon S \to Y$, and we wish to show that $\eta'$ factors through $h^{(m)}$. Now, since $\eta$ factors through $g$, there exists, by definition, a faithfully flat map $p \colon T' \to S$ with $T'$ affine and a closed immersion $i \colon T \to T'$ defined by a quasi-coherent ideal, whose $(m+1)$th power is zero, such that $\eta|_T$ factors through $g \simeq j \circ h$. Since $j$ is $(-1)$-truncated, this implies that $\eta'|_T$ factors through $h$, so $\eta'$ factors, by definition, through $h^{(m)}$ as we wanted to prove.
\end{proof}

We now specialize to the case of the $(-1)$-truncated map of Dirac stacks given by the zero section of a pointed formal hyperplane.

\begin{corollary}\label{corollary:properties_of_formal_neighbourhoods_of_pointed_hyperplane}
Let $q \colon Y \to X$ be a formal hyperplane pointed by $\eta \colon X \to Y$.
\begin{enumerate}
\item[{\rm (1)}]The Dirac stack $Y$ is formally complete along $\eta$.
\item[{\rm (2)}]The composite map $q_m \simeq q \circ \eta^{(m)} \colon X_{\eta}^{(m)} \to X$ is finite affine for all $m \geq 0$, and it is is flat for $0 \leq m \leq 1$.
\end{enumerate}
\end{corollary}

\begin{proof}The statement~(1) is a special case of \cref{proposition:formal_completion_is_formally_complete}. The statement~(2) can be verified after base-change along any map from an affine Dirac scheme, so it suffices to consider the case, where $q \colon Y \simeq \affinespace_S(\mathcal{E}) \to S$ is the formal affine space associated with a locally free $\mathcal{O}_S$-module $\mathcal{E}$, and where $\eta \colon S \to Y$ is the zero section. In this case, \cref{proposition:formal_completion_is_formally_complete} identifies
$$\xymatrix@C=10mm{
{ S^{(m)} \simeq \Spec(\Sym_{\mathcal{O}_S}(\mathcal{E}) / \Sym_{\mathcal{O}_S}^{> m}(\mathcal{E})) } \ar[r]^-{\eta^{(m)}} &
{ \affinespace_S(\mathcal{E}) \simeq \Spf(\Sym_{\mathcal{O}_S}(\mathcal{E})), } \cr
}$$
from which~(2) follows.
\end{proof}

\begin{warning}We note that, in the presence of sections of half-integer spin, the symmetric powers $\Sym^d_{\mathcal{O}_S}(\mathcal{E})$ of a locally free $\mathcal{O}_S$-module of finite rank $\mathcal{E}$ are not necessarily flat for $d \geq 2$. For instance,
$$\Sym_{\mathbb{Z}}^d(\mathbb{Z}(\nicefrac{1}{2})) \simeq \begin{cases}
\mathbb{Z}(\nicefrac{d}{2}) & \text{for $0 \leq d \leq 1$,} \cr
\mathbb{Z}/2(\nicefrac{d}{2}) & \text{for $d \geq 2$,} \cr
\end{cases}$$
is not a flat $\mathbb{Z}$-module for $d \geq 2$.
\end{warning}

\begin{remark}
\label{remark:recovering_locally_free_module_from_pointed_hyperplane_in_affine_case}
If $q \colon Y \to S$ is the formal affine space associated with a locally free $\mathcal{O}_S$-module $\mathcal{E}$ of finite rank pointed by the zero section $\eta \colon S \to Y$, then, by the proof of \cref{corollary:properties_of_formal_neighbourhoods_of_pointed_hyperplane}, there is a cofiber sequence of quasi-coherent $\mathcal{O}_S$-modules
$$\xymatrix{
{ \mathcal{O}_S } \ar[r] &
{ q_{1*}\mathcal{O}_{S_{\eta}^{(1)}} } \ar[r] &
{ \mathcal{E}. } \cr
}$$
In particular, the cofiber of the left-hand map is a locally free $\mathcal{O}_S$-module.
\end{remark}

The following definition globalizes this observation.

\begin{definition}
\label{definition:conormal_sheaf_of_pointed_hyperplane}
Let $q \colon Y \to X$ be a formal hyperplane pointed by $\eta \colon X \to Y$.
\begin{enumerate}
\item[(1)]The conormal sheaf of $q \colon Y \to X$ at $\eta \colon X \to Y$ is the cofiber
$$\mathcal{N}_{X/Y} \simeq \cof(\mathcal{O}_X \to q_{1*}\mathcal{O}_{X_{\eta}^{(1)}}) \in \QCoh(X).$$
\item[(2)]The tangent space of $q \colon Y \to X$ at $\eta \colon X \to Y$ is the vector bundle
$$T_{Y/X,\eta} \simeq V_X(\mathcal{N}_{X/Y}) \in \Vect(X).$$
\end{enumerate}
\end{definition}

The assignment of the conormal sheaf defines a functor
$$\xymatrix{
{ \Hyp_*(X)^{\op} } \ar[r] &
{ \QCoh(X), } \cr
}$$
or equivalently, the assignment of the tangent space defines a functor
$$\xymatrix{
{ \Hyp_*(X) } \ar[r] &
{ \Vect(X). } \cr
}$$

\begin{proposition}
\label{proposition:conormal_sheaf_of_pointed_hyperplane_satisfies_base_change}
Let $f \colon X' \to X$ be a map of Dirac stacks, let $q \colon Y \to X$ be a formal hyperplane pointed by $\eta \colon X \to Y$, and let $q' \colon Y' \to X'$ be the base-change of $q$ along $f$ pointed by $\eta' \simeq (\eta,\id)$. In this situation, the canonical map
$$\xymatrix{
{ f^*\mathcal{N}_{X/Y} } \ar[r] &
{ \mathcal{N}_{X'/Y'} } \cr
}$$
is an equivalence.
\end{proposition}

\begin{proof}
Since $f^* \colon \QCoh(X) \rightarrow \QCoh(X')$ is exact and $f^* \mathcal{O}_{X} \simeq \mathcal{O}_{X'}$, it will suffice to show that in the cartesian diagram
$$\xymatrix@C=10mm{
{ X_{\eta'}^{(1)} } \ar[r]^-{f'} \ar[d]^-(.45){q_1'} &
{ X_{\eta}^{(1)} } \ar[d]^-{q_1} \cr
{ X' } \ar[r]^{f} &
{ X, } \cr
}$$
the base-change map $f^*q_{1*} \to q_{1*}'f'{}^*$ is an equivalence. But this follows from \cref{prop:pushforward_along_affine_flat_satisfies_base_change}, since $q_1$ is affine flat by \cref{corollary:properties_of_formal_neighbourhoods_of_pointed_hyperplane}.
\end{proof}

\begin{corollary}
If $q \colon Y \to X$ is a formal hyperplane pointed at $\eta \colon X \to Y$, then the conormal sheaf $\mathcal{N}_{X/Y}$ is a locally free $\mathcal{O}_X$-module of finite rank.
\end{corollary}

\begin{proof}
By \cref{proposition:conormal_sheaf_of_pointed_hyperplane_satisfies_base_change}, we can assume that $X$ is affine, in which case the statement follows from \cref{remark:recovering_locally_free_module_from_pointed_hyperplane_in_affine_case}.
\end{proof}

If $q \colon Y \to X$ is a formal hyperplane, then so is its base-change $q' \colon Y' \to Y$ along itself, and moreover, the formal hyperplane $q' \colon Y' \to Y$ is canonically pointed by the diagonal $\Delta \colon Y \to Y'$.

\begin{definition}
\label{definition:tangent_bundle}
Let $q \colon Y \to X$ be a formal hyperplane.
\begin{enumerate}
\item[(1)]The $\mathcal{O}_Y$-module of differentials of $q \colon Y \to X$ is the conormal sheaf
$$\Omega_{Y/X}^1 \simeq \mathcal{N}_{Y/Y'} \in \QCoh(Y)$$
of the base-change $q' \colon Y' \to Y$ of $q$ along itself pointed by the diagonal.
\item[(2)]The tangent bundle of $q \colon Y \to X$ is the vector bundle
$$T_{Y/X} \simeq V_Y(\Omega_{Y/X}^1) \in \Vect(Y).$$
\end{enumerate}
\end{definition}

\begin{remark}
\label{remark:tangent_bundle}
Let $q \colon Y \to X$ be a formal hyperplane, and let $q' \colon Y' \to Y$ be the base-change of $q$ along itself, pointed by the diagonal $\Delta \colon Y \to Y'$. If $\eta \colon X \to Y$ is a point of $q \colon Y \to X$, then the base-change of $q' \colon Y' \to Y$ pointed by $\Delta \colon Y \to Y'$ along $\eta \colon X \to Y$ is canonically identified with $q \colon Y \to X$ pointed by $\eta \colon X \to Y$,
$$\xymatrix@C=12mm{
{ Y } \ar[r] \ar[d]_-{q} &
{ Y' } \ar[r] \ar[d]_-(.46){q'} &
{ Y\phantom{,} } \ar[d]^-{q} \cr
{ X } \ar[r]^-{\eta} \ar@/_.7pc/[u]_-{\eta} &
{ Y } \ar[r]^-{q} \ar@/_.7pc/[u]_-{\Delta} &
{ X, } \cr
}$$
so  \cref{proposition:conormal_sheaf_of_pointed_hyperplane_satisfies_base_change} shows that $\mathcal{N}_{X/Y} \simeq \eta^*\Omega_{Y/X}^1$ and $T_{Y/X,\eta} \simeq \eta^*T_{Y/X}$.
\end{remark}

\begin{remark}
\label{remark:conormal_sheaf_of_pointed_hyperplane}
Our definition of the tangent bundle is provisional. In~\cite{diracgeometry3}, we will define the cotangent complex $L_{Y/X}$ of a formal hyperplane $q \colon Y \to X$ and prove that it is a locally free $\mathcal{O}_Y$-module isomorphic to $\Omega_{Y/X}^1$.
\end{remark}

If $q \colon Y \to X$ is a formal hyperplane pointed by $\eta \colon X \to Y$, then the rank of the conormal sheaf $\mathcal{N}_{X/Y}$ is locally constant in the sense that for every map $f \colon S \to X$ from a Dirac scheme, the rank the locally free $\mathcal{O}_S$-module $f^*\mathcal{N}_{X/Y}$ is a locally constant function on the underlying topological space $|S|$.

\begin{definition}
\label{definition:dimension_of_pointed_formal_hyperplane}
The dimension of a pointed formal hyperplane is the rank of its conormal sheaf.
\end{definition}

\begin{warning}
It is tempting also define the spin of a pointed formal hyperplane as the spin of local generators of its conormal sheaf. But said spin is not well-defined, so this is not possible. For example, if $k$ is the Dirac field $\mathbb{F}_2[x^{\pm1}]$ with $x$ of spin $\nicefrac{1}{2}$, then as a $k$-vector space, $k$ is generated by $x^d$ for any integer $d$.
\end{warning}

\begin{remark}
If $X$ is an affine Dirac scheme, then \cref{remark:recovering_locally_free_module_from_pointed_hyperplane_in_affine_case} shows that every pointed formal hyperplane over $X$ is isomorphic to the formal affine space of its conormal sheaf pointed by the zero section. However, such an isomorphism is highly non-canonical, and if $X$ is not affine, then it may not exist.
\end{remark}

\subsection{Automorphisms of pointed formal hyperplanes and descent} 
\label{subsection:automorphisms_of_pointed_formal_hyperplanes_and_descent}

Given a formal hyperplane $q \colon Y \to X$, a point $\eta \colon X \to Y$ determines by \cref{corollary:properties_of_formal_neighbourhoods_of_pointed_hyperplane} an exhaustive increasing filtration of $Y$ by the infinitesimal neighborhoods of $X$ in $Y$. This filtration, in turn, induces a filtration of the automorphism group, which we proceed to describe. We first promote the automorphism group to a group stack.

We show in \cref{theorem:existence_of_additional_right_adjoint} that the $\infty$-categories of accessible sheaves of anima have the same functoriality as do $\infty$-topoi. So, given a map $f \colon Y \to X$ in $\sheaves(\Aff)$, we have adjoint functors
\vspace{-1mm}
$$\begin{xy}
(0,0)*+{ \sheaves(\Aff)_{/X} };
(33,0)*+{ \sheaves(\Aff)_{/Y} };
(8.5,0)*+{}="1";
(24,0)*+{}="2";
(0,-5)*+{};
{ \ar@<.5ex>@/^.8pc/_-{f_*} "1";"2";};
{ \ar^-{f^*} "2";"1";};
{ \ar@<-1ex>@/_.8pc/_-{f_!} "1";"2";};  
\end{xy}$$
with $f_!$ given by composition with $f$ and $f^*$ given by base-change along $f$. Moreover, the base-change formula holds in the sense that given a cartesian square
$$\xymatrix@C=9mm{
{ Y' } \ar[r]^-{j'} \ar[d]^-{f'} &
{ Y } \ar[d]^{f} \cr
{ X' } \ar[r]^-{j} &
{ X } \cr
}$$
in $\sheaves(\Aff)$, the base-change maps $f_!'\,j'{}^* \to j^*f_!$ and $j^*f_* \to f_*'\,j'{}^*$ are equivalences. This has the following consequence:

\begin{proposition}
\label{proposition:self-enriched}
The slice $\infty$-category $\sheaves(\Aff)_{/X}$ is cartesian closed. 
\end{proposition}

\begin{proof}
Suppose that $f \colon Y \to X$ and $g \colon Z \to X$ are two objects of the slice category, and let $\mathrm{id}_{X} \colon X \to X$ be the terminal object. We claim that internal mapping object $\Hom(f, g)$ is given by $f_*f^*g_!g^*\mathrm{id}_{X} $. To see this, note that for any $h \colon R \rightarrow X$,
$$\Map(h, f_*f^*g_!g^*\mathrm{id}_{X} ) \simeq \Map(f_!f^*h, g_{!}g^{*}\mathrm{id}_{X} ) \simeq \Map(f \times h, g),$$
where the mapping spaces and the product $f \times h$ are calculated in the slice category $\sheaves(\Aff)_{/X}$. 
\end{proof}

\begin{remark}
\label{remark:classifying_stack_for_automorphisms_of_object_in_slice_category}
Given an object $f \colon Y \to X$ in $\sheaves(\Aff)_{/X}$, the internal mapping object $\Hom(f,f)$ canonically promotes to an $\mathbf{E}_1$-monoid $\End(f)$ in $\sheaves(\Aff)_{/X}$. We denote its maximal sub-$\mathbf{E}_1$-group by $\Aut(f)$ and refer to it as the classifying stack for automorphisms of $f$.
\end{remark}

\begin{remark}
\label{remark:classifying_stack_for_automorphisms_of_pointed_object_in_slice_category}
Suppose that $\sigma \in (\sheaves(\Aff)_{/X})_*$ is a pointed object in $\sheaves(\Aff)_{/X}$. In a similar manner, we have a classifying stack
$$\Aut(\sigma) \in \Grp_{\mathbf{E}_1}((\sheaves(\Aff)_{/X})_*) \simeq \Grp_{\mathbf{E}_1}(\sheaves(\Aff)_{/X})$$
for automorphisms of $\sigma$ as a pointed object. Writing $\sigma$ as a diagram
$$\xymatrix@C=6mm{
{} &
{ Y } \ar[dr]^-{f} &
{} \cr
{ X } \ar[ur]^-{\eta} \ar[rr]^{\id_X} &&
{ X } \cr
}$$
in $\sheaves(\Aff)$, we also use the notation $\Aut_{\eta}(Y/X)$ instead of $\Aut(\sigma)$.
\end{remark}

Let $q \colon Y \to X$ be a formal hyperplane pointed by $\eta \colon X \to Y$. The infinitesimal neighborhoods of $X$ in $Y$ define a filtered object in $(\sheaves(\Aff)_{/X})_*$ with
$$\xymatrix@C=5mm{
{} &
{ X_{\eta}^{(m)} } \ar[dr]^-{q_m} &
{} \cr
{ X } \ar[ur]^-{\eta_m} \ar[rr]^-{\id_X} &&
{ X } \cr
}$$
as its $m$th term. 

\begin{definition}
\label{definition:canonical_filtration_of_automorphism_group_of_pointed_formal_hyperplane}
Let $q \colon Y \to X$ be a formal hyperplane pointed by $\eta \colon X \to Y$, and let $G \simeq \Aut_{\eta}(Y/X)$ be its classifying stack for automorphisms. The canonical filtration of $G$ is the filtered object in $\Grp_{\mathbf{E}_1}(\sheaves(\Aff)_{/X})$ with
$$G^{\hspace{.5pt}m} \simeq \Aut_{\eta_m}(X_{\eta}^{(m)}/X).$$
\end{definition}

\begin{remark}
We recall from \cref{corollary:properties_of_formal_neighbourhoods_of_pointed_hyperplane} that, in the situation of \cref{definition:canonical_filtration_of_automorphism_group_of_pointed_formal_hyperplane}, the Dirac stack $Y$ is formally complete along $\eta \colon X \to Y$. Thus, the canonical filtration of $G$ is complete in the sense that
$$\textstyle{ G \simeq \varprojlim_m\,G^{\hspace{.5pt}m}. }$$
Since $q_m \colon X_{\eta}^{(m)} \to X$ is affine, it is in particular $0$-truncated. This implies that also the $\mathbf{E}_1$-group $G^{\hspace{.5pt}m}$ is $0$-truncated, or equivalently, that it represents a functor taking values in ordinary groups.
\end{remark}

\begin{remark}
\label{remark:classifying_stack_for_automorphisms_of_quasicoherent_module}
Given $\mathcal{F} \in \QCoh(X)$, we similarly have a classifying stack
$$\Aut(\mathcal{F}) \in \Grp_{\mathbf{E}_1}(\sheaves(\Aff)_{/X})$$ that represents the functor $\sheaves(\Aff)_{/X} \to \Grp_{\mathbf{E}_1}(\spaces)$ that to $j \colon X' \to X$ assigns the $\mathbf{E}_1$-group of automorphisms of $j^*(\mathcal{F}) \in \QCoh(X')$. Indeed, this functor descends along effective epimorphisms, because $\QCoh(X)$ does so, and one may show that it is accessible.
\end{remark}

\begin{remark}
\label{remark:classifying_stack_for_automorphisms_of_vector_bundle}
If $\mathcal{E} \in \QCoh(X)^{\loc}$ is a locally free $\mathcal{O}_X$-module with corresponding vector bundle $V \simeq V_X(\mathcal{E}) \in \Vect(X)$, then the $\mathbf{E}_1$-group
$$\Aut(V) \simeq \Aut(\mathcal{E})^{\op} \in \Grp_{\mathbf{E}_1}(\sheaves(\Aff)_{/X})$$
represents the functor $\sheaves(\Aff)_{/X} \to \Grp_{\mathbf{E}_1}(\spaces)$ that to the map $j \colon X' \to X$ assigns the $\mathbf{E}_1$-group of automorphisms of $j^*(V) \in \Vect(X')$.
\end{remark}

Let $q \colon Y \to X$ be a formal hyperplane pointed by $\eta \colon X \to Y$. Taking the tangent space at $\eta$ defines a map of classifying stacks for automorphisms
$$\xymatrix{
{ G \simeq \Aut_{\eta}(Y/X) } \ar[r] &
{ \Aut(\,T_{Y/X,\eta}) } \cr
}$$
that to an automorphism $f$ assigns its derivative $df_{\eta}$ at $\eta$.

\begin{lemma}
\label{lemma:automorphism_of_inf_1_same_as_those_of_conormal_sheaf}
The derivative at $\eta$ factors through an equivalence
$$\xymatrix{
{ G^1 \simeq \Aut_{\eta_1}(X_{\eta}^{(1)}/X) } \ar[r] &
{ \Aut(\,T_{Y/X,\eta}). } \cr
}$$
\end{lemma}

\begin{proof}We must show that the induced map of functors represented by the map in the statement is an equivalence pointwise for all $j \colon X' \to X$. It will suffice to consider the case, where $X' \simeq \Spec(A)$ is affine. Indeed, we show in \cref{appendix:accessible_sheaves} that the Yoneda embedding extends to an equivalence $\sheaves(\Aff_{/X}) \to \sheaves(\Aff)_{/X}$. In this case, the statement follows from the well-known fact that the map
$$\xymatrix{
{ \CAlg(\Ab)_{A/-/A} } \ar[r] &
{ \Mod_A(\Ab) } \cr
}$$
that to an $A$-algebra $q_1 \colon A \to B$ augmented by $\eta_1 \colon B \to A$ assigns the cokernel of $q_1$ is an equivalence of categories.
\end{proof}

\begin{lemma}
\label{lemma:higher_filtration_quotiones}
For all $m \geq 2$, the following hold:
\begin{enumerate}
\item[{\rm (1)}]The canonical map $G^{\hspace{.5pt}m} \to G^{\hspace{.5pt}m-1}$ is an effective epimorphism.
\item[{\rm (2)}]If $j \colon X' \to X$ has $X' \simeq \Spec(A)$ affine, then
$$\ker(G^{\hspace{.5pt}m} \to G^{\hspace{.5pt}m-1})(X' \to X) \simeq \Hom_A(P,\Sym_A^m(P))$$
with $P \simeq \mathcal{N}_{Y/X}(X' \to X)$. In particular, the kernel is an abelian group.
\end{enumerate}
\end{lemma}

\begin{proof}(1)~We claim that the map is a $\pi_0$-surjection pointwise for every $j \colon X' \to X$ with $X'$ affine. Granting this, the lemma follows, because the Yoneda embedding extends to an equivalence $\sheaves(\Aff_{/X}) \to \sheaves(\Aff)_{/X}$ by \cref{appendix:accessible_sheaves}. To prove the claim, we write $X' \simeq \Spec(A)$ and $Y' \simeq \Spf(\Sym_A(P))$ for some finitely generated projective $A$-module $P$. We let $q \colon A \to \Sym_A(P)$ and $\eta \colon \Sym_A(P) \to A$ be the structure map and the augmentation, respectively, and let $I \subset \Sym_A(P)$ be the augmentation ideal. We wish to prove that the canonical map
$$\xymatrix{
{ \Aut(\Sym_A(P)/I^{m+1}) } \ar[r] &
{ \Aut(\Sym_A(P)/I^m) } \cr
}$$
is surjective. By the universal property of the symmetric algebra, an endomorphism of the augmented $A$-algebra $\Sym_A(P)/I^{m+1}$ determines and is determined by an $A$-linear map $P \to I/I^{m+1}$. Since $P$ is projective, any $A$-linear map $P \to I/I^m$ factors through the projection $I/I^{m+1} \to I/I^m$, so
$$\xymatrix{
{ \End(\Sym_A(P)/I^{m+1}) } \ar[r] &
{ \End(\Sym_A(P)/I^m) } \cr
}$$
is surjective. But an endomorphism of the augmented $A$-algebra $\Sym_A(P)/I^{m+1}$ is an automorphism if and only if the induced map of the associated graded for the $I$-adic filtration is an isomorphism if and only if the induced $A$-linear map
$$\xymatrix{
{ P } \ar[r] &
{ I/I^{m+1} } \ar[r] &
{ I/I^2 } \cr
}$$
is an isomorphism. This proves the claim.

\noindent(2)~The proof of~(1) shows that
$$\ker(G^{\hspace{.5pt}m} \to G^{\hspace{.5pt}m-1})(X' \to X) \simeq \Hom_A(P,I^m/I^{m+1}),$$
and $I^m/I^{m+1} \simeq \Sym_A^m(P)$, as we wanted to show.
\end{proof}

We return to our purpose in this section of showing that the $\infty$-category of pointed formal hyperplanes descends along effective epimorphisms.

\begin{lemma}
\label{lemma:formal_pointed_hyperplanes_descend_for_fpqc_maps_of_affines}
If $f \colon S' \to S$ is a faithfully flat map of affine Dirac stacks, then
$$\xymatrix{
{ \Hyp_*(S) } \ar[r]^-{f^*} \ar[d] &
{ \Hyp_*(S') } \ar[d] \cr
{ (\sheaves(\Aff)_{/S})_* } \ar[r]^-{f^*} &
{ (\sheaves(\Aff)_{/S'})_* } \cr
}$$
is a cartesian square of $\infty$-categories.
\end{lemma}

\begin{proof}
By definition, pointed formal hyperplanes are preserved under base-change, so the square exists. Thus, to complete the proof, we must show that if $\sigma$ is a pointed Dirac stack over $S$, whose base-change $\sigma'$ along $f \colon S' \to S$ is a pointed formal hyperplane, then $\sigma$ is itself a pointed formal hyperplane. We write $q \colon Y \to S$ and $\eta \colon S \to Y$ for the structure map and point of $\sigma$, and similarly for $\sigma'$.

We claim that $Y$ is formally complete along $S$. Indeed, the canonical map
$$\xymatrix{
{ \varinjlim_{\,m}S^{(m)} } \ar[r] &
{ Y } \cr
}$$
is an equivalence, because its base-change
$$\xymatrix{
{ \varinjlim_{\,m}S'{}^{(m)} } \ar[r] &
{ Y' } \cr
}$$
along the effective epimorphism $f$ is so. We also claim that $q_1 \colon S^{(1)} \to S$ is affine flat. Indeed, its base-change $q_1'$ along $f$ is affine flat by \cref{corollary:properties_of_formal_neighbourhoods_of_pointed_hyperplane}, and the property of being affine flat descends along effective epimorphisms. Now, the cofiber
$$\mathcal{E} \simeq \cof(
\mathcal{O}_S \to q_{1*}\mathcal{O}_{S^{(1)}} ) \in \QCoh(S)$$
is a candidate for the conormal sheaf $\mathcal{N}_{S/Y}$ of $\sigma$ as in  \cref{definition:conormal_sheaf_of_pointed_hyperplane}, except that we do not yet know that $\sigma$ is a pointed formal hyperplane. Since $$f^*\mathcal{E} \simeq \mathcal{N}_{S'/Y'}$$
is a locally free $\mathcal{O}_{S'}$-module of finite rank, and since the property of being locally free of finite rank descends along effective epimorphisms, we conclude that $\mathcal{E}$ is a locally free $\mathcal{O}_S$-module of finite rank.

To complete the proof, we will show that there exists an equivalence
$$\xymatrix{
{ Z \simeq \affinespace_S(\mathcal{E}) } \ar[r]^-{h} &
{ Y } \cr
}$$
of pointed Dirac stacks over $S$ from the formal affine space associated with $\mathcal{E}$ pointed by the zero section $\zeta \colon S \to Z$. Since both $Y$ and $Z$ are formally complete along $S$, it will suffice to show that there exists a compatible sequence of equivalences
$$\xymatrix@C=10mm{
{ S_{\zeta}^{(m)} } \ar[r]^-{h^{(m)}} &
{ S_{\eta}^{(m)} } \cr
}$$
for all $m \geq 0$. The case $m = 0$ is trivial, and the case $m = 1$ follows from the definition of $\mathcal{E}$. So we let $m \geq 2$ and show that an equivalence $h^{(m-1)}$ admits an extension to an equivalence $h^{(m)}$. To this end, we let $\mathcal{F}^{\,m}$ be the sheaf on $\Aff_{/S}^{\hspace{.5pt}\flat}$ that to $g \colon U \to S$ affine flat assigns the anima of diagrams
$$\xymatrix@C=14mm{
{ U \times_SS_{\zeta}^{(m-1)} } \ar[r]^-{h^{(m-1)}|_U} \ar[d] &
{ U \times_SS_{\eta}^{(m-1)} } \ar[d] \cr
{ U \times_SS_{\zeta}^{(m)} } \ar[r]^-{h_U^{(m)}} &
{ U \times_SS_{\eta}^{(m)} } \cr
}$$
with $h_U^{(m)}$ an equivalence. The sheaf $\mathcal{F}^{\,m}$ is naturally a torsor for the $\mathbf{E}_1$-group
$$\mathcal{K}^{\hspace{.5pt}m} \simeq \ker(G^{\hspace{.5pt}m} \to G^{\hspace{.5pt}m-1})$$
of automorphisms of $S_{\eta}^{(m)}$ that extend the identity on $S_{\eta}^{(m-1)}$. The isomorphism classes of such torsors are parametrized by the sheaf cohomology group $H^1(S,\mathcal{K}^{\hspace{.5pt}m})$. But \cref{lemma:higher_filtration_quotiones} identifies this group with a coherent cohomology group of the affine Dirac scheme $S$, and therefore, it is zero. Hence, the torsor $\mathcal{F}^{\,m}$ is trivial, so it admits a global section, which is the desired equivalence $h^{(m)}$.
\end{proof}

\begin{theorem}
\label{theorem:pointed_hyperplanes_takes_colimits_of_stacks_to_catinfty_limits}
The functor $\Hyp_{*} \colon \sheaves(\Aff)^{\op} \to \verylargecat_{\infty}$ preserves small limits.
\end{theorem}

\begin{proof}
Suppose that $X \simeq \varinjlim_{\alpha}X_{\alpha}$ is a colimit of a small diagram in $\sheaves(\Aff)$. We consider the diagram
$$\xymatrix{
{ \Hyp_*(X) } \ar[r] \ar[d] &
{ \varprojlim_{\alpha} \Hyp_*(X_{\alpha}) } \ar[d] \cr
{ (\sheaves(\Aff)_{/X})_* } \ar[r] &
{ \varprojlim_{\alpha} (\sheaves(\Aff)_{/X_{\alpha}})_* } \cr
}$$
and wish to prove that the top horizontal map is an equivalence. We claim that the bottom horizontal map is an equivalence. By~\cref{theorem:slice_infty_categories_of_stacks_satisfy_descent}, the map
$$\xymatrix{
{ \sheaves(\Aff)_{/X} } \ar[r] &
{ \varprojlim_{\alpha} \sheaves(\Aff)_{/X_{\alpha}} } \cr
}$$
is an equivalence, so if suffices to show that the canonical map
$$\xymatrix{
{ (\varprojlim_{\alpha} \sheaves(\Aff)_{/X_{\alpha}})_* } \ar[r] &
{ \varprojlim_{\alpha} (\sheaves(\Aff)_{/X_{\alpha}})_* } \cr
}$$
is an equivalence, which follows from the adjunctions
$$\xymatrix{
{ \verylargecat_{\infty} } \ar@<.7ex>[r] &
{ \verylargecat\hspace{-.5pt}{}_{\infty}^1 } \ar@<.7ex>[l] \ar@<-.7ex>[r] &
{ \verylargecat\hspace{-.5pt}{}_{\infty}^{\operatorname{ptd}} } \ar@<-.7ex>[l] \cr
}$$
Here, the bottom left-hand map is the inclusion of the (non-full) subcategory spanned by the $\infty$-categories that admit a final object and the maps that preserve final objects, and the top right-hand map is the inclusion of the full subcategory spanned by the pointed $\infty$-categories. The former is a right adjoint, by (the dual of) \cite[Corollary~5.3.6.10]{luriehtt}, and therefore, it preserves limits. Moreover, the latter admits a right adjoint that to an $\infty$-category $\mathcal{C}$ that admits a final object assigns the $\infty$-category $\mathcal{C}_* \simeq \mathcal{C}_{1/}$ of pointed objects in $\mathcal{C}$. This proves the claim.

Now, in the diagram at the beginning of the proof, the vertical maps are fully faithful, so the claim shows that the top horizontal map is fully faithful. It remains to prove that it is essentially surjective. So we fix $\sigma \in (\sheaves(\Aff)_{/X})_*$ such that $i_{\alpha}^*(\sigma) \in \Hyp_*(X_{\alpha})$ for all $\alpha$ and proceed to show that $\sigma \in \Hyp_*(X)$.

First, we consider the special case, where $T \simeq \coprod_{\alpha}T_{\alpha}$ is a finite coproduct of affine Dirac schemes. In this case, the desired statement follows from that the fact that for a locally free $\mathcal{O}_T$-module of finite rank $\mathcal{E}$, the $\mathcal{O}_{T_{\alpha}}$-modules $\mathcal{E}_{\alpha} \simeq i_{\alpha}^*(\mathcal{E})$ are locally free of finite rank, and the diagram
$$\xymatrix{
{ \coprod_{\alpha} \affinespace_{T_{\alpha}}(\mathcal{E}_{\alpha}) } \ar[r] \ar[d] &
{ \affinespace_T(\mathcal{E}) } \ar[d] \cr
{ \coprod_{\alpha} T_{\alpha} } \ar[r] &
{ T } \cr
}$$
is cartesian with the horizontal maps equivalences.

Next, in the general case, we must show that for $\eta \colon S \to X$ with $S$ affine, $\eta^*(\sigma) \in \Hyp_*(S)$. By \cref{lemma:formal_pointed_hyperplanes_descend_for_fpqc_maps_of_affines}, if suffices to show that there exists a faithfully flat map of affine Dirac schemes $g \colon T \to S$ such that $g^*\eta^*(\sigma) \in \Hyp_*(T)$. To this end, we consider the canonical map of presheaves
$$\xymatrix{
{ \varinjlim_{\alpha} \iota(X_{\alpha}) } \ar[r] &
{ \iota(\varinjlim_{\alpha} X_{\alpha}). } \cr
}$$
It becomes an equivalence, and therefore, an effective epimorphism after sheafification, so it follows from \cref{proposition:characterization_of_effective_epimorphisms_of_sheaves} that there exists a diagram of presheaves
$$\xymatrix@C=13mm{
{ \coprod_{i \in I}h(T_i) } \ar[r]^-{\sum_i\tilde{\eta}_i} \ar[d]^-{\sum_ih(g_i)} &
{ \varinjlim_{\alpha} \iota(X_{\alpha}) } \ar[d] \cr
{ h(S) } \ar[r]^-{\eta} &
{ \iota(\varinjlim_{\alpha} X_{\alpha}) } \cr
}$$
with $(g_i \colon T_i \to S)_{i \in I}$ a finite covering family. Moreover, since colimits of presheaves are calculated pointwise, each $\tilde{\eta}_i$ factors through some $\iota(X_{\alpha_i})$, so
$$g_i^*\eta^*(\sigma) \simeq \tilde{\eta}_i^*i_{\alpha_i}^*(\sigma) \in \Hyp_*(T_i)$$
for all $i \in I$. Now, the coproduct $g \colon T \to S$ of the $g_i \colon T_i \to S$ is a faithfully flat map of affine Dirac schemes, and the special case that we considered above shows that
$g^*\eta^*(\sigma) \in \Hyp_*(T)$ as desired. This completes the proof.
\end{proof}

\begin{corollary}
\label{cor:formal_hyperplanes_descent_along_effective_epis}
If $f \colon X' \to X$ is an effective epimorphism of Dirac stacks, then
$$\xymatrix{
{ \Hyp_*(X) } \ar[r] &
{ \varprojlim_{\Delta} \Hyp_*(X'{}^{\times_X[-]}) } \cr
}$$
is an equivalence of $\infty$-categories. 
\end{corollary}

\begin{proof}
This is a special case of \cref{theorem:pointed_hyperplanes_takes_colimits_of_stacks_to_catinfty_limits}.
\end{proof}

\begin{corollary}
\label{corollary:pointed_hyperplanes_form_essentially_small_category}
For every Dirac stack $X$, $\Hyp_*(X)$ is essentially small.
\end{corollary}

\begin{proof}
If $X \simeq \Spec(R)$ is affine, then the cardinality of the set of equivalence classes of formal hyperplanes over $X$ is equal, by \cref{remark:recovering_locally_free_module_from_pointed_hyperplane_in_affine_case}, to that of isomorphism classes of finitely generated projective graded $R$-modules, which is small. In the general, we can write $X \simeq \varinjlim X_{\alpha}$ as a small colimit of affine Dirac schemes, and since $\Hyp_{*}(X) \simeq \varprojlim \Hyp_{*}(X_{\alpha})$ by \cref{theorem:pointed_hyperplanes_takes_colimits_of_stacks_to_catinfty_limits}, the corollary follows.
\end{proof}

\subsection{Formal groups}

In this section, we define the groupoid of formal groups over a Dirac stack and show that it descends along effective epimorphisms.

\begin{definition}
\label{definition:formal_groups}
Let $\Lat$ be the category of finitely generated free abelian groups. If $\mathcal{C}$ is an $\infty$-category, which admits finite products, then the $\infty$-category of abelian group objects in $\mathcal{C}$ is the full subcategory
$$\Ab(\mathcal{C}) \simeq \Fun^{\Sigma}(\Lat^{\op},\mathcal{C}) \subset \Fun(\Lat^{\op},\mathcal{C})$$
spanned by the functors that preserve finite products.
\end{definition}

\begin{remark}
We note that the full subcategory $\Ab(\mathcal{C}) \subset \Fun(\Lat^{\op}, \mathcal{C})$ is closed under sifted colimits and limits by~\cite[Proposition~5.5.8.10]{luriehtt}. Therefore, it follows from \cite[Corollary~5.1.2.3]{luriehtt} that small limits and small sifted colimits in $\Ab(\mathcal{C})$ are calculated pointwise. 
\end{remark}

By  \cref{corollary:pointed_hyperplanes_form_essentially_small_category}, the functor $\Hyp_*$ takes values in the $\infty$-category $\Cat_{\infty}$ of small $\infty$-categories. In fact, it takes values in the subcategory $\Cat_{\infty}^{\Pi} \subset \Cat_{\infty}$ spanned by the $\infty$-categories that admit finite products and the finite product preserving functors, and hence, we may compose it with the functor
$$\xymatrix@C=11mm{
{ \Cat_{\infty}^{\Pi} } \ar[r]^-{\phantom{,}\Ab\phantom{,}} &
{ \Cat_{\infty} } \cr
}$$
that to $\mathcal{C}$ assigns $\Ab(\mathcal{C}) \simeq \Fun^{\Sigma}(\Lat^{\op},\mathcal{C})$. 

\begin{definition}[Formal groups]
\label{def:formalgroup}
The functor that to a Dirac stack $X$ assigns the $\infty$-category $\FGroup(X)$ of formal groups over $X$ is the composite functor
$$\xymatrix@C=11mm{
{ \sheaves(\Aff)^{\op} } \ar[r]^-{\Hyp_*} &
{ \Cat_{\infty}^{\Pi} } \ar[r]^-{\phantom{,}\Ab\phantom{,}} &
{ \Cat_{\infty}. } \cr
}$$
\end{definition}

\begin{example}[Formal additive groups]
\label{ex:additiveformalgroup}
Let $X$ be a Dirac stack, and let $\mathcal{E}$ be a locally free $\mathcal{O}_X$-module of finite rank. The functor $\smash{ \additivegroup(\mathcal{E}) \colon \Lat^{\op} \to \Hyp(X) }$
that to $L$ assigns $q_L \colon \affinespace_X(L \otimes_{\mathbb{Z}} \mathcal{E}) \to X$ is a formal group over $X$, which we call the formal additive group associated with $\mathcal{E}$. Its dimension is equal to the rank of $\mathcal{E}$.
\end{example}

\begin{proposition}
\label{proposition:formal_groups_takes_colimits_of_stacks_to_catinfty_limits}
The functor $\FGroup \colon \sheaves(\Aff)^{\op} \to \Cat_{\infty}$ takes small colimits of Dirac stacks to limits of $\infty$-categories.
\end{proposition}

\begin{proof}
Indeed, this is a consequence of \cref{theorem:pointed_hyperplanes_takes_colimits_of_stacks_to_catinfty_limits} and from the fact, which is proved in (the dual of)~\cite[Corollary~5.3.6.10]{luriehtt}, that $\Cat_{\infty}^{\Pi} \subset \Cat_{\infty}$ is closed under limits.
\end{proof}

\begin{remark}[Formal groups and formal group laws]
\label{rem:formalgroup}
We spell out the definition of a formal group in more concrete terms. So let $\mathcal{G} \colon \Lat^{\op} \to \Hyp(X)$ be a formal group over $X$, and let $q \colon Y \to X$ be its underlying formal hyperplane $\mathcal{G}(\mathbb{Z})$. The canonical inclusions $\mathbb{Z} \to \mathbb{Z} \oplus \mathbb{Z}$ induce an equivalence
$$\xymatrix{
{ \mathcal{G}(\mathbb{Z} \oplus \mathbb{Z}) } \ar[r] &
{ \mathcal{G}(\mathbb{Z}) \times \mathcal{G}(\mathbb{Z}), } \cr
}$$
and a composition of an inverse of this equivalence and the map induced by the diagonal map $\mathbb{Z} \to \mathbb{Z} \oplus \mathbb{Z}$ defines the addition
$$\xymatrix{
{ Y \times_XY } \ar[r]^-{+} &
{ Y } \cr
}$$
in the formal group. Similarly, the unique map $\mathbb{Z} \to 0$ defines the zero section
$$\xymatrix{
{ X } \ar[r]^-{\eta} &
{ Y. } \cr
}$$
Conversely, a formal hyperplane $q \colon Y \to X$ together with  maps $+ \colon Y \times_XY \to Y$ and $\eta \colon X \to Y$ satisfying the abelian group axioms determine a formal group $\mathcal{G} \colon \Lat^{\op} \to \Hyp(X)$, which is unique up to contractible choice.

Suppose that $X \simeq \Spec(R)$ and $Y \simeq \Spf(A)$ with $I \subset A$ the ideal that defines the zero section $\eta \colon X \to Y$. The map $+ \colon Y \times_XY \to Y$ is a point of the anima
$$\textstyle{ \Map_X(Y \times_XY,Y) \simeq \varprojlim_{m,n}\varinjlim_p \Map_R(A/I^{p+1},A/I^{m+1} \otimes_RA/I^{n+1}) }$$
which is $0$-truncated. Moreover, Zariski locally on $X$, we may choose homogeneous coordinates $x_1,\dots,x_d \in A$ such that the family $(x_1,\dots,x_d)$ generates $I \subset A$ and such that the family $(x_1+I^2,\dots,x_d+I^2)$ is a basis of the $R$-module $I/I^2$. In terms of these coordinates, the map $+ \colon Y \times_XY \to Y$ determines and is determined by a family $f(y,z) = (f_s(y,z))$ of $d$ power series in $2d$ variables,
$$f_s = f_s(y,z) = f_s(y_1,\dots,y_d,z_1,\dots,z_d) \in R[[y_1,\dots,y_d,z_1,\dots,z_d]],$$
such that
$$\deg(f_s) = \deg(x_s) = \deg(y_s) = \deg(z_s).$$
Moreover, the map $+ \colon Y \times_XY \to Y$ defines a structure of formal group on the formal hyperplane $q \colon Y \to X$ with zero section $\eta \colon X \to Y$ defined by $I \subset A$ if and only if the following identities among families of power series hold.
\begin{enumerate}
 \item[(a)]$f(x,0) = x = f(0,x)$
 \item[(b)]$f(f(x,y),z) = f(x,f(y,z))$
 \item[(c)]$f(x,y) = f(y,x)$
\end{enumerate}
A family of power series satisfying these conditions is classically known as a \emph{formal group law}; the above discussion implies that any formal group is locally described (in a non-canonical way) by a formal group law. 

We remark that the above three axioms imply the existence and uniqueness of an inverse with respect to addition. Indeed, a map $- \colon Y \to Y$ determines and is determined by a family $i(x) = (i_s(x))$ of $d$ power series in $d$ variables,
$$i_s = i_s(x) = i_s(x_1,\dots,x_d) \in R[[x_1,\dots,x_d]],$$
such that $\deg(i_s) = \deg(x_s)$, and the map $- \colon Y \to Y$ is an inverse with respect to addition if and only if the following identity holds.
\begin{enumerate}
\item[(d)]$f(x,i(x)) = 0 = f(i(x),x)$
\end{enumerate}
But it follows from~(a) that this system of equations admits a unique solution.
\end{remark}

\begin{example}
\label{ex:additiveformalgroupagain}
Let $R$ be a Dirac ring, and let $M$ be a free $R$-module. If we choose a basis $(x_1,\dots,x_d)$ of $M$ that consists of homogeneous elements, then the additive formal group $\additivegroup(M)$ associated with $M$ determines and is determined by the family of power series $f(y,z) = y + z$ with $f_s(y,z) = y_s+z_s$.
\end{example}

\begin{definition}
\label{def:liealgebra}
Let $X$ be a Dirac stack. The Lie algebra is the functor
$$\xymatrix@C=10mm{
{ \FGroup(X) } \ar[r]^-{\phantom{,}\Lie\phantom{,}} &
{ \Vect(X) } \cr
}$$
that to a formal group $\mathcal{G}$ assigns the tangent space $\Lie(\mathcal{G}) \simeq T_{Y/X,\eta}$ of its underlying formal hyperplane $q \colon Y \to X$ at the zero section $\eta \colon X \to Y$.
\end{definition}

Since $\mathcal{G}$ is abelian, it serves no purpose to define a Lie bracket on the vector bundle $\Lie(\mathcal{G})$, because it will turn out to be zero.

The functor $\FGroup \colon \sheaves(\Aff)^{\op} \to \Cat_{\infty}$ classifies a cartesian fibration
$$\xymatrix{
{ \FGroup } \ar[r]^-{p} &
{ \sheaves(\Aff), } \cr
}$$
and we refer to its domain, which by abuse of notation we also denote by $\FGroup$, as the $\infty$-category of formal groups. Its objects are formal groups $\mathcal{G}$ over varying base Dirac stacks. A map $f \colon \mathcal{G}' \to \mathcal{G}$ in $\FGroup$ determines a diagram
$$\xymatrix@C=10mm{
{ Y' } \ar[r] \ar[d]^-{q'} &
{ Y\phantom{,} } \ar[d]^-{q} \cr
{ X' } \ar[r] &
{ X } \cr
}$$
in $\sheaves(\Aff)$, where $q$ and $q'$ are the underlying formal hyperplanes of $\mathcal{G}$ and $\mathcal{G}'$, respectively, and the map $f$ is $p$-cartesian if and only if this diagram of Dirac stacks is cartesian. We recall that a functor $\mathcal{G} \colon K \to \FGroup$ is said to be cartesian if it takes every map $\alpha \colon k' \to k$ in $K$ to a cartesian map in $\FGroup$.

\begin{proposition}
\label{prop:formalgroups}
Let $\mathcal{G} \colon K \to \FGroup$ be a cartesian diagram in the $\infty$-category of formal groups. There exists an essentially unique pair of a cartesian diagram 
$$\xymatrix@C=10mm{
{ K^{\triangleright} } \ar[r]^-{\bar{\mathcal{G}}} &
{ \FGroup } \cr
}$$
and an equivalence $\bar{\mathcal{G}}|_K \simeq \mathcal{G}$. Moreover, the diagram $\bar{\mathcal{G}}(\mathbb{Z}) \colon K^{\triangleright} \to \sheaves(\Aff)^{\Delta^1}$ of underlying formal hyperplanes is a colimit diagram.
\end{proposition}

\begin{proof}
We consider the following diagram of $\infty$-categories
$$\xymatrix@C=10mm{
{ \FGroup_K } \ar[r]^-{X'} \ar[d]^-(.45){p_K} &
{ \FGroup } \ar[d]^-{p} \cr
{ K } \ar[r]^-(.6){X} \ar[ur]^-{\mathcal{G}} \ar@/^1pc/[u]^-{\mathcal{G}_K} &
{ \sheaves(\Aff), } \cr
}$$
where $p_K$ is the base-change of $p$ along $X \simeq p \circ \mathcal{G}$. Since $p$ is a cartesian fibration, so is $p_K$, and since $\mathcal{G}$ is a cartesian functor, so is the induced section $\mathcal{G}_K$ of $p_K$. We fix a pair of a colimit diagram $\bar{X} \colon K^{\triangleright} \to \sheaves(\Aff)$ and an equivalence $\bar{X}|_K \simeq X$ and consider the diagram
$$\xymatrix@C=10mm{
{ \FGroup_K } \ar[r] \ar[d]^-(.45){p_K} &
{ \FGroup_{K^{\triangleright}} } \ar[r]^-{\bar{X}'} \ar[d]^-{p_{K^{\triangleright}}} &
{ \FGroup } \ar[d]^-{p} \cr
{ K } \ar[r] \ar@/^1pc/[u]^-{\mathcal{G}_K} &
{ K^{\triangleright} } \ar[r]^-{\bar{X}} &
{ \sheaves(\Aff), } \cr
}$$
where $p_{K^{\triangleright}}$ is the base-change of $p$ along $\bar{X}$. We proved in \cref{proposition:formal_groups_takes_colimits_of_stacks_to_catinfty_limits} that the functor $\FGroup \colon \sheaves(\Aff)^{\op} \to \Cat_{\infty}$ takes colimits of Dirac stacks to limits of $\infty$-categories. Therefore, it follows from~\cite[Proposition~3.3.3.1]{luriehtt} that cartesian section $\mathcal{G}_K$ of $p_K$ extends essentially uniquely to a cartesian section $\mathcal{G}_{K^{\triangleright}}$ of $p_{K^{\triangleright}}$. This determines a cartesian functor $\smash{ \bar{\mathcal{G}} \colon K^{\triangleright} \to \FGroup }$ and an equivalence $\smash{ \bar{\mathcal{G}}|_K \simeq \mathcal{G} }$, as desired.

It remains to prove that $\bar{\mathcal{G}}(\mathbb{Z}) \colon K^{\triangleright} \to \sheaves(\Aff)^{\Delta^1}$ is a colimit diagram. Let us write $q_k \colon Y_k \to X_k$ for the value of $\bar{\mathcal{G}}(\mathbb{Z})$ at $k \in K^{\triangleright}$, and let $\infty \in K^{\triangleright}$ be the cone point. Since $\bar{\mathcal{G}} \colon K^{\triangleright} \to \FGroup$ is a cartesian functor, the diagram
$$\xymatrix{
{ Y_k } \ar[r] \ar[d]^-{q_k} &
{ Y_{\infty} } \ar[d]^-{q_{\infty}} \cr
{ X_k } \ar[r] &
{ X_{\infty} } \cr
}$$
is cartesian, for all $k \in K$, and since colimits in $\sheaves(\Aff)$ are universal, so is
$$\xymatrix{
{ \varinjlim_{k \in K}Y_k } \ar[r] \ar[d] &
{ Y_{\infty}\phantom{.} } \ar[d]^-{q_{\infty}} \cr
{ \varinjlim_{k \in K}X_k } \ar[r] &
{ X_{\infty}. } \cr
}$$
But the lower horizontal map is an equivalence, and hence, so is the top horizontal map.
\end{proof}

\section{Stable homotopy theory}

The purpose of this final section is to demonstrate that Dirac geometry is the natural geometric framework in which to formulate various statements in stable homotopy theory.

\subsection{The descent spectral sequence}

Suppose that $\eta \colon k \to E$ be a map of commutative algebras in spectra.\footnote{\,We assume that $\eta \colon k \to E$ is a map of commutative algebras in spectra for ease of exposition, but less will do: If $k$ is an $\mathbf{E}_2$-algebra, then $\Mod_k(\Sp)$ promotes to an $\mathbf{E}_1$-monoidal $\infty$-category, and we can ask that $E$ be an $\mathbf{E}_1$-algebra therein. We refer to \cite[Construction~2.7]{mathew2017nilpotence} for the definition of the cobar construction in this generality. Since the graded ring $\smash{ \pi_*(E) }$ is not automatically Dirac, we must assume this to be the case, too.} We may form the augmented cosimplicial object 
in the $\infty$-category of $k$-modules in spectra
$$\xymatrix@C=12mm{
{ \Delta_+ } \ar[r]^-{E^{\otimes_k[-]}} &
{ \Mod_k(\Sp), } \cr
}$$
which is known as the cobar construction of $\eta$. We recall that $[n] = \{0,1,\dots,n\}$ so that $E^{\otimes_k[n]}$ is an $(n+1)$-fold tensor power of $E$ over $k$. Given a $k$-module in spectra $V$, the restriction
$$\xymatrix{
{ V \simeq \varprojlim_{\Delta_+} V \otimes_k E^{\otimes_k[-]}} \ar[r] &
{ \widehat{V} \simeq \varprojlim_{\Delta} V \otimes_k E^{\otimes_k[-]} } \cr
}$$
along the inclusion $\Delta \subset \Delta_+$ is known as the $\eta$-nilpotent completion of $V$. It is typically not a particularly well-behaved operation, but it has the advantage that its target can be analyzed by means of the descent spectral sequence, which we now explain.

We let $\Delta^{\leq n} \subset \Delta$ be the full subcategory spanned by the objects $[m]$ with $m \leq n$ and recall that Lurie's $\infty$-categorical Dold--Kan correspondence
$$\xymatrix{
{ \Fun(\Delta,\Mod_k(\Sp)) } \ar[r] &
{ \Fun(\mathbb{Z}_{\geq 0}^{\op},\Mod_k(\Sp)) } \cr
}$$
from \cite[Theorem~1.2.4.1]{lurieha}
takes the cosimplicial $k$-module in spectra $V \otimes_kE^{\otimes_k[-]}$ to the filtered $k$-module in spectra $\Fil_{\eta}(V)$ that to $n$ assigns
$$\textstyle{ \Fil_{\eta}^n(V) \simeq \varprojlim_{\Delta^{\leq n}} V \otimes_kE^{\otimes_k[-]}. }$$
The descent spectral sequence is the associated spectral sequence
$$\xymatrix@C=5mm{
{ E_{i,j}^1 \simeq \pi_{i+j}(\gr_{\eta}^{-i}(V)) } \ar@{=>}[r] &
{ \pi_{i+j}(\widehat{V}), } \cr
}$$
which converges conditionally in the sense of Boardman \cite{boardman} to the homotopy groups of the $\eta$-nilpotent completion of $V$. We recall the following result concerning the convergence of the descent spectral sequence.

\begin{lemma}
Let $\eta \colon k \to E$ be a map of commutative algebra in spectra, let $I$ be the $k$-module in spectra given by the fiber of $\eta$, and let $V$ be any $k$-module in spectra. In this situation, there is a fiber sequence of $k$-modules in spectra
$$\xymatrix{
{ V \otimes_k I^{\otimes_k[n]} } \ar[r] &
{ V } \ar[r] &
{ \Fil_\eta^n(V) } \cr
}$$
for all $n \geq 0$. In particular, if $\eta$ is $1$-connective and if $V$ is right bounded, then the completion map $V \to \widehat{V}$ is an equivalence.
\end{lemma}

\begin{proof}
This is proved in~\cite[Proposition~2.14]{mathew2017nilpotence}.
\end{proof}

The proof of the $\infty$-categorical Dold--Kan correspondence identifies
$$\xymatrix@C=10mm{
{ \cdots } &
{ E_{i-1,j}^1 } \ar[l]_-(.4){d^1} &
{ E_{i,j}^1 } \ar[l]_-(.4){d^1} &
{ \cdots } \ar[l]_-(.4){d^1} &
{ E_{-1,j}^1 } \ar[l]_-(.4){d^1} &
{ E_{0,j}^1 } \ar[l]_-(.4){d^1} \cr
}$$
with the normalized cochain complex associated with the cosimplicial $\pi_0(k)$-module $\smash{ \pi_j(V \otimes_kE^{\otimes_k[-]}) }$. It follows that the descent spectral sequence takes the form
$$\xymatrix@C=5mm{
{ E_{i,j}^2 \simeq \pi^{-i}(\pi_j(V \otimes_kE^{\otimes_k[-]})) } \ar@{=>}[r] &
{ \pi_{i+j}(\widehat{V}), } \cr
}$$
and we proceed to identify the $E^2$-term in Dirac geometric terms.

The map $\eta \colon k \to E$ of commutative algebra in spectra gives rise to the groupoid
$$S \simeq \Spec(\pi_*(E^{\otimes_k[-]}))$$
in affine Dirac schemes, so we may form the Dirac stack given by the colimit
$$\textstyle{ X \simeq \varinjlim_{\Delta^{\op}} S \in \sheaves(\Aff) }$$
in Dirac stacks. We henceforth assume that the groupoid $S$ has flat face maps so that, by \cref{theorem:geometricity_in_terms_of_groupoids}, the Dirac stack $X$ is geometric.\footnote{\,Equivalently, we assume that $\eta$ be descent-flat in the terminology of \cite{burklund2022adams}, where it is proved that this property of $\eta$ is not stable under composition.}

\begin{construction}
\label{construction:quasi_coherent_sheaf_of_coefficients}
In this situation, we define a pair $(\mathcal{F},\partial)$ of a functor
$$\xymatrix@C=10mm{
{ \Mod_k(\Sp) } \ar[r]^-{\mathcal{F}} &
{ \QCoh(X)^{\heartsuit} } \cr
}$$
that to a $k$-module in spectra $V$ assigns a quasi-coherent $\mathcal{O}_X$-module $\mathcal{F}(V)$ and a natural isomorphism of quasi-coherent $\mathcal{O}_X$-modules
$$\xymatrix@C=10mm{
{ \mathcal{F}(\Sigma V) } \ar[r]^-{\partial_V} &
{ \mathcal{F}(V)(-\nicefrac{1}{2}), }
}$$
which is a homology theory in the sense of~\cite[Definition~2.8]{patchkoria2021adams}. To this end, we let
$$\xymatrix@C=10mm{
{ \QCoh^{\heartsuit} } \ar[r]^-{p} &
{ \Aff } \cr
}$$
be the cartesian fibration classified by the functor $\Aff^{\op} \to \verylargecat_{\infty}$ that to an affine Dirac scheme $T$ assigns the abelian category $\QCoh(T)^{\heartsuit}$. Now, we have a lift
$$\xymatrix@C=15mm{
{} & 
{ \QCoh^{\heartsuit} } \ar[d]^-{p} \cr
{ \Delta^{\op} } \ar@<.5ex>[ur]^-{\widetilde{S}(V)\;} \ar[r]^-{S} &
{ \Aff } \cr
}$$
defined by
$$\widetilde{S}(V) \simeq (S,\pi_*(V \otimes_kE^{\otimes_k[-]}))$$
and it takes injective maps in $\Delta$ to $p$-cartesian maps in $\QCoh^{\heartsuit}$ by the assumption that $S$ takes injective maps in $\Delta$ to flat maps in $\Aff$. So let $\Delta_s \subset \Delta$ be the subcategory spanned by the injective maps. It follows from \cref{cor:effectiveepimorphismqcoh} and by \cite[Corollary~3.3.3.2]{luriehtt} that the canonical maps
$$\xymatrix{
{ \QCoh(X)_{\geq 0} } \ar[r] &
{ \varprojlim_{\Delta} \QCoh(S)_{\geq 0} } \ar[r] &
{ \varprojlim_{\Delta_s} \QCoh(S)_{\geq 0} } \cr
}$$
are equivalences and from \cref{theorem:qcoh_on_geometric_dirac_stack_is_grothendieck_prestable} that the composite map restricts to an equivalence
$$\xymatrix{
{ \QCoh(X)^{\heartsuit} } \ar[r] &
{ \varprojlim_{\Delta_s} \QCoh(S)^{\heartsuit}. } \cr
}$$
By~\cite[Corollary~3.3.3.2]{luriehtt}, the cartesian lift $\widetilde{S}(V)$ defines a point on the right-hand side, and we define $\mathcal{F}(V) \in \QCoh(X)^{\heartsuit}$ to be the essentially unique corresponding point on the left-hand side.

It is clear that the cartesian lift $\widetilde{S}(V)$, and hence, the point $\mathcal{F}(V)$ is functorial in $V$, so we obtain the desired functor $\mathcal{F}$. Moreover, the natural isomorphism
$$\xymatrix{
{ \pi_*(\Sigma V \otimes_kE^{\otimes_k[-]}) } \ar[r] &
{ \pi_*(V \otimes_kE^{\otimes_k[-]})(-\nicefrac{1}{2}) } \cr
}$$
from~\cite[Example~2.4]{diracgeometry1} induces the natural isomorphism $\partial_V \colon \mathcal{F}(\Sigma V) \to \mathcal{F}(V)(-\nicefrac{1}{2})$.
\end{construction}

\begin{remark}
\label{remark:quasi_coherent_sheaf_of_coefficients}
Let $\eta \colon k \to E$ be a map of commutative algebras in spectra with the property that the groupoid $S \simeq \Spec(\pi_*(E^{\otimes_k[-]}))$ has flat face maps, and let $X \simeq |S|$. Since $f_n \colon S_n \to X$ is a submersion, \cref{theorem:qcoh_on_geometric_dirac_stack_is_grothendieck_prestable} implies that $f_n^*$ is $t$-exact, so we may consider the composition
$$\xymatrix{
{ \Mod_k(\Sp) } \ar[r]^-{\mathcal{F}} &
{ \QCoh(X)^{\heartsuit} } \ar[r]^-{f_n^*} &
{ \QCoh(S_n)^{\heartsuit} \simeq \Mod_{\pi_*(E^{\otimes_k[n]})}(\Ab). } \cr
}$$
We note that the definition of $\mathcal{F}$ in  \cref{construction:quasi_coherent_sheaf_of_coefficients} provides an equivalence between this composite functor and the functor that to $V$ assigns $\pi_*(V \otimes_kE^{\otimes_k[n]})$.
\end{remark}

\begin{definition}
\label{definition:associated_homology_theory}
Let $\eta \colon k \to E$ be a map of commutative algebras in spectra such that the groupoid $S \simeq \Spec(\pi_*(E^{\otimes_k[-]}))$ has flat face maps. The homology theory associated with $\eta$ is the pair $(\mathcal{F},\partial)$ defined in \cref{construction:quasi_coherent_sheaf_of_coefficients}.
\end{definition}

\begin{theorem}
\label{theorem:e_2_term_of_descent_spectral_sequence_in_terms_of_cohomology}
Let $\eta \colon k \to E$ be a map of commutative algebras in spectra such that the groupoid in affine Dirac schemes $S \simeq \Spec(\pi_*(E^{\otimes_k[-]}))$ has flat face maps, and let $X \simeq |S|$ be the geometric Dirac stack given by its geometric realization. The descent spectral sequence for a $k$-module in spectra $V$ takes the form
$$\xymatrix@C=5mm{
{ E_{i,j}^2 \simeq H^{-i}(X,\mathcal{F}(V)(\nicefrac{j}{2})) } \ar@{=>}[r] &
{ \pi_{i+j}(\widehat{V}), } \cr
}$$
where $(\mathcal{F},\partial)$ is the homology theory associated with $\eta$.
\end{theorem}

\begin{proof}
Let $p \colon X \to 1$ be the unique map to the final Dirac stack $1 \simeq \Spec(\mathbb{Z})$. The coherent cohomology of any $\mathcal{F} \in \QCoh(X)$ is defined by
$$H^{-i}(X,\mathcal{F}) \simeq \pi_0\Map(\mathcal{O}_1,\Sigma^{-i}p_*(\mathcal{F})),$$
and we first explain how to use the affine submersion $f_0 \colon S_0 \to X$ to understand its structure. To this end, we consider the diagrams
$$\xymatrix@C=5mm{
{ S_n } \ar[rr]^-{f_n} \ar@<-.3ex>[dr]_-(.4){p_n} &&
{ X } \ar@<.3ex>[dl]^-(.4){p} \cr
{} &
{ 1 } &
{} \cr
}$$
and observe that the canonical map
$$\xymatrix{
{ \mathcal{F} } \ar[r] &
{ \varprojlim_{[n] \in \Delta} f_{n*}f_n{}^{\hspace{-4pt}*}(\mathcal{F}) } \cr
}$$
is an equivalence. Indeed, for any $\mathcal{G} \in \QCoh(X)$, we have a diagram of anima
$$\xymatrix{
{ \Map(\mathcal{G},\mathcal{F}) } \ar[r] \ar[d] &
{ \Map(\mathcal{G},\varprojlim_{[n] \in \Delta} f_{n*}f_n^*(\mathcal{F}))\phantom{,} } \ar[d] \cr
{ \varprojlim_{[n] \in \Delta}\Map(f_n^*(\mathcal{G}),f_n^*(\mathcal{F})) } \ar[r] &
{ \varprojlim_{[n] \in \Delta}\Map(\mathcal{G},f_{n*}f_n^*(\mathcal{F})), } \cr
}$$
where the left-hand vertical map is an equivalence, because the canonical map
$$\xymatrix{
{ \QCoh(X) } \ar[r] &
{ \varprojlim_{[n] \in \Delta} \QCoh(S_n) } \cr
}$$
is an equivalence of $\infty$-categories. Since also the right-hand vertical map and the lower horizontal maps are equivalences, we conclude that the top horizontal map is an equivalence, as desired. It follows that the canonical map
$$\xymatrix{
{ p_*(\mathcal{F}) } \ar[r] &
{ \varprojlim_{[n] \in \Delta} p_*f_{n*}f_n^*(\mathcal{F}) \simeq \varprojlim_{[n] \in \Delta} p_{n*}f_n^*(\mathcal{F}) } \cr
}$$
is an equivalence. Hence, we find that
$$\begin{aligned}
H^{-i}(X,\mathcal{F})
{} & \simeq \pi_0\Map(\mathcal{O}_1,\Sigma^{-i}p_*(\mathcal{F})) \simeq \pi_i\map(\mathcal{O}_1,p_*(\mathcal{F})) \cr
{} & \simeq \textstyle{
\pi_i(\varprojlim_{[n] \in \Delta}\map(\mathcal{O}_1,p_{n*}f_n^*(\mathcal{F})). } \cr
\end{aligned}$$
Now, suppose that $\mathcal{F} \in \QCoh(X)^{\heartsuit}$. Since $f_n \colon S_n \to X$ is a submersion from an affine Dirac stack, \cref{theorem:qcoh_on_geometric_dirac_stack_is_grothendieck_prestable} shows that $f_n^*(\mathcal{F}) \in \QCoh(S_n)^{\heartsuit}$, and since $p_n$ is affine flat, \cref{prop:pushforward_along_affine_flat_satisfies_base_change} further shows that $p_{n*}f_n^*(\mathcal{F}) \in \QCoh(1)^{\heartsuit}$. Thus, the mapping spectrum $\map(\mathcal{O}_1,p_{n*}f_n^*(\mathcal{F}))$ is the Eilenberg--MacLane spectrum given by the degree $0$ part $p_{n*}f_n^*(\mathcal{F})_0$  of the graded abelian group $p_{n*}f_n^*(\mathcal{F})$. It follows that the Bousfield--Kan spectral sequence
$$\begin{aligned}
E_{s,t}^2 & = \pi^{-s}([n] \mapsto \pi_t\map(\mathcal{O}_1,p_{n*}f_n^*(\mathcal{F}))) \cr
{} & \Rightarrow \textstyle { \pi_{s+t}(\varprojlim_{[n] \in \Delta} \map(\mathcal{O}_1,p_{n*}f_n^*(\mathcal{F}))) } \cr
\end{aligned}$$
is concentrated on the line $t = 0$, and therefore, it collapses to an isomorphism
$$\textstyle{ \pi^{-i}([n] \mapsto p_{n*}f_n^*(\mathcal{F})_0) \simeq \pi_i(\varprojlim_{[n] \in \Delta} \map(\mathcal{O}_1,p_{n*}f_n^*(\mathcal{F}))). }$$
Hence, for $\mathcal{F} \in \QCoh(X)^{\heartsuit}$, its coherent cohomology is given by
$$\textstyle{ H^{-i}(X,\mathcal{F}) \simeq \pi^{-i}([n] \mapsto p_{n*}f_n^*(\mathcal{F})_0), }$$
and more generally, the coherent cohomology of $\mathcal{F}(\nicefrac{j}{2})$ is given by
$$\textstyle{ H^{-i}(X,\mathcal{F}(\nicefrac{j}{2})) \simeq \pi^{-i}([n] \mapsto p_{n*}f_n^*(\mathcal{F})_j). }$$
Finally, if $\mathcal{F} \simeq \mathcal{F}(V)$, then $f_n^*(\mathcal{F}) \in \QCoh(S_n)^{\heartsuit}$ is the graded $\pi_*(E^{\otimes_k[n]})$-module $\pi_*(V \otimes_kE^{\otimes_k[n]})$, and $p_{n*}f_n^*(\mathcal{F}) \in \QCoh(1)^{\heartsuit}$ is its underlying graded abelian group, and therefore, we conclude that
$$H^{-i}(X,\mathcal{F}(V)(\nicefrac{j}{2})) \simeq \pi^{-i}([n] \mapsto \pi_j(V \otimes_kE^{\otimes_k[n]})) \simeq E_{i,j}^2,$$
as we wanted to prove.
\end{proof}

\begin{remark}
\label{remark:e_2_term_of_descent_spectral_sequence_in_terms_of_cohomology}
In the situation of \cref{theorem:e_2_term_of_descent_spectral_sequence_in_terms_of_cohomology}, the descent spectral sequence does not depend on the map of commutative algebras in spectra $\eta \colon k \to E$, but only on its associated homology theory $(\mathcal{F},\eta)$. We refer to \cite[Construction~2.24]{patchkoria2021adams} for the construction in these terms, but we show below that the homology theory $(\mathcal{F},\partial)$ is indeed adapted in the sense of~\cite[Definition~2.19]{patchkoria2021adams}, as the construction requires.
\end{remark}

In the situation of \cref{theorem:e_2_term_of_descent_spectral_sequence_in_terms_of_cohomology}, the counit of the adjunction
$$\xymatrix@C=10mm{
{ \Mod_k(\Sp) } \ar@<.7ex>[r]^-{\eta^*} &
{ \Mod_E(\Sp) } \ar@<.7ex>[l]^-{\eta_*} \cr
}$$
given by the extension of scalars along $\eta$ and the restriction of scalars along $\eta$ induces a map of graded $\pi_*(E)$-modules
$\pi_*(\eta^*\eta_*(W)) \to \pi_*(W)$, or equivalently, a map $f_0^*\mathcal{F}(\eta_*(W)) \to \mathcal{G}(W)$ in $\QCoh(S_0)^{\heartsuit}$, the mate of which is a natural map
$$\xymatrix@C=10mm{
{ \mathcal{F}(\eta_*(W)) } \ar[r] &
{ f_{0*}(\mathcal{G}(W)) } \cr
}$$
in $\QCoh(X)^{\heartsuit}$.

\begin{lemma}
\label{lemma:homology_of_an_e_module_is_extended}
In the situation of \cref{theorem:e_2_term_of_descent_spectral_sequence_in_terms_of_cohomology}, the canonical map
$$\xymatrix{
{ \mathcal{F}(\eta_*(W)) } \ar[r] &
{ f_{0*}(\mathcal{G}(W)) } \cr
}$$
is a natural isomorphism.
\end{lemma}

\begin{proof}
Since $\QCoh(X) \simeq \varprojlim_{[n] \in \Delta} \QCoh(S_n)$, it suffices to show that the map
$$\xymatrix@C=10mm{
{ f_n^*(\mathcal{F}(\eta_*(W))) } \ar[r] &
{ f_n^*f_{0*}(\mathcal{G}(W)) } \cr
}$$
induced by the map in the statement is an equivalence for all $n \geq 0$. In fact, it suffices to consider the case $n = 0$, but the argument becomes more transparent by considering all $n \geq 0$. We consider the diagram of finite ordinals
$$\xymatrix{
{ [n+1] } &
{ [0]*[n] } \ar@{=}[l] &
{ [0] } \ar[l] \cr
{ [n] } \ar[u]_-{d^0} &
{ \emptyset*[n] } \ar@{=}[l] \ar[u]_-{d^0 * [n]} &
{ \emptyset. } \ar[u]_-{d^0} \ar[l] \cr
}$$
The induced diagram of commutative algebras in spectra
$$\xymatrix@C=10mm{
{ E^{\otimes_k[n+1] } } &
{ E } \ar[l]_-{\eta_n'} \cr
{ E^{\otimes_k[n]} } \ar[u]_-(.45){\eta'} &
{ k } \ar[u]_-(.4){\eta} \ar[l]_-{\eta_n} \cr
}$$
is cocartesian, and the induced diagram of Dirac stacks
$$\xymatrix@C=12mm{
{ S_{n+1} } \ar[r]^-{f_n'} \ar[d]^-{f_0'} &
{ S_0 } \ar[d]^-{f_0} \cr
{ S_n } \ar[r]^-{f_n} &
{ X } \cr
}$$
is cartesian, because $S$ is a groupoid. Since $f_0$ is affine flat, \cref{prop:pushforward_along_affine_flat_satisfies_base_change} shows that the base-change map associated with the latter square is an equivalence, so it suffices to show that the composition
$$\xymatrix{
{ f_n^*f_{0*}(\mathcal{G}(W)) } \ar[r] &
{ f_n^*f_{0*}(\mathcal{G}(W)) } \ar[r] &
{ f_{0*}'f_n'{}^{\hspace{-1pt}*}(\mathcal{G}(W)) } \cr
}$$
of the map in question with the base-change map is an equivalence for all $n \geq 0$. Changing to graded module notation, we have a diagram
$$\xymatrix@C=-4mm{
{ \pi_*(\eta_n^*\eta_*(W)) } \ar[rr] \ar[dr] &&
{ f_{0*}'f_n'{}^{\hspace{-1pt}*}(\pi_*(W)) } \cr
{} &
{ \pi_*(\eta_*'\eta_n'{}^{\hspace{-2pt}*}(W)) } \ar[ur] &
{} \cr
}$$
of graded $\pi_*(E^{\otimes_k[n]})$-modules, and we wish to prove that the top map is an equivalence. The left-hand slanted map is the base-change map associated with the above cocartesian square of commutative algebras in spectra, so it is an isomorphism, and the right-hand slanted map is the K\"{u}nneth map, which is an isomorphism, because $\pi_*(\eta_n')$ is a flat map of Dirac rings by the assumption that $S$ has flat face maps.
\end{proof}

\begin{proposition}
\label{proposition:associated_homology_theory_is_adapted}
In the situation of \cref{theorem:e_2_term_of_descent_spectral_sequence_in_terms_of_cohomology}, the homology theory $(\mathcal{F},\partial)$ associated with $\eta \colon k \to E$ is adapted.
\end{proposition}

\begin{proof}
We recall that, by \cite[Definition~2.19]{patchkoria2021adams}, the homotopy theory $(\mathcal{F},\partial)$ is adapted if for every injective object $\mathcal{I} \in \QCoh(X)^{\heartsuit}$, there is an isomorphism $\phi \colon \mathcal{F}(I) \to \mathcal{I}$ with $I \in \Mod_k(\Sp)$ such that for every $V \in \Mod_k(\Sp)$, the composite map
$$\xymatrix{
{ \pi_0\Map(V,I) } \ar[r] &
{ \Map(\mathcal{F}(V),\mathcal{F}(I)) } \ar[r] &
{ \Map(\mathcal{F}(V),\mathcal{I}) } \cr
}$$
is an isomorphism. We notice that the class of injective objects $\mathcal{I} \in \QCoh(X)^{\heartsuit}$ for which such an isomorphism $\phi \colon \mathcal{F}(I) \to \mathcal{I}$ exists is closed under retracts.

Now, since $f \simeq f_0 \colon S_0 \to X$ is affine flat, both functors in the adjunction
$$\xymatrix@C=10mm{
{ \QCoh(X) } \ar@<.7ex>[r]^-{f^*} &
{ \QCoh(S_0) } \ar@<.7ex>[l]^-{f_*} \cr
}$$
are $t$-exact, so they restrict to an adjunction
$$\xymatrix@C=10mm{
{ \QCoh(X)^{\heartsuit} } \ar@<.7ex>[r]^-{f^*} &
{ \QCoh(S_0)^{\heartsuit} } \ar@<.7ex>[l]^-{f_*} \cr
}$$
with both functors exact. Moreover, the unit $\eta \colon \mathcal{F} \to f_*f^*(\mathcal{F})$ is a monomorphism, because $f \colon S_0 \to X$ is an effective epimorphism. Thus, every injective object on the left-hand side is a retract of an injective object of the essential image of $f_*$.

So we let $\mathcal{I} \simeq f_*(\mathcal{J})$ with $\mathcal{J} \in \QCoh(S_0)^{\heartsuit}$ injective and proceed to find an isomorphism $\phi \colon \mathcal{F}(I) \to \mathcal{I}$ with $I \in \Mod_k(\Sp)$ with the required property. It follows from Brown representability, \cite[Theorem~1.4.1.2]{lurieha}, that the functor
$$\xymatrix@C=20mm{
{ h\Mod_E(\Sp) } \ar[r]^-{\Map(\mathcal{G}(-),\mathcal{J})} &
{ \Set } \cr
}$$
with $\mathcal{G} \colon \Mod_E(\Sp) \to \QCoh(S_0)^{\heartsuit}$ as in \cref{lemma:homology_of_an_e_module_is_extended} is representable. Therefore, we obtain a map $\psi \colon \mathcal{G}(J) \to \mathcal{J}$ with $J \in \Mod_E(\Sp)$ such that the composite map
$$\xymatrix{
{ \pi_0\Map(-,J) } \ar[r] &
{ \Map(\mathcal{G}(-),\mathcal{G}(J)) } \ar[r] &
{ \Map(\mathcal{G}(-),\mathcal{J}) } \cr
}$$
is a natural isomorphism, and since $\QCoh(S_0)^{\heartsuit}$ has a family of generators contained in the essential image of $\mathcal{G}$, we conclude that the map $\psi$ is an isomorphism. Accordingly, we let $I \simeq \eta_*(J)$ and define $\phi \colon \mathcal{F}(I) \to \mathcal{I} \simeq f_*(\mathcal{J})$ to be the composition
$$\xymatrix{
{ \mathcal{F}(I) } \ar[r] &
{ f_*(\mathcal{G}(J)) } \ar[r] &
{ f_*(\mathcal{J}) } \cr
}$$
of the isomorphism provided by \cref{lemma:homology_of_an_e_module_is_extended} and the isomorphism induced by $\psi$. Since the composite
$$\xymatrix{
{ f^*\mathcal{F}(V) } \ar[r] &
{ f^*\mathcal{F}(\eta_*\eta^*(V)) } \ar[r] &
{ f^*f_*\mathcal{G}(\eta^*(V)) } \ar[r] &
{ \mathcal{G}(\eta^*(V)) } \cr
}$$
is an isomorphism for all $V \in \Mod_k(\Sp)$, we conclude that the composite map
$$\xymatrix{
{ \pi_0\Map(V,I) } \ar[r] &
{ \Map(\mathcal{F}(V),\mathcal{F}(I)) } \ar[r] &
{ \Map(\mathcal{F}(V),\mathcal{I}) } \cr
}$$
is an isomorphism, as desired. This completes the proof.
\end{proof}

\begin{remark}
\label{remark:synthetic_deformation}
We briefly discuss the second author's categorification of the descent spectral sequence associated with $\eta \colon k \to E$. By \cref{proposition:associated_homology_theory_is_adapted}, the homology theory $(\mathcal{F},\partial)$ associated with $\eta \colon k \to E$ is adapted, so \cite[Theorem~6.40]{patchkoria2021adams} exhibits the homological functor $\mathcal{F}$ as the $0$th truncation
$$\xymatrix@C=0mm{
{} & 
{ \mathcal{D}(\mathcal{F})_{\geq 0} } \ar@{->>}[dr]^-{\pi_0^{\heartsuit}} &
{} \cr
{ \Mod_k(\Sp) } \ar@{^{(}->}[ur]^-{\nu} \ar[rr]^-{\phantom{,}\mathcal{F}\phantom{}} &&
{ \QCoh(X)^{\heartsuit} } \cr
}$$
of a fully faithful embedding as a reflexsive subcategory in a Grothendieck prestable $\infty$-category.\footnote{\,As for the derived $\infty$-category of an Grothendieck abelian category, there are several variants of the derived $\infty$-category attached to $(\mathcal{F},\partial)$; compare \cite[Remarks 6.38 and 6.39]{patchkoria2021adams}. We consider the unseparated variant only.} Moreover, if we exhibit $\mathcal{D}(\mathcal{F})_{\geq 0}$ as the connective part of a presentable stable $\infty$-category $\mathcal{D}(\mathcal{F})$ with a $t$-structure compatible with filtered colimits, then~\cite[Theorem~5.60]{patchkoria2021adams} identifies the descent spectral sequence in \cref{theorem:e_2_term_of_descent_spectral_sequence_in_terms_of_cohomology} with the spectral sequence associated to the filtered spectrum\footnote{\,The slice filtration on the stable motivic $\infty$-category encodes the motivic spectral sequence $\smash{ E_{i,j}^2 = H^{-i}(X,\mathbb{Z}(\nicefrac{j}{2})) \Rightarrow K_{i+j}(X) }$ in a similar manner. In the case of the descent spectral sequence, the role of the slice filtration is played by the $t$-structure on $\mathcal{D}(\mathcal{F})$, as the tower
$$\xymatrix@C=7mm{
{ \cdots } \ar[r]^-{\tau} &
{ \Sigma \nu (\Sigma^{-1} V) } \ar[r]^-{\tau} &
{ \nu(V) } \ar[r]^-{\tau} &
{ \Sigma^{-1} \nu (\Sigma V) } \ar[r]^-{\tau} &
{ \cdots }
}$$
is a Postnikov tower with respect to this $t$-structure.}
$$F^{\star}(V) \simeq \map(\nu(k),\Sigma^{\star}\nu(\Sigma^{-\star}V)) \in \Fun(\mathbb{Z}^{\op},\Sp)$$
with structure maps given by the colimit interchange maps $\tau \colon \Sigma\nu(W) \to \nu(\Sigma W)$. However, the Grothendieck prestable $\infty$-category $\mathcal{D}(\mathcal{F})_{\geq 0}$ encodes a great deal more information than does the descent spectral sequence, whence its usefulness.
\end{remark}

\begin{remark}
\label{remark:module_with_descent_data}
In \cref{construction:quasi_coherent_sheaf_of_coefficients}, we used the equivalences
$$\xymatrix{
{ \QCoh(X)^{\heartsuit} } \ar[r] &
{ \varprojlim_{\Delta} \QCoh(S)^{\heartsuit} } \ar[r] &
{ \varprojlim_{\Delta_s} \QCoh(S)^{\heartsuit} } \cr
}$$
where, in the middle term, it is important that we use the truncated extension of scalars maps $\theta_{\heartsuit}^*$ associated with $\theta \colon S_m \to S_n$, because not every such map $\theta$ is flat. Since the diagram $\QCoh(S)^{\heartsuit}$ is a diagram of $1$-categories, we conclude from~\cite[Proposition~A.1]{diracgeometry1} that the restriction map
$$\xymatrix{
{ \varprojlim_{\Delta}\QCoh(S)^{\heartsuit} } \ar[r] &
{ \varprojlim_{\Delta^{\leq 2}}\QCoh(S)^{\heartsuit} } \cr
}$$
is an equivalence. An object of the right-hand side, in turn, is a pair $(\mathcal{F}_0,\epsilon)$ of an object $\mathcal{F}_0 \in \QCoh(S_0)^{\heartsuit}$ and a ``descent structure'' map\footnote{\,The map $\epsilon$ is necessarily an isomorphism.}
$$\xymatrix{
{ d_1^*(\mathcal{F}_0) } \ar[r]^-{\epsilon} &
{ d_0^*(\mathcal{F}_0) } \cr
}$$
in $\QCoh(S_1)^{\heartsuit}$ that makes the diagrams
$$\begin{xy}
(-2,12)*+{ s_{0,\heartsuit}^*d_1^*(\mathcal{F}_0) }="1";
(32,12)*+{ s_{0,\heartsuit}^*d_0^*(\mathcal{F}_0) }="2";
(15,0)*+{ \mathcal{F}_0 }="3";
{ \ar^-{s_{0,\heartsuit}^*(\epsilon)} "2";"1";};
{ \ar@{=} "3";"1";};
{ \ar@{=} "3";"2";};
\end{xy}$$
$$\begin{xy}
(0,12)*+{ d_1^*d_1^*(\mathcal{F}_0) }="11";
(30,12)*+{ d_1^*d_0^*(\mathcal{F}_0) }="12";
(-15,0)*+{ d_2^*d_1^*(\mathcal{F}_0) }="21";
(45,0)*+{ d_0^*d_0^*(\mathcal{F}_0) }="22";
(0,-12)*+{ d_2^*d_0^*(\mathcal{F}_0) }="31";
(30,-12)*+{ d_0^*d_1^*(\mathcal{F}_0) }="32";
{ \ar^-{d_1^*(\epsilon)} "12";"11";};
{ \ar@{=} "21";"11";};
{ \ar@{=} "12";"22";};
{ \ar_-{d_2^*(\epsilon)} "31";"21";};
{ \ar_-{d_0^*(\epsilon)} "22";"32";};
{ \ar@{=} "32";"31";};
\end{xy}$$
in $\QCoh(S_0)^{\heartsuit}$ and $\QCoh(S_2)^{\heartsuit}$ commute. The equality signs indicate the various unique isomorphisms. We say that $(\mathcal{F}_0,\epsilon)$ is a quasi-coherent $\mathcal{O}_{S_0}$-module with descent structure along $f \colon S_0 \to X$, and we refer to the fact that the map $\epsilon$ makes these diagrams commute by saying that it satisfies the cocycle condition. A map $\alpha \colon (\mathcal{F}_0,\epsilon) \to (\mathcal{F}_0',\epsilon')$ of quasi-coherent $\mathcal{O}_{S_0}$-modules with descent structure along the map $f \colon S_0 \to X$ is a map $\alpha \colon \mathcal{F}_0 \to \mathcal{F}_0'$ in $\QCoh(S_0)^{\heartsuit}$ that makes the diagram
$$\xymatrix@C+=12mm{
{ d_1^*(\mathcal{F}_0) } \ar[r]^-{\epsilon} \ar[d]^-{d_1^*(\alpha)} &
{ d_0^*(\mathcal{F}_0) } \ar[d]^-{d_0^*(\alpha)} \cr
{ d_1^*(\mathcal{F}_0') } \ar[r]^-{\epsilon'} &
{ d_0^*(\mathcal{F}_0') } \cr
}$$
in $\QCoh(S_1)^{\heartsuit}$ commute.
\end{remark}

\begin{remark}
\label{remark:module_with_flat_connection}
Suppose that $(\mathcal{F}_0,\epsilon)$ is a quasi-coherent $\mathcal{O}_{S_0}$-module with descent data along $f \colon S_0 \to X$ as in \cref{remark:module_with_descent_data}. The map of quasi-coherent $\mathcal{O}_{S_1}$-modules
$$\xymatrix@C=10mm{
{ d_1^*(\mathcal{F}_0) } \ar[r]^-{\epsilon} &
{ d_0^*(\mathcal{F}_0) } \cr
}$$
determines and is determined by the map of quasi-coherent $\mathcal{O}_{S_0}$-modules
$$\xymatrix@C=10mm{
{ \mathcal{F}_0 } \ar[r]^-{\nabla_{\epsilon}} &
{ d_{1*}d_0^*(\mathcal{F}_0) } \cr
}$$
given by its mate. The cocycle condition for $\epsilon$ translates to a Leibniz rule for the map $\nabla_{\epsilon}$, which behaves as a connection.\footnote{\,The map $\nabla_{\epsilon}$ is a comodule structure for the descent comonad $f^*\hspace{-.5pt}f_* \simeq d_{1*}d_0^*$.}
\end{remark}

\subsection{Quillen's theorem}
\label{subsection:quillens_thm}

Given a commutative algebra in spectra $E$ and a group in anima $G$, we defined, in the introduction, a formal Dirac scheme
$$\xymatrix{
{ Y \simeq Y_{E,G} } \ar[r]^-{q} &
{ S \simeq S_E } \cr
}$$
over $S \simeq \Spec(R)$ with $R \simeq \pi_*(E)$. We consider the case, where
$$G \simeq \widehat{L} \simeq \Hom(L,U(1))$$
is the Pontryagin dual of a finitely generated free abelian group $L$. The Postnikov filtration of $E$ defines the Atiyah--Hirzebruch spectral sequence converging to
$$A \simeq \pi_*(p_*p^*(E)) \simeq \pi_*(\map(\Sigma_+^{\infty}BG,E)),$$
and since the integral cohomology groups of $BG$ are free, it takes the form
$$\xymatrix@C=5mm{
{ E^2 \simeq \Sym_R(L \otimes_{\mathbb{Z}}R(1)) } \ar@{=>}[r] &
{ A. } \cr
}$$
We define $I \subset A$ to be the kernel of the edge homomorphism $\theta \colon A \to R$ and recall that $E$ is said to be complex orientable, if the canonical map
$$\xymatrix{
{ E_{-2,*}^{\infty} \simeq I/I^2 } \ar[r] &
{ E_{-2,*}^2 \simeq L \otimes_{\mathbb{Z}}R(1) } \cr
}$$
is an isomorphism. In this case, all differentials in the spectral sequence are necessarily zero, so we conclude that $q \colon Y \to S$ is a formal hyperplane.

\begin{definition}
\label{definition:quillen_formal_group}
Let $E$ be a complex orientable commutative algebra in spectra. The Quillen formal group associated with $E$ is the functor
$$\xymatrix@C=12mm{
{ \Lat^{\op} } \ar[r]^-{\mathcal{G}_E^Q} &
{ \Hyp(S) } \cr
}$$
that to $L$ assigns the formal hyperplane $q \colon Y_{E,\widehat{L}} \to S$.
\end{definition}

\begin{remark}
\label{remark:quillen_formal_group_has_canonically_trivialized_lie_algebra}
Let $E$ be a complex orientable commutative algebra in spectra. The canonical isomorphism $I/I^2 \to R(1)$ provided by the Atiyah--Hirzebruch spectral sequence in the case $L = \mathbb{Z}$ defines a canonical isomorphism
$$\xymatrix@C=12mm{
{ \mathbb{A}_S(1) } \ar[r]^-{\phi_E^Q} &
{ \Lie(\mathcal{G}_E^Q) } \cr
}$$
of line bundles over $S$. This isomorphism provides a canonical trivialization of the Lie algebra, which we defined in \cref{def:liealgebra}, of the Quillen formal group. In particular, the Quillen formal group $\mathcal{G}_E^Q$ is $1$-dimensional and spin-$1$ over $S$.
\end{remark}

We now let $E \simeq \operatorname{MU}^{\otimes [n]}$ be a tensor power of the commutative algebra in spectra representing complex cobordism. It is proved in~\cite[Lemma~4.5 and Theorem~8.1]{adams1} that the homotopy groups of $E$ are concentrated in even degrees, so $E$ is complex orientable. Therefore, the Quillen formal group provides the simplicial formal group
$$\xymatrix@C+=18mm{
{ \Delta^{\op} } \ar[r]^-{\mathcal{G}_{\operatorname{MU}^{\otimes[-]}}^Q } &
{ \FGroup, } \cr
}$$
whose value at $L \in \Lat$ is the simplicial formal hyperplane
$$\xymatrix{
{ \Spf(\pi_*(\map(\Sigma_+^{\infty}B\widehat{L},\operatorname{MU}^{\otimes[-]}))) } \ar[r] &
{ \Spec(\pi_*(\operatorname{MU}^{\otimes[-]})). } \cr
}$$
It satisfies the hypotheses of \cref{prop:formalgroups}, so the geometric realization
$$\mathcal{G}^Q \simeq \lvert \mathcal{G}_{\operatorname{MU}^{\otimes[-]}}^Q \rvert \in \FGroup$$
exists and its value at $L \in \Lat$ is the formal hyperplane
$$\xymatrix{
{ |\Spf(\pi_*(\map(\Sigma_+^{\infty}B\widehat{L},\operatorname{MU}^{\otimes[-]})))| } \ar[r] &
{ |\Spec(\pi_*(\operatorname{MU}^{\otimes[-]}))|. } \cr
}$$
We write $X$ for the Dirac stack on the right-hand side so that $\mathcal{G}^Q \in \FGroup(X)$. By the descent of modules along effective epimorphisms, the canonical trivializations of \cref{remark:quillen_formal_group_has_canonically_trivialized_lie_algebra} determine a canonical isomorphism of line bundles over $X$,
$$\xymatrix@C=10mm{
{ \mathbb{A}_X(1) } \ar[r]^-{\phi^Q}  &
{ \Lie(\mathcal{G}^Q). } \cr
}$$
We now use Quillen's theorem to show that $(\mathcal{G}^Q,\phi^Q)$ is the universal $1$-dimensional spin-$1$ formal group equipped with a trivialization of its Lie algebra in the following sense.

\begin{theorem}
\label{thm:quillen}
Let $T$ be a Dirac stack. The functor that to $f \colon T \to X$ assigns the pair $(f^*\mathcal{G}^Q,f^*\phi^Q)$ is an equivalence from $\Map(T, X)$ to the groupoid of pairs $(\mathcal{G},\phi)$ of a formal group $\mathcal{G}$ over $T$ and an isomorphism $\phi \colon \mathbb{A}_T(1) \to \Lie(\mathcal{G})$ of line bundles over $T$.
\end{theorem}

\begin{proof}
We let $p_n \colon Y_n \to S_n$ be the underlying formal hyperplane of the Quillen formal group associated with $\operatorname{MU}^{\otimes[n]}$ and note that a complex orientation of $\operatorname{MU}$ determines and is determined by an isomorphism
$$\begin{xy}
(0,7)*+{ Y_0 }="11";
(22,7)*+{ \affinespace_{S_0}(1) }="12";
(0,-7)*+{ S_0 }="21";
(22,-7)*+{ S_0 }="22";
{ \ar^-{h} "12";"11";};
{ \ar_-{p} "21";"11";};
{ \ar_-{q} "22";"12";};
{ \ar@/_.7pc/_-{0} "11";"21";};
{ \ar@/_.7pc/_-{0} "12";"22";};
{ \ar@{=} "22";"21";};
\end{xy}$$
of pointed formal hyperplanes over $S_0$ such that the diagram
$$\xymatrix@C=1mm{
{} &
{ \mathbb{A}_{S_0}(1) } \ar[dl]_-{\phi_{\operatorname{MU}}^Q} \ar[dr]^-{\alpha} &
{} \cr
{ \;\;\; T_{Y_0/S_0,0} } \ar[rr]^{dh_0} &&
{ T_{\affinespace_{S_0}(1)/S_0,0} } \cr
}$$
of line bundles over $S_0$ commutes. Here $\phi_{\operatorname{MU}}^Q$ and $\alpha$ are the canonical trivializations provided by \cref{remark:quillen_formal_group_has_canonically_trivialized_lie_algebra} and \cref{proposition:formal_completion_is_formally_complete}. Moreover, under the identification provided by a choice of complex orientation $h$, the Quillen formal group $\mathcal{G}_{\operatorname{MU}}^Q$ is encoded by a formal group law $f(y,z) \in \pi_*(MU)[[y,z]]$, which, in turn, determines and is determined by a map of affine Dirac schemes
$$\xymatrix{
{ \Spec(\pi_*(\operatorname{MU})) } \ar[r]^-{\theta} &
{ \Spec(L) } \cr
}$$
with $L$ the Lazard ring. Quillen shows in~\cite[Theorem~2]{quillen6} that $\theta$ is an isomorphism; see also~\cite[Theorem~II.8.2]{adams1} for a purely homotopy theoretic proof of this fact.

Next, we note that the face maps $d_0,d_1 \colon Y_1 \to Y_0$ and the complex orientation $h \colon Y \to \affinespace_{S_0}(1)$ of $\operatorname{MU}$ give rise to two complex orientations of $\operatorname{MU}^{\otimes[1]}$,
$$\begin{xy}
(0,0)*+{ {}^{\phantom{1}}Y_1 }="1";
(21,0)*+{ \affinespace_{S_1}(1). }="2";
{ \ar@<.7ex>^-{h_0} "2";"1";};
{ \ar@<-.7ex>_-{h_1} "2";"1";};
\end{xy}$$
The composite map $g \simeq h_0 \circ h_1^{-1}$ is an automorphism of the pointed spin-$1$ formal affine line $\affinespace_{S_1}(1)$ that induces the identity automorphism on conormal sheaves. It defines an $S_1$-valued point $g$ of the affine group scheme $G'$ of automorphisms of the pointed spin-$1$ formal affine line $\affinespace_{\mathbb{Z}}(1)$ which induce the identity map of tangent spaces at the zero section.\footnote{\,We have studied $G'$ earlier. Indeed, by \cref{lemma:automorphism_of_inf_1_same_as_those_of_conormal_sheaf}, we have $G' \simeq \mathrm{ker}(G \to G^1)$.} Quillen shows that the map of affine Dirac schemes
$$\xymatrix@C=12mm{
{ S_1 } \ar[r]^-{(g,d_1)} &
{ G' \times S_0 } \cr
}$$
is an isomorphism. Moreover, there is are commutative diagrams
$$\begin{xy}
(0,12)*+{ S_1 }="11";
(24,12)*+{ G' \times S_0 }="13";
(12.5,0)*+{ S_0 }="22";
{ \ar^-{(g,d_1)} "13";"11";};
{ \ar_-(.45){d_0} "22";"11";};
{ \ar^-(.45){\mu} "22";"13";};
(45,12)*+{ S_1 }="14";
(69,12)*+{ G' \times S_0 }="16";
(57.5,0)*+{ S_0 }="25";
{ \ar^-{(g,d_1)} "16";"14";};
{ \ar_-(.45){d_1} "25";"14";};
{ \ar^-(.45){p} "25";"16";};
\end{xy}$$
where the action map $\mu$ is given as follows. Let $T = \Spec(B)$ be any affine Dirac scheme. A $T$-valued point of $S_0$ determines and is determined by a formal group law $f = f(y,z) \in B[[y,z]]$, and a $T$-valued point of $G'$ determines and is determined by a power series $g = g(t) \in B[[t]]$ with $g(0) = 0$ and $g'(0) = 1$. Here, the generators $y$, $z$, and $t$ are all spin-$1$, and $f$ and $g$ are also required to be spin-$1$. Now, the image of $(f,g)$ by $\mu$ is the formal group law $g(f(g^{-1}(y),g^{-1}(z)))$. It follows, in particular, that the groupoid $S \colon \Delta^{\op} \to \sheaves(\Aff)$ has flat face maps, so by \cref{theorem:geometricity_in_terms_of_groupoids}, this its geometric realization $X \simeq \lvert S \rvert$ is a geometric stack.

Let $X'(T)$ be the groupoid consisting of pairs $(\mathcal{G},\phi)$ of a formal group $\mathcal{G}$ over $T$ and an isomorphism $\phi \colon \mathbb{A}_T(1) \to \Lie(\mathcal{G})$ of vector bundles over $T$. We wish to prove that the functor
$$\xymatrix{
{ X(T) } \ar[r] &
{ X'(T) } \cr
}$$
that to $f \colon T \to X$ assigns $(f^*\mathcal{G}^Q,f^*\phi^Q)$ is an equivalence. To this end, we consider the colimit
$$\mathcal{G}'' \simeq |\mathcal{G}_{\operatorname{MU}^{\otimes [-]}}| \in \FGroup(X'')$$
calculated in presheaves, so that $X$ is the sheafification of $X''$. As colimits of presheaves are calculated pointwise, we have $X''(T) \simeq \lvert S(T) \rvert \in \spaces$, and Quillen's theorem shows that the functor
$$\xymatrix{
{ X''(T) } \ar[r] &
{ X'(T) } \cr
}$$
that to $f \colon T \to X''$ assigns $(f^*\mathcal{G}'',f^*\phi'')$ is fully faithful and that its essential image is the subgroupoid spanned by the pairs $(\mathcal{G},\phi)$ such that $\mathcal{G}$ admits a global coordinate in the sense that $\mathcal{G}(\mathbb{Z})$ is isomorphic as a formal hyperplane over $T$ to the spin-$1$ formal affine line. Hence, to prove that $X(T) \to X'(T)$ is an equivalence, we must show that the presheaf $X'$ is a sheaf. But this is precisely flat descent for formal groups, which we have proved in \cref{proposition:formal_groups_takes_colimits_of_stacks_to_catinfty_limits}.
\end{proof}

\subsection{Milnor's theorem}

We next let $p$ be an odd prime number, and let 
$$\xymatrix{
{ G \simeq \widehat{L}[p] \simeq \Hom(L,C_p) } \ar[r]^-{i} &
{ \widehat{L} \simeq \Hom(L,U(1)) } \cr
}$$
is the $p$-torsion subgroup of the Pontryagin dual of a finitely generated free abelian group $L$. We claim that if $\eta \colon \mathbb{F}_p \to E$ is a map of commutative algebras in spectra,\footnote{\,By the Hopkins--Mahowald theorem \cite[Theorem~4.16]{mathew2015nilpotence}, if $E$ is an $\mathbf{E}_2$-algebra in spectra and if $p = 0$ in $\pi_*(E)$, then there exists a map of $\mathbf{E}_2$-algebras in spectra $\eta \colon \mathbb{F}_p \to E$.} then the formal Dirac scheme
$$\xymatrix{
{ Y \simeq Y_{E,G} } \ar[r]^-{q} &
{ S \simeq S_E } \cr
}$$
over $S \simeq \Spec(R)$ with $R \simeq \pi_*(E)$, which we defined in the introduction, is a formal hyperplane. Since $p$ is odd, the $\mathbb{F}_p$-cohomology of $BC_p$ is a free Dirac $\mathbb{F}_p$-algebra on a spin-$\nicefrac{1}{2}$ generator $e$ and a spin-$1$ generator $\gamma$. It follows that the Atiyah--Hirzebruch spectral sequence converging to
$$A \simeq \pi_*(p_*p^*(E)) \simeq \pi_*(\map(\Sigma_+^{\infty}BG,E))$$
takes the form
$$\xymatrix@C=5mm{
{ E^2 \simeq \Sym_R(L \otimes_{\mathbb{Z}}R(e,\gamma)) } \ar@{=>}[r] &
{ A, } \cr
}$$
and we claim that all differentials are zero. Indeed, the spectral sequence is a spectral sequence of $R$-algebras, and the generators of the symmetric Dirac $R$-algebra are all defined over $\mathbb{F}_p$, where the differentials vanish for degree-reasons. This proves our claim that $q \colon Y \to S$ is a formal hyperplane.

\begin{definition}
\label{definition:milnor_formal_group}
Let $E$ be a commutative $\mathbb{F}_p$-algebra in spectra, where $p$ is an odd prime number. The Milnor formal group associated with $E$ is the functor
$$\xymatrix@C=12mm{
{ \Lat^{\op} } \ar[r]^-{\mathcal{G}_E^M} &
{ \Hyp(S) } \cr
}$$
that to $L$ assigns the formal hyperplane $q \colon Y_{E,\widehat{L}[p]} \to S$.
\end{definition}

\begin{remark}
\label{remark:milnor_formal_group_is_additive}
The Milnor formal group associated with a commutative $\mathbb{F}_p$-algebra in spectra $E$ is additive in the sense that, as a formal group over $S$, it is non-canonically isomorphic to an additive formal group.
\end{remark}

If $E$ is a commutative $\mathbb{F}_p$-algebra in spectra with $p$ odd, then $i \colon \widehat{L}[p] \to \widehat{L}$ induces a map of formal groups from the $2$-dimensional Milnor formal group to the $1$-dimensional Quillen formal group,
$$\xymatrix{
{ \mathcal{G}_E^M } \ar[r] &
{ \mathcal{G}_E^Q. } \cr
}$$
This map is surjective and its kernel
$$\xymatrix{
{ \Fil^1\mathcal{G}_E^M } \ar[r] &
{ \Fil^0\mathcal{G}_E^M \simeq \mathcal{G}_E^M } \cr
}$$
is again a formal group over $S$. This filtration gives rise to a grading of the Lie algebra $\Lie(\mathcal{G}_E^M)$, which we refer to as the grading by charge. The Atiyah--Hirzebruch spectral sequence determines a canonical isomorphism of graded vector bundles
$$\xymatrix@C+=10mm{
{ \mathbb{A}_S(e,\gamma) } \ar[r]^-{\phi^M} &
{ \gr(\Lie(\mathcal{G}_E^M)), } \cr
}$$
where $e$ has charge $1$ and $\gamma$ has charge $0$. As in \cref{remark:classifying_stack_for_automorphisms_of_object_in_slice_category}, we have a classifying stack $$\Aut(\Fil\mathcal{G}_E^M) \in \Grp_{\mathbf{E}_1}(\sheaves(\Aff)_{/S})$$
for automorphisms of the Milnor formal group associated with $E$ that preserve the canonical filtration and induce the identity map on $\gr(\Lie(\mathcal{G}_E^M))$.

\begin{proposition}
\label{prop:milnor}
Let $S_0 \simeq \Spec(\mathbb{F}_p)$ with $p$ prime. The classifying stack
$$\Aut(\Fil\mathcal{G}_{\mathbb{F}_p}^M) \in \Grp_{\mathbf{E}_1}(\sheaves(\Aff)_{S_0})$$
is an affine Dirac group $S_0$-scheme represented by the Hopf algebra
$$C = \mathbb{F}_p[\tau_0,\tau_1,\tau_2,\dots,\xi_1,\xi_2,\dots]$$
with $\tau_i$ and $\xi_i$ generators of spin $p^i - \nicefrac{1}{2}$ and $p^i - 1$, respectively, and with comultiplication given by
$$\begin{aligned}
\psi(\tau_i) & = \textstyle{ \tau_i \otimes 1 + \sum_{0 \leq j \leq i} \xi_{i-j}^{p^j} \otimes \tau_j } \cr
\psi(\xi_i) & = \textstyle{ \sum_{0 \leq j \leq i} \xi_{i-j}^{p^j} \otimes \xi_j, } \cr
\end{aligned}$$
 where $\xi_0 = 1$ by convention.
\end{proposition}

\begin{proof}Let $\mathcal{G}_0$ be the Milnor formal group associated with $\mathbb{F}_p$, and let $q_0 \colon Y_0 \to S_0$ be its underlying hyperplane $\mathcal{G}_0(\mathbb{Z})$. So $Y_0 \simeq \Spf(A_0)$ with $A_0 = \mathbb{F}_p[e,\gamma]$ and with the zero section defined by the graded ideal $I_0 = (e,\gamma)$. Now, given any Dirac $\mathbb{F}_p$-algebra $R$, we have
$$\begin{aligned}
\Hom_{S_0}(Y_0,Y_0)(R) & \simeq \Hom_R(\Spf(A_0 \otimes_{\mathbb{F}_p}R),\Spf(A_0 \otimes_{\mathbb{F}_p}R)) \cr {} & \simeq \Hom_R(R[e,\gamma],R[[e,\gamma]]),
\end{aligned}$$
so $f \in \Hom_{S_0}(Y_0,Y_0)(R)$ determines and is determined by the two power series
$$\begin{aligned}
f(e) & = \textstyle{ \sum_{i \geq 0}(a_i + c_ie)\gamma^i } \cr
f(\gamma) & = \textstyle{ \sum_{i \geq 0}(b_i + d_ie)\gamma^i } \cr
\end{aligned}$$
with $a_i \in R_{2i-1}$, $b_i \in R_{2i-2}$, $c_i \in R_{2i}$, and $d_i \in R_{2i-1}$. Now, let $T \simeq \Spec(R)$ and $Y_{0,T} \simeq Y_0 \times_{S_0}T$. The statement that $f$ is a map of formal groups is equivalent to the statement that the diagram
$$\xymatrix{
{ Y_{0,T} \times_T Y_{0,T} } \ar[r]^-{+} \ar[d]^-{f \times f} &
{ Y_{0,T} } \ar[d]^-{f} \cr
{ Y_{0,T} \times_T Y_{0,T} } \ar[r]^-{+} &
{ Y_{0,T} } \cr
}$$
commutes. Each of the two composites in the diagram are elements of
$$\Hom_{S_0}(Y_0 \times_SY_0,Y_0)(R) \simeq \Hom_R(R[e,\gamma],R[[e_1,\gamma_1,e_2,\gamma_2]]),$$
and by comparing the images of $e$ and $\gamma$ by the two composites, we see that the diagram commutes if and only if the following equations hold:
$$\begin{aligned}
\textstyle{ \sum_{i \geq 0} (a_i + c_i(e_1 + e_2))(\gamma_1 + \gamma_2)^i }
& = \textstyle{ \sum_{i \geq 0} (a_i + c_ie_1)\gamma_1^i + \sum_{i \geq 0}(a_i + c_ie_2)\gamma_2^i, } \cr
\textstyle{ \sum_{i \geq 0} (b_i + d_i(e_1 + e_2))(\gamma_1 + \gamma_2)^i }
& = \textstyle{ \sum_{i \geq 0} (b_i + d_ie_1)\gamma_1^i + \sum_{i \geq 0}(b_i + d_ie_2)\gamma_2^i. } \cr
\end{aligned}$$
These are satisfied if and only if $c_i$ and $d_i$ are zero for all $i > 0$ and $a_i$ and $b_i$ are zero for all $i \geq 0$ that are not a power of $p$. Moreover, the map $f$ preserves the canonical filtration of the Milnor formal group if and only if $d_0 = 0$, and, finally, the induced map of $\gr(\Lie(\mathcal{G}_0))$ is the identity map if and only if $c_0 = 1$ and $b_1 = 1$. Hence, a choice of $e$ and $\gamma$ determines an identification of the set of $R$-valued points of $\Aut(\Fil\mathcal{G}_0)$ and the set of pairs of power series of the form
$$\begin{aligned}
f(e) & = \textstyle{ e + \sum_{i \geq 0}a_i\gamma^{p^i} } \cr
f(\gamma) & = \gamma + \textstyle{ \sum_{i \geq 1}b_i\gamma^{p^i} } \cr
\end{aligned}$$
with $a_i \in R_{2p^i-1}$ and $b_i \in R_{2p^i-2}$. This shows that the underlying $S_0$-stack of the $\mathbf{E}_1$-group $\Aut(\Fil\mathcal{G}_0)$ is represented by the stated Dirac $\mathbb{F}_p$-algebra $C$, so in particular, we conclude that $\Aut(\Fil\mathcal{G}_0)$ is an affine Dirac group $S_0$-scheme.

The multiplication in the Dirac group $S_0$-scheme $\Aut(\Fil\mathcal{G}_0)$ is given, on $R$-valued points, by the composition of maps
$$\xymatrix{
{ \Hom_T(Y_{0,T},Y_{0,T}) \times_T \Hom_T(Y_{0,T},Y_{0,T}) } \ar[r]^-{\circ} &
{ \Hom_T(Y_{0,T},Y_{0,T}). } \cr
}$$
So let $(f,g)$ be a pair of maps with $f$ given by
$$\begin{aligned}
f(e) & = \textstyle{ e + \sum_{i \geq 0}a_i\gamma^{p^i} } \cr
f(\gamma) & = \textstyle{ \gamma + \sum_{i \geq 1}b_i\gamma^{p^i} } \cr
\end{aligned}$$
and with $g$ given by
$$\begin{aligned}
g(e) & = \textstyle{ e + \sum_{i \geq 0}c_i\gamma^{p^i} } \cr
g(\gamma) & = \textstyle{ \gamma + \sum_{i \geq 1}d_i\gamma^{p^i}. } \cr
\end{aligned}$$
Since $\Spf$ is contravariant, we have
$$\begin{aligned}
(f \circ g)(e) & = \textstyle{ g(f(e)) = g(e + \sum_{i \geq 0}a_i\gamma^{p^i}) = g(e) + g(\sum_{i \geq 0}a_i\gamma^{p^i}) } \cr
{} & = \textstyle{ e + \sum_{i \geq 0}c_i\gamma^{p^i} + \sum_{i,j \geq 0} a_i^{p^j}d_j\gamma^{p^{i+j}} } \cr
(f \circ g)(\gamma) & = \textstyle{ g(\sum_{i \geq 0}b_i\gamma^{p^i}) = \sum_{i,j \geq 0} b_i^{p^j}d_j\gamma^{p^{i+j}} } \cr
\end{aligned}$$
with $b_0 = d_0 = 1$. This shows that the comultiplication on $C$ is as stated.
\end{proof}

We recall Milnor's calculation of the dual Steenrod algebra at an odd prime. To this end, we consider the groupoid $S \colon \Delta^{\op} \to \Aff_{/S_0}$ defined by
$$S \simeq \Spec(\pi_*(\mathbb{F}_p^{\,\otimes[-]})).$$
Since its face maps are flat, the geometric realization
$$X \simeq \lvert S \rvert \in \sheaves(\Aff)_{/S_0}$$
is a geometric Dirac $S_0$-stack. Moreover, the face maps $d_0,d_1 \colon S_1 \to S_0$ are necessarily equal, because $\mathbb{F}_p$ is a prime field, so we conclude that $X$ is the classifying stack
$$X \simeq BG$$
of the affine Dirac $S_0$-group $G$ corepresented by the dual Steenrod algebra $\smash{ A \simeq \pi_*(\mathbb{F}_p^{\,\otimes[1]}) }$. The geometric realization of the simplicial Milnor formal group
$$\mathcal{G}^M \simeq \lvert \mathcal{G}_{\mathbb{F}_p^{\,\otimes[-]}}^M \rvert \in \FGroup(X)$$
comes equipped with a canonical filtration and a canonical isomorphism of graded vector bundles
$$\xymatrix{
{ \mathbb{A}_X(e,\gamma) } \ar[r]^-{\phi^M} &
{ \gr(\Lie(\mathcal{G}^M)). } \cr
}$$
We now use Milnor's theorem to determine the structure of the pair $(\Fil(\mathcal{G}^M),\phi^M)$.

\begin{theorem}\label{thm:milnor}Let $S_0 = \Spec(\mathbb{F}_p)$ with $p$ an odd prime, and let $T$ be a Dirac $S_0$-stack. The functor that to $f \colon T \to X$ assigns the pair $(f^*\Fil(\mathcal{G}^M),f^*\phi^M)$ is an equivalence from $X(T)$ to the groupoid $X'(T)$ of pairs $(\Fil(\mathcal{G}),\phi)$ of a filtered formal group $\Fil(\mathcal{G})$ over $T$ that is locally isomorphic to a filtered additive formal group and an isomorphism $\phi \colon \mathbb{A}_X(e,\gamma) \to \gr(\Lie(\mathcal{G}))$
of graded vector bundles over $T$.
\end{theorem}

\begin{proof}By definition, we have $X' \simeq BG'$, where $G' \simeq \Aut(\Fil\mathcal{G}_0)$. Moreover, the map of Dirac $S_0$-stacks $h \colon X \to X'$ in the statement determines and is determined by a map of Dirac $S_0$-groups $\rho \colon G \to G'$. It follows from \cref{prop:milnor} that the latter map determines and is determined by the power series
$$\rho(e),\rho(\gamma) \in \pi_*(\mathbb{F}_p^{\otimes[1]})[[e,\gamma]]$$
of the form
$$\begin{aligned}
\rho(e) & = \textstyle{ e + \sum_{i \geq 0} a_i\gamma^{p^i} } \cr
\rho(\gamma) & = \textstyle{ \gamma + \sum_{i \geq 1} b_i\gamma^{p^i} } \cr
\end{aligned}$$
with $a_i \in \pi_{2p^i-1}(\mathbb{F}_p^{\otimes[1]})$ and $b_i \in \pi_{2p^i-2}(\mathbb{F}_p^{\otimes[1]})$. Now, the statement that $\rho$ is an isomorphism of Dirac $S_0$-group schemes is equivalent to the statement that all $a_i$ and $b_i$ are nonzero. But this is precisely the statement of~\cite[Lemma~6]{milnor}.
\end{proof}

\appendix

\section{Accessible presheaves}
\label{appendix:accessible_sheaves}

In the body of the paper, we have occasion to consider presheaves of anima on the $\infty$-category $\Aff$ of affine Dirac schemes. The $\infty$-category $\Aff$ is not small, but it is coaccessible in the sense that $\Aff^{\op} \simeq \CAlg(\Ab)$ is accessible. In this context, it is natural to only consider accessible presheaves. We will show that, in general, the $\infty$-category $\presheaves(\mathcal{C})$ of accessible presheaves of anima on a coaccessible $\infty$-category $\mathcal{C}$ has the following properties:
\begin{enumerate}
\item[(1)]The Yoneda embedding restricts to a fully faithful functor $h \colon \mathcal{C} \to \presheaves(\mathcal{C})$, which exhibits $\presheaves(\mathcal{C})$ as the cocompletion of $\mathcal{C}$.
\item[(2)] The $\infty$-category $\presheaves(\mathcal{C})$ satisfies the Giraud axioms, with the exception that, in general, it is not presentable.
\item[(3)] If $\mathcal{C}$ admits pullbacks, then there are adjoint functors
$$\begin{xy}
(0,0)*+{ \presheaves(\mathcal{C})_{/X} };
(28,0)*+{ \presheaves(\mathcal{C})_{/Y} };
(5.5,0)*+{}="1";
(22,0)*+{}="2";
(0,-5)*+{};
{ \ar@<.5ex>@/^.8pc/_-{f_*} "1";"2";};
{ \ar^-{f^*} "2";"1";};
{ \ar@<-1ex>@/_.8pc/_-{f_!} "1";"2";};  
\end{xy}$$
associated with every map $f \colon Y \to X$ in $\presheaves(\mathcal{C})$, where $f_!$ is given by restriction along $f$, and where $f^*$ given by base-change along $f$.
\end{enumerate}
This appendix does not go as far as one would want. We expect that there is a larger and more systematic theory of $\infty$-categories which behave as $\infty$-categories of accessible presheaves on coaccessible $\infty$-categories. However, this is not the place to develop such a theory, but to encourage work in this direction, we end with a list of properties that we expect such a theory to have; see \cref{question:properties_of_macro_categories}.

We now proceed to prove the properties~(1)--(3) of accessible presheaves listed above. We recall from~\cite[Definition~A.2.6.3]{luriehtt} that the relation $\lambda \ll \kappa$ among regular cardinals indicates that for every $\kappa_0 < \kappa$ and $\lambda_0 < \lambda$, we have $\smash{ \kappa_0^{\lambda_0} < \kappa }$.

\begin{lemma}
\label{lemma:filteredness_of_overcategories_of_accessible_cat}
Let $\kappa \gg \lambda$ be regular cardinals. If $\mathcal{D}$ is a $\lambda$-accessible $\infty$-category, then the slice category $(\mathcal{D}^{\kappa})_{/d}$ is $\kappa$-filtered for all $d \in \mathcal{D}$.
\end{lemma}

\begin{proof}
Given a diagram $p \colon K \to (\mathcal{D}^{\kappa})_{/d}$ with the property that $K$ is $\kappa$-small, we must produce a diagram $\bar{p} \colon K^{\triangleright} \to (\mathcal{D}^{\kappa})_{/d}$ and an equivalence $\bar{p}|_K \simeq p$. We identify 
$$\mathcal{D} \simeq \Ind_{\lambda}(\mathcal{D}^{\lambda}) \subset \presheaves(\mathcal{D}^{\lambda})$$
with the full subcategory of presheaves of anima on $\mathcal{D}^{\lambda}$ generated under $\lambda$-filtered colimits by the full subcategory of representable presheaves. In these terms, we can identify $p$ with a pair $(q,g)$ of a diagram $q \colon K \to \mathcal{D}^{\kappa} \subset \mathcal{P}(\mathcal{D}^{\lambda})$ and a map
$$\xymatrix{
{ \varinjlim_K q } \ar[r]^-{g} &
{ d } \cr
}$$
from the colimit calculated in $\presheaves(\mathcal{D}^{\lambda})$, and we must show that $g$ factors through an object of $\mathcal{D}^{\kappa}$. We recall from~\cite[Proposition~5.4.2.11]{luriehtt} that, since $\kappa \gg \lambda$, every $\lambda$-accessible $\infty$-category is $\kappa$-accessible. So we can write $d \in \mathcal{D}$ as a $\kappa$-filtered colimit
$$\textstyle{ d \simeq \varinjlim_{\alpha}d_{\alpha} }$$
with $d_{\alpha} \in \mathcal{D}^{\kappa} \subset \mathcal{D}$, and the proof of loc.cit.\ shows that we can arrange that each of the $d_{\alpha}$ be a colimit of a $\kappa$-small $\lambda$-filtered diagram in $\mathcal{D}^{\lambda} \subset \mathcal{D}$. Accordingly, each of the $d_{\alpha}$ is also $\kappa$-compact as an object of $\presheaves(\mathcal{D}^{\lambda})$. Now, the inclusion $\mathcal{D} \subset \presheaves(\mathcal{D}^{\lambda})$ preserves $\lambda$-filtered colimits, so, in particular, it preserves $\kappa$-filtered colimits. Therefore, we also have
$$\textstyle{ d \simeq \varinjlim_{\alpha}d_{\alpha} }$$
with the colimit calculated in $\presheaves(\mathcal{D}^{\lambda})$. Since the domain of the map $g$ in question is a $\kappa$-small colimit of $\kappa$-compact objects in $\presheaves(\mathcal{D}^{\lambda})$, it is itself $\kappa$-compact, and therefore, we conclude that $g$ factors through some $d_{\alpha}$, as we wanted to prove.
\end{proof}

\begin{notation}
If $\kappa$ is a regular cardinal and if $\mathcal{C}$ is a $\kappa$-coaccessible $\infty$-category, then we write
$\mathcal{C}_{\kappa} \subset \mathcal{C}$ for the full subcategory spanned by the $\kappa$-cocompact objects. The $\infty$-category $\mathcal{C}_{\kappa}$ is essentially small, and $\mathcal{C} \simeq \Pro_{\kappa}(\mathcal{C}_{\kappa})$ is its completion under small $\kappa$-cofiltered limits.
\end{notation}

\begin{proposition}
\label{prop:equivalent_conditions_for_accessability}Let $\mathcal{C}$ be a coaccessible $\infty$-category. For a functor $X \colon \mathcal{C}^{\op} \to \spaces$, the following conditions are equivalent:
\begin{enumerate}
\item[{\rm (1)}]The functor $X \colon \mathcal{C}^{\op} \to \spaces$ is accessible.
\item[{\rm (2)}]The functor $X \colon \mathcal{C}^{\op} \to \spaces$ is the left Kan extension of a functor $Y \colon (\mathcal{C}_{\kappa})^{\op} \to \spaces$ along the canonical inclusion $i \colon (\mathcal{C}_{\kappa})^{\op} \to \mathcal{C}^{\op}$ for some regular cardinal $\kappa$.
\item[{\rm (3)}]The functor $X \colon \mathcal{C}^{\op} \to \spaces$ is a colimit in $\Fun(\mathcal{C}^{\op},\spaces)$ of a small diagram of representable functors.
\end{enumerate}
\end{proposition}

\begin{proof}To prove that~(1) implies~(2), we first choose a regular carcinal $\lambda$ such $\mathcal{C}$ is $\lambda$-coaccessible and such that $X \colon \mathcal{C}^{\op} \to \spaces$ is $\lambda$-accessible. We next choose a regular cardinal $\kappa$ such that $\kappa \gg \lambda$. We claim that the counit map
$$\xymatrix{
{ i_!i^*(X) } \ar[r] &
{ X } \cr
}$$
is an equivalence. The value of $i_!i^*(X)$ at $c \in \mathcal{C}^{\op}$ is the colimit of the diagram
$$\xymatrix{
{ ((\mathcal{C}_{\kappa})^{\op})_{/c} } \ar[r]^-{p} &
{ (\mathcal{C}_{\kappa})^{\op} } \ar[r]^-{i} &
{ \mathcal{C}^{\op} } \ar[r]^-{X} &
{ \spaces. } \cr
}$$
Now, by  \cref{lemma:filteredness_of_overcategories_of_accessible_cat}, the slice category on the left is $\kappa$-filtered, and since $X$ preserves $\lambda$-filtered, and hence, $\kappa$-filtered colimits, we conclude that the counit map
$$\xymatrix{
{ i_!i^*(X)(c) } \ar[r] &
{ X(c) } \cr
}$$
is an equivalence, which shows that~(2) holds with $Y \simeq i^*(X)$.

To prove that~(2) implies~(3), we assume that $X \simeq i_!(Y)$ with $Y \in \Fun(\mathcal{C}_{\kappa}^{\op},\spaces)$. Since $\mathcal{C}_{\kappa}$ is essentially small, we can identify $Y$ with a small colimit of representable functors. But left Kan extension preserves colimits and takes representable functors to representable functors, so we conclude that~(3) holds.

Finally, to see that~(3) implies~(1), we note that, since $\mathcal{C}^{\op}$ is accessible, every representable functor $X \colon \mathcal{C}^{\op} \to \spaces$ is accessible. Hence, given a small diagram of representable functors $X_{\alpha}$, we can find a small cardinal $\lambda$ such that each of the $X_{\alpha}$ preserves $\lambda$-filtered colimits. Since colimits commute with colimits, their levelwise colimit again preserves $\lambda$-filtered colimits, and hence, is accessible. 
\end{proof}

\begin{definition}
\label{definition:accessible_presheaves}
Let $\mathcal{C}$ be a coaccessible $\infty$-category.
The $\infty$-category of accessible presheaves of anima on $\mathcal{C}$ is the full subcategory
$$\presheaves(\mathcal{C}) \simeq \Fun_{\mathcal{A}}(\mathcal{C}^{\op},\spaces) \subset \Fun(\mathcal{C}^{\op},\spaces)$$
spanned by the functors that satisfy the equivalent conditions of \cref{prop:equivalent_conditions_for_accessability}.
\end{definition}

\begin{remark}[Universal property of accessible presheaves] 
\label{remark:universal_property_of_accessible_presheaves_as_cocompletion}
If $\mathcal{C}$ is a coaccessible $\infty$-category, then the Yoneda embedding factors through a functor
$$\xymatrix{
{ \mathcal{C} } \ar[r]^-{h} &
{ \presheaves(\mathcal{C}) } \cr
}$$
which, by abuse of language, we also refer to as the Yoneda embedding. It follows from \cref{prop:equivalent_conditions_for_accessability} that this functor is characterized by the following universal property: If $\mathcal{D}$ is any $\infty$-category, which admits small colimits, then composition with the Yoneda embedding gives rise to an equivalence
$$\xymatrix{
{ \Fun^L(\presheaves(\mathcal{C}),\mathcal{D}) } \ar[r] &
{ \Fun(\mathcal{C},\mathcal{D}) } \cr
}$$
from the full subcategory $\Fun^L(\presheaves(\mathcal{C}),\mathcal{D}) \subset \Fun(\presheaves(\mathcal{C}),\mathcal{D})$ spanned by the functors that preserve small colimits. An inverse takes $F \colon \mathcal{C} \to \mathcal{D}$ to its left Kan extension $h_!(F) \colon \presheaves(\mathcal{C}) \to \mathcal{D}$, which exists because every accessible presheaf is a small colimit of representables and $\mathcal{D}$ is cocomplete; see also \cite[Remark~5.3.5.9]{luriehtt}.
\end{remark}

\begin{warning}
\label{warning:abusivenotation}
Let $\mathcal{C}$ be a coaccessible $\infty$-category. We write $h \colon \mathcal{C} \to \presheaves(\mathcal{C})$ for its cocompletion. If $\mathcal{C}$ is small, then \cite[Corollary~5.4.3.6]{luriehtt} shows that its cocompletion agrees with the Yoneda embedding $h \colon \mathcal{C} \to \Fun(\mathcal{C}^{\op},\spaces)$ to the functor $\infty$-category. In general, however, the functor $\infty$-category is larger and pathological. For example, we show in \cref{lemma:infinity_category_of_accessible_presheaves_is_locally_small} that $\presheaves(\mathcal{C})$ is locally small, but this is generally not true for the functor $\infty$-category.
\end{warning}

\begin{remark}
\label{remark:idempotent_completion}
We recall from \cite[Corollary~5.4.3.6]{luriehtt} that a small $\infty$-category $\mathcal{C}$ is accessible if and only if it is idempotent complete. However, if $j \colon \mathcal{C} \to \mathcal{C}'$ is an idempotent completion of $\mathcal{C}$, then $\mathcal{C}'$ is again small and $j_! \colon \presheaves(\mathcal{C}) \to \presheaves(\mathcal{C}')$ is an equivalence. Thus, the theory of accessible presheaves on coaccessible $\infty$-categories is a strict generalization of the theory of presheaves on small $\infty$-categories.
\end{remark}

\begin{lemma}
\label{lemma:infinity_category_of_accessible_presheaves_is_locally_small}
Let $\mathcal{C}$ be coaccessible $\infty$-category. The $\infty$-category $\presheaves(\mathcal{C})$ of accessible presheaves of anima on $\mathcal{C}$ is locally small. 
\end{lemma}

\begin{proof}
Given $X,Y \in \presheaves(\mathcal{C})$, the mapping anima $\Map(Y,X)$ is potentially large. We must show that it is small. By \cref{prop:equivalent_conditions_for_accessability}, we may write $Y \simeq \varinjlim_{\alpha} h(c_{\alpha})$ as a small colimit of representable functors. Thus, the anima
$$\textstyle{ \Map(Y,X) \simeq \Map(\varinjlim_{\alpha} h(c_{\alpha}),X) \simeq \varprojlim_{\alpha}X(c_{\alpha}) }$$
is a small limit of small anima, and hence, it is indeed small. 
\end{proof}

\begin{lemma}
\label{lemma:left_kan_extension_along_cardinality_change_preserves_certain_limits}
Suppose that $\mathcal{C}$ is $\lambda$-coaccessible and let $\lambda \ll \kappa \leq \kappa'$ be infinite regular cardinals. In this situation, the left Kan extension
$$\xymatrix@C=8mm{
{ \presheaves(\mathcal{C}_{\kappa}) } \ar[r]^-{i_!} &
{ \presheaves(\mathcal{C}_{\kappa'}) } \cr
}$$
along the inclusion $i \colon (\mathcal{C}_{\kappa})^{\op} \to (\mathcal{C}_{\kappa'})^{\op}$ preserves all colimits and $\kappa$-small limits.
\end{lemma}

\begin{proof}
Since left Kan extension is a left adjoint, it preserves colimits. For limits, we note that if $X \in \presheaves(\mathcal{C}_{\kappa})$ and $c \in (\mathcal{C}_{\kappa'})^{\op}$, then $i_!(X)(c)$ is the colimit of the diagram
$$\xymatrix{
{ ((\mathcal{C}_{\kappa})^{\op})_{/c} } \ar[r]^-{p} &
{ (\mathcal{C}_{\kappa})^{\op} } \ar[r]^-{X} &
{ \spaces, } 
}$$
where $p$ is the projection. Thus, the statement follows from \cref{lemma:filteredness_of_overcategories_of_accessible_cat} and from the fact proved in~\cite[Proposition~5.3.3.3]{luriehtt} that $\kappa$-filtered colimits and $\kappa$-small limits of anima commute.
\end{proof}

\begin{proposition}
\label{proposition:accessible_presheaves_are_complete_cocomplete}
Let $\mathcal{C}$ be a coaccessible $\infty$-category. 
The $\infty$-category $\presheaves(\mathcal{C})$ of accessible presheaves of anima on $\mathcal{C}$ admits all small limits and colimits, and both are calculated pointwise. 
\end{proposition}

\begin{proof}
We claim that $\presheaves(\mathcal{C}) \subset \Fun(\mathcal{C}^{\op},\spaces)$ is closed under small limits and colimits. This is clear from condition~(3) of \cref{prop:equivalent_conditions_for_accessability} in the case of colimits. To prove the claim for limits, we fix a diagram $p \colon K \to \presheaves(\mathcal{C})$ with $K$ small. We can choose a cardinal $\kappa$ such that 
\begin{enumerate}
\item $\mathcal{C}$ is $\kappa$-coaccessible, 
\item for every $k \in K$, $p(k) \colon \mathcal{C}^{op} \to \spaces$ preserves $\kappa$-filtered colimits, and 
\item $K$ is $\kappa$-small. 
\end{enumerate}
With this choice, the limit of $p$ calculated in $\Fun(\mathcal{C}^{\op},\spaces)$ preserves $\kappa$-filtered colimits, and therefore, condition~(1) of \cref{prop:equivalent_conditions_for_accessability} shows that it belongs to $\presheaves(\mathcal{C})$.
\end{proof}

\begin{theorem}
\label{theorem:giraud_axioms_for_accessible_presheaves}
The $\infty$-category $\presheaves(\mathcal{C})$ of accessible presheaves on a coaccessible $\infty$-category $\mathcal{C}$ enjoys the following properties: 
\begin{enumerate}[leftmargin=8mm]
\item[{\rm (i)}] It is cocomplete, complete, and generated under small colimits by the essential image of $h \colon \mathcal{C} \to \presheaves(\mathcal{C})$, which consists of compact objects and is coaccessible.
\item[{\rm (ii)}]Colimits in $\presheaves(\mathcal{C})$ are universal.
\item[{\rm (iii)}]Coproducts in $\presheaves(\mathcal{C})$ are disjoint.
\item[{\rm (iv)}]Every groupoid in $\presheaves(\mathcal{C})$ is
      effective.
\end{enumerate}
\end{theorem}

\begin{proof}
Since colimits in $\presheaves(\mathcal{C})$ are pointwise, we conclude that all representable presheaves are compact. Moreover, \cref{prop:equivalent_conditions_for_accessability} shows that $\presheaves(\mathcal{C})$ is generated under small colimits by the representable presheaves, and \cref{proposition:accessible_presheaves_are_complete_cocomplete} shows that $\presheaves(\mathcal{C})$ is complete and cocomplete. This proves~(i). The statements (ii)--(iv) concern the interaction of limits and colimits in $\presheaves(\mathcal{C})$, and since both are calculated pointwise, these statements follow from the corresponding statements about the $\infty$-category $\spaces$ of anima, which are proved in~\cite[Theorem~6.1.0.6]{luriehtt}.
\end{proof}

We proceed to discuss localization. We recall that if $\mathcal{C}$ is small, then for every presheaf $X \colon \mathcal{C}^{\op} \to \spaces$, there is a canonical equivalence
$$\xymatrix{
{ \presheaves(\mathcal{C}_{/X}) } \ar[r] &
{ \presheaves(\mathcal{C})_{/X} } \cr
}$$
with $j \colon \mathcal{C}_{/X} \to \mathcal{C}$ the base-change of $p \colon \presheaves(\mathcal{C})_{/X} \to \presheaves(\mathcal{C})$ along $h \colon \mathcal{C} \to \mathcal{P}(\mathcal{C})$. There is an essentially unique diagram
$$\xymatrix{
{ \mathcal{C}_{/X} } \ar[r]^-{j} \ar[d]^-{h} &
{ \mathcal{C} } \ar[d]^-{h} \cr
{ \presheaves(\mathcal{C}_{/X}) } \ar[r]^-{j_!} &
{ \presheaves(\mathcal{C}) } \cr
}$$
with $j_!$ cocontinuous, and $j_!$ factors through an equivalence as stated. We wish to show that there is an analogous picture in the case, where $\mathcal{C}$ is coaccessible.

\begin{lemma}
\label{lemma:accessible_undercategory_of_accessible_is_also_accessible}
Suppose that $\mathcal{C}$ is a coaccessible $\infty$-category.  Let $X \in \presheaves(\mathcal{C})$, and let $j \colon \mathcal{C}_{/X} \to \mathcal{C}$ be the base-change of $p \colon \presheaves(\mathcal{C})_{/X} \to \presheaves(\mathcal{C})$ along $h \colon \mathcal{C} \to \presheaves(\mathcal{C})$. In this situation, the slice $\infty$-category $\mathcal{C}_{/X}$ and the functor $j \colon \mathcal{C}_{/X} \to \mathcal{C}$ are coaccessible.
\end{lemma}

\begin{proof}
The functor $j^{\op} \colon (\mathcal{C}_{/X})^{\op} \to \mathcal{C}^{\op}$ is the left fibration classified by the functor $X \colon \mathcal{C}^{\op} \to \spaces$. So there is a cartesian diagram of $\infty$-categories
$$\xymatrix{
{ (\mathcal{C}_{/X})^{\op} } \ar[r] \ar[d]^-{j^{\op}} &
{ \spaces_* } \ar[d]^-{p} \cr
{ \mathcal{C}^{\op} } \ar[r]^-{\phantom{,}X\phantom{,}} &
{ \spaces } \cr
}$$
with the right-hand vertical map given by the universal left fibration, which is the canonical projection $p \colon \spaces_* \simeq \spaces_{/1} \to \spaces$. Since $p$ and $X$ both are accessible functors between accessible $\infty$-categories, we conclude from \cite[Proposition~5.4.6.6]{luriehtt} that the same is true for $j^{\op}$.
\end{proof}

\begin{remark}
\label{remark:cofiltered_limits_in_overcategories_are_computed_in_underlying_category}
In the context of \cref{lemma:accessible_undercategory_of_accessible_is_also_accessible}, $j \colon \mathcal{C}_{/X} \to \mathcal{C}$ is a right fibration. Since every $\kappa$-cofiltered $\infty$-category is weakly contractible, it follows from \cite[\href{https://kerodon.net/tag/02KS}{Tag 02KS}]{kerodon} that if $\mathcal{C}$ is $\kappa$-coaccessible, then $\mathcal{C}_{/X}$ admits $\kappa$-cofiltered limits and $j$ preserves and reflects them. This does not necessarily mean that $\mathcal{C}_{/X}$ is $\kappa$-coaccessible, since it may not be generated under $\kappa$-cofiltered limits by $\kappa$-cocompact objects. 
\end{remark}

\begin{proposition}
\label{proposition:accessible_presheaves_over_an_overcategory_is_an_overcategory}
Let $\mathcal{C}$ be a coaccessible $\infty$-category, and $X \in \presheaves(\mathcal{C})$. In this situation, the Yoneda embedding induces an equivalence
$$\xymatrix{
{ \presheaves(\mathcal{C}_{/X}) } \ar[r] &
{ \presheaves(\mathcal{C})_{/X}. } \cr
}$$
\end{proposition}

\begin{proof}Since $\mathcal{C}_{/X}$ is coaccessible by \cref{lemma:accessible_undercategory_of_accessible_is_also_accessible}, and since $\presheaves(\mathcal{C})_{/X}$ admits small colimits by \cref{proposition:accessible_presheaves_are_complete_cocomplete}, the map $h_{/X} \colon \mathcal{C}_{/X} \to \presheaves(\mathcal{C})_{/X}$ induced by the Yoneda embedding extends essentially uniquely to the map in the statement. We may view this map as the restriction of the corresponding map
$$\xymatrix{
{ \widehat{\presheaves}(\mathcal{C}_{/X}) } \ar[r] &
{ \widehat{\presheaves}(\mathcal{C})_{/X} } \cr
}$$
between the $\infty$-categories of presheaves of large anima. The latter is an equivalence, because $\mathcal{C}$ is small with respect to this larger universe, so we conclude that the map in the statement is fully faithful. It remains to show that it is essentially surjective. Given $(Y,f \colon Y \to X)$ in $\presheaves(\mathcal{C})_{/X}$, we can find a diagram
$$\xymatrix{
{ K } \ar[r]^-{p} \ar[d] &
{ \mathcal{C} } \ar[d]^-{h} \cr
{ K^{\triangleright } } \ar[r]^-{\bar{p}} &
{ \presheaves(\mathcal{C}) } \cr
}$$
with $K$ small such that $\bar{p}$ is a colimit diagram, whose value at the cone point is $Y$. Thus, the maps $\bar{p}$ and $f$ determine a diagram $q \colon K \to \mathcal{C}_{/X}$, and the map in the statement takes the colimit of the composite map
$$\xymatrix{
{ K } \ar[r]^-{q} &
{ \mathcal{C}_{/X} } \ar[r]^-{h} &
{ \presheaves(\mathcal{C}_{/X}) } \cr
}$$
to the given object $(Y,f \colon Y \to X)$, as desired.
\end{proof}

We recall from~\cite[Lemma~6.1.1.1]{luriehtt} that a map $f \colon Y \to X$ in an $\infty$-category $\mathcal{D}$ that admits pullbacks gives rise to an adjoint pair of functors
$$\xymatrix{
{ \mathcal{D}_{/Y} } \ar@<.7ex>[r]^-{f_!} &
{ \mathcal{D}_{/X} } \ar@<.7ex>[l]^-{f^*} \cr
}$$
and from~\cite[Lemma~4.4.2.1]{luriehtt} that given a cartesian diagram
$$\xymatrix@C=10mm{
{ Y' } \ar[r]^-{g'} \ar[d]^-{f'} &
{ Y } \ar[d]^-{f} \cr
{ X' } \ar[r]^-{g} &
{ X } \cr
}$$
in $\mathcal{D}$, the composite ``base-change'' map
$$\xymatrix{
{ f'_!g'{}^* } \ar[r] &
{ g^*g_!f_!'g'{}^* \simeq g^*f_!g_!'g'{}^* } \ar[r] &
{ g^*f_! } \cr
}$$
is an equivalence. If $\mathcal{D} \simeq \presheaves(\mathcal{C})$ with $\mathcal{C}$ coaccessible, then \cref{theorem:giraud_axioms_for_accessible_presheaves} shows that $f^*$ preserves small colimits. However, since $\presheaves(\mathcal{C})$ may not be presentable, we cannot invoke the adjoint functor theorem to conclude that $f^*$ admits a right adjoint $f_*$. We will prove that this is nevertheless the case, under the assumption that $\mathcal{C}$ admits pullbacks. We do not know if this assumption is necessary.

\begin{theorem}
\label{theorem:existence_of_additional_right_adjoint}
Let $\mathcal{C}$ be a coaccessible $\infty$-category, and let $f \colon Y \to X$ be a map in the $\infty$-category $\presheaves(\mathcal{C})$ of accessible presheaves on $\mathcal{C}$. If $\mathcal{C}$ admits pullbacks, then the functor $f^* \colon \presheaves(\mathcal{C})_{/X} \to \presheaves(\mathcal{C})_{/Y}$ admits a right adjoint $f_* \colon \presheaves(\mathcal{C})_{/Y} \to \presheaves(\mathcal{C})_{/X}$.
\end{theorem}

Before we begin the the proof of \cref{theorem:existence_of_additional_right_adjoint}, which will occupy the remainder of this appendix, we give an example to demonstrate the subtleties of situation.

\begin{example}
\label{example:existence_of_additional_right_adjoint}Under the equivalence of  \cref{proposition:accessible_presheaves_over_an_overcategory_is_an_overcategory}, we can identify the functor $f^*$ in \cref{theorem:existence_of_additional_right_adjoint} with the functor $f^* \colon \presheaves(\mathcal{C}_{/X}) \to \presheaves(\mathcal{C}_{/Y})$ given by  restriction along $f_! \colon \mathcal{C}_{/Y} \to \mathcal{C}_{/X}$. Thus, one might expect that for every coaccessible functor $f \colon \mathcal{D} \to \mathcal{C}$ between coaccessible $\infty$-categories, the functor
$$\xymatrix@C=10mm{
{ \presheaves(\mathcal{C}) } \ar[r]^{f^*} &
{ \presheaves(\mathcal{D}) } \cr
}$$
given by restriction along $f$ admit a right adjoint. Surprisingly, this is not the case in general! For an explicit counterexample, we consider the functor
$$\xymatrix@C=10mm{
{ 1 } \ar[r]^-{f} &
{ \spaces^{\op} } \cr
}$$
that to the unique object assigns the set $\{0,1\}$. In this case, the functor
$$\xymatrix{
{ \presheaves(\spaces^{\op}) } \ar[r]^-{f^*} &
{ \spaces } \cr
}$$
takes an accessible functor $X \colon \spaces \to \spaces$ to its value at $\{0,1\}$. We claim that this functor $f^*$ does not admit a right adjoint $f_*$. If it did, then we would have
$$f_*(\{0,1\})(-) \simeq \Map(\Map(-,\{0,1\}),\{0,1\}),$$
but this functor is not accessible, since any accessible functor $X \colon \spaces \to \spaces$ preserves $\kappa$-compact objects for some $\kappa$.
\end{example}

Let $\mathcal{C}$ be a coaccessible $\infty$-category, and let $f \colon Y \to X$ be a map of accessible presheaves on $\mathcal{C}$. We consider the composite functor
$$\xymatrix@C=10mm{
{ \mathcal{C}_{/X} } \ar[r]^-{h} &
{ \presheaves(\mathcal{C})_{/X} } \ar[r]^-{f^*} &
{ \presheaves(\mathcal{C})_{/Y} } \ar[r]^-{q_!} &
{ \presheaves(\mathcal{C}) } \cr
}$$
and compose it with the functor represented by $Z \in \presheaves(\mathcal{C})$ to obtain a functor
$$\xymatrix@C=12mm{
{ (\mathcal{C}_{/X})^{\op} } \ar[r]^-{F_{f,Z}} &
{ \spaces, } \cr
}$$
which takes values in the $\infty$-category of small anima, because $\presheaves(\mathcal{C})$ is locally small, as we proved in \cref{lemma:infinity_category_of_accessible_presheaves_is_locally_small}.

\begin{lemma}
\label{lemma:accessability_bound_on_pullback_along_a_map_into_representable}
Let $\kappa$ be a regular cardinal, and let $\mathcal{C}$ be a $\kappa$-coaccessible $\infty$-category which admits pullbacks. Let $f \colon Y \to X$ be a map of accessible presheaves on $\mathcal{C}$ with $X$ a representable presheaf and with $Y$ a $\kappa$-small colimit of representable presheaves, and let $Z$ be a $\kappa$-accessible presheaf on $\mathcal{C}$. In this situation, the functor
$$\xymatrix@C=12mm{
{ (\mathcal{C}_{/X})^{\op} } \ar[r]^-{F_{f,Z}} &
{ \spaces } \cr
}$$
preserves $\kappa$-filtered colimits.
\end{lemma}

\begin{proof}
We note that by \cref{remark:cofiltered_limits_in_overcategories_are_computed_in_underlying_category}, the slice $\infty$-category $\mathcal{C}_{/X}$ admits $\kappa$-cofiltered limits.
If we write $\smash{ Y \simeq \varinjlim_{\alpha}Y_{\alpha} }$ as a $\kappa$-small colimit of representable presheaves and let $f_{\alpha} \colon h(b_{\alpha}) \to X$ be the restriction of $f$ along $Y_{\alpha} \to Y$, then
$$\begin{aligned}
{ F_{f,Z}(c,h(c) \to X) } & \simeq \Map(Y \times_Xh(c),Z) \simeq \textstyle{ \Map(\varinjlim_{\alpha} Y_{\alpha} \times_Xh(c),Z) } \cr
{} & \simeq \textstyle{ \varprojlim_{\alpha} \Map(Y_{\alpha} \times_X h(c),Z) \simeq \varprojlim_{\alpha} F_{f_{\alpha},Z}(c,h(c) \to X), } \cr
\end{aligned}$$
and since $\kappa$-small limits and $\kappa$-filtered colimits of anima commute, we may assume that $Y$ is representable.

So we write $X \simeq h(x)$ and $Y \simeq h(y)$ and suppose that $(c,h(c) \to X)$ is the limit of a $\kappa$-cofiltered diagram $(c_{\alpha},h(c_{\alpha}) \to X)$ in $\mathcal{C}_{/X}$. We conclude from
\cref{remark:cofiltered_limits_in_overcategories_are_computed_in_underlying_category} that $c$ is the limit of the underlying diagram $(c_{\alpha})$ in $\mathcal{C}$. Since the Yoneda embedding preserves all limits that exist in $\mathcal{C}$, we have
$$\textstyle{ Y \times_Xh(c) \simeq h(y) \times_{h(x)}h(c) \simeq h(y \times_xc) \simeq h(y \times_x\varprojlim_{\alpha}c_{\alpha}) \simeq h(\varprojlim_{\alpha} y \times_xc_{\alpha}). }$$
Therefore, since $Z$ is $\kappa$-accessible, we have
$$\begin{aligned}
{ F_{f,Z}(c,h(c) \to X) } & \simeq \Map(Y \times_Xh(c),Z) \simeq  \textstyle{ Z(\varprojlim_{\alpha} y\times_xc_{\alpha}) } \cr {} & \simeq \textstyle{ \varinjlim_{\alpha} Z(y\times_xc_{\alpha}) \simeq \varinjlim_{\alpha} F_{f,Z}(c_{\alpha},h(c_{\alpha}) \to X) } \cr
\end{aligned}$$
as we wanted to prove.
\end{proof}

We will need the following notion.

\begin{definition}
\label{definition:relatively_kappa_compact_map_of_accessible_presheaves}
Let $\mathcal{C}$ be a coaccessible $\infty$-category, and let $\kappa$ be a regular cardinal. A map $f \colon Y \to X$ of accessible presheaves on $\mathcal{C}$ is $\kappa$-compact if for every map $g \colon h(c) \to X$ from a representable presheaf on $\mathcal{C}$, the domain $Y'$ of the base-change $f' \colon Y' \to h(c)$ of $f$ along $g$ can be written as a $\kappa$-small colimit of representable presheaves on $\mathcal{C}$.
\end{definition}

\begin{lemma}
\label{lemma:upper_bound_on_compactness_of_pullbacks}
If $\mathcal{C}$ is a coaccessible $\infty$-category which admits pullbacks, then every map $f \colon Y \to X$ of accessible presheaves on $\mathcal{C}$ is $\kappa$-compact for some $\kappa$.
\end{lemma}

\begin{proof}
Suppose first that $X \simeq h(x)$ is representable. We can write $Y \simeq \varinjlim_{\alpha} h(y_{\alpha})$ as a $\kappa$-small colimit of representable presheaves for some regular cardinal $\kappa$. Now, for any map $g \colon h(c) \to X$ from a representable presheaf, we find that
$$\textstyle{ Y \times_Xh(c) \simeq \varinjlim_{\alpha} h(y_{\alpha}) \times_{h(x)}h(c) \simeq \varinjlim_{\alpha} h(y_{\alpha}\times_xc), }$$
which is a $\kappa$-small colimit of representable presheaves on $\mathcal{C}$.

In general, we write $\smash{ X \simeq \varinjlim_{\alpha}X_{\alpha} }$ as a small colimit of representable presheaves, and let $f_{\alpha} \colon Y_{\alpha} \to X_{\alpha}$ be the base-change of $f$ along $X_{\alpha} \to X$. Since colimits in $\presheaves(\mathcal{C})$ are calculated pointwise, any map $g \colon h(c) \to X$ from a representable presheaf will factor through $X_{\alpha} \to X$ for some $\alpha$. Moreover, by the special case that we considered above, we find that
$$Y \times_Xh(c) \simeq Y_{\alpha} \times_{X_{\alpha}}h(c)$$
is $\kappa_{\alpha}$-compact for some regular cardinal $\kappa_{\alpha}$. Thus, if we choose a regular cardinal $\kappa$ such that $\kappa \geq \kappa_{\alpha}$ for all $\alpha$, then $f$ is $\kappa$-compact.
\end{proof}

We now improve \cref{lemma:accessability_bound_on_pullback_along_a_map_into_representable} to the general case. 

\begin{lemma}
\label{lemma:absolute_right_adjoint_to_pullback_exists}
Let $\mathcal{C}$ be a coaccessible $\infty$-category which admits pullbacks, let $f \colon Y \to X$ be a map of accessible presheaves on $\mathcal{C}$, and let $Z$ be an accessible presheaf on $\mathcal{C}$. In this situation, the functor
$$\xymatrix@C=12mm{
{ (\mathcal{C}_{/X})^{\op} } \ar[r]^-{F_{f,Z}} &
{ \spaces } \cr
}$$
is accessible.
\end{lemma}

\begin{proof}
Using \cref{lemma:upper_bound_on_compactness_of_pullbacks}, we can choose a regular cardinal $\kappa$ such that $\mathcal{C}$ is $\kappa$-coaccessible, $f$ is $\kappa$-compact, and $Z$ is $\kappa$-accessible. We proceed to show that for this $\kappa$, the functor $F_{j, Z}$ preserves $\kappa$-filtered colimits. Here we use \cref{remark:cofiltered_limits_in_overcategories_are_computed_in_underlying_category} to conclude that $\mathcal{C}_{/X}$ admits $\kappa$-cofiltered limits.

So we fix a diagram $p \colon K \to \mathcal{C}_{/X}$, where $K$ is $\kappa$-cofiltered, and we wish to prove that the canonical map of anima
$$\xymatrix{
{ \varinjlim_{K^{\op}} (F_{f,Z} \circ p^{\op}) } \ar[r] &
{ F_{f,Z}(\varprojlim_K p) } \cr
}$$
is an equivalence. For any $k \in K$, the projection $j \colon K_{/k} \to K$ is a $\varprojlim$-equivalence, and therefore, by replacing $p$ by $pj$, if necessary, we may assume that $K$ has a final object. So let $1 \in K$ be final, and let $(d,g \colon h(d) \to X)$ be its value under $p$. In this situation, the functor $p \colon K \to \mathcal{C}_{/X}$ factors canonically as
$$\xymatrix@C=6mm{
{} &
{ K } \ar[dl]_-{p'} \ar[dr]^-{p} &
{} \cr
{ \mathcal{C}_{/d} } \ar[rr]^-{g_!} &
{} &
{ \mathcal{C}_{/X} } \cr
}$$
and, by base-change, we have a diagram
$$\xymatrix@C=4mm{
{ (\mathcal{C}_{/d})^{\op} } \ar[rr]^-{(g_!)^{\op}} \ar[dr]_-(.4){F_{f',Z}} &&
{ (\mathcal{C}_{/X})^{\op} } \ar[dl]^-(.4){F_{f,Z}} \cr
{} &
{ \spaces } &
{} \cr
}$$
with $f' \colon Y' \to h(c)$ the base-change of $f$ along $g$. Since $g_!$ preserves $\kappa$-cofiltered limits by \cref{remark:cofiltered_limits_in_overcategories_are_computed_in_underlying_category}, it suffices to show that the canonical map
$$\xymatrix{
{ \varinjlim_{K^{\op}} (F_{f',Z} \circ p'{}^{\op}) } \ar[r] &
{ F_{f',Z}(\varprojlim_K p') } \cr
}$$
is an equivalence. But $Y'$ is a $\kappa$-small colimit of representable presheaves, because $f$ is $\kappa$-compact, so this follows from \cref{lemma:accessability_bound_on_pullback_along_a_map_into_representable}.
\end{proof}

We now give the proof of the main result of this section. 

\begin{proof}[{Proof of \cref{theorem:existence_of_additional_right_adjoint}}]By \cite[\href{https://kerodon.net/tag/02FV}{Tag 02FV}]{kerodon}, it will suffice to show that for every object $g \colon Z \to Y$ in $\presheaves(\mathcal{C})_{/Y}$, the functor
$$\xymatrix@C=10mm{
{ (\presheaves(\mathcal{C})_{/X})^{\op} } \ar[r]^-{R_g} &
{ \spaces } \cr
}$$
given by the composition of $f^*$ and the functor represented by $g$ is representable. As we already remarked, the functor $f^*$ preserves colimits, because colimits in $\presheaves(\mathcal{C})$ are universal, so $R_g$ preserves limits. Hence, since the canonical map
$$\xymatrix{
{ \presheaves(\mathcal{C}_{/X}) } \ar[r] &
{ \presheaves(\mathcal{C})_{/X} } \cr
}$$
is an equivalence by \cref{proposition:accessible_presheaves_over_an_overcategory_is_an_overcategory}, it suffices to show that the restriction
$$\xymatrix@C=9mm{
{ (\mathcal{C}_{/X})^{\op} } \ar[r]^-{r_g} &
{ \spaces } \cr
}$$
of $R_g$ along the Yoneda embedding is an accessible functor. Indeed, if so, then the image of $r_g$ by the equivalence above is the object $f_*(g)$ that represents $R_g$.

Let $q \colon Y \to 1$ be the unique map to the final object. The adjunction
$$\xymatrix{
{ \presheaves(\mathcal{C})_{/Y} } \ar@<.7ex>[r]^-{q_!} &
{ \presheaves(\mathcal{C}) } \ar@<.7ex>[l]^-{q^*} \cr
}$$
is comonadic, and hence, the cobar construction of $g \colon Z \to Y$,
which is constructed in \cite[Example~4.7.2.7]{lurieha}, is a limit diagram
$$\xymatrix@C=10mm{
{ \Delta_+ } \ar[r]^-{C_g} &
{ \presheaves(\mathcal{C})_{/Y} } \cr
}$$
whose value at $[n]$ is the $(n+1)$th iterate of $q^*q_!$ applied to $g$. It follows that
$$\textstyle{ R_g(-) \simeq \Map(f^*(-),g) \simeq \varprojlim_{\Delta} \Map(f^*(-),C_g([-])) \simeq \varprojlim_{\Delta} R_{C_g([-])}(-), }$$
and since accessible presheaves are closed under small limits, it suffices to show that $r_g$ is accessible in the case, where $g \simeq q^*(W)$ for some $W \in \presheaves(\mathcal{C})$. But
$$r_{q^*(W)}(-) \simeq \Map(f^*(h(-)),q^*(W)) \simeq \Map(q_!f^*(h(-)),W) \simeq F_{f,W}(-),$$
which is accessible by \cref{lemma:absolute_right_adjoint_to_pullback_exists}. This completes the proof.
\end{proof}

\begin{question}
\label{question:properties_of_macro_categories}
Let us say that an $\infty$-category $\mathcal{D}$ is:\footnote{\,The notions of class-accessible and class-presentable $\infty$-categories from~\cite{chornyrosicky} are more general and include $\infty$-categories of small presheaves on not necessarily coaccessible $\infty$-categories.}
\begin{enumerate}[label=\alph*)]
\item macroaccessible, if it can be written as $\mathcal{D} \simeq \Ind_{\kappa}(\mathcal{C})$ for some regular cardinal $\kappa$ and some coaccessible $\infty$-category $\mathcal{C}$.
\item macropresentable, if it is macroaccessible and cocomplete.
\item a macrotopos, if it is macropresentable and satisfies Giraud's axioms~(ii)--(iv).
\end{enumerate}
We expect that the following questions have affirmative answers:
\begin{enumerate}
\item If $\mathcal{D}$ is macroaccessible, then for every regular cardinal $\kappa$, is the full subcategory $\mathcal{D}^{\kappa} \subset \mathcal{D}$ spanned by the $\kappa$-compact objects coaccessible?
\item If $\mathcal{D}$ is macropresentable, then is it automatically complete?
\item Is it possible to formula an adjoint functor theorem that characterizes left adjoint functors between macropresentable $\infty$-categories?
\item Does there exist a good theory of localizations of macroaccessible categories, and if so, then can the properties of being macropresentable and of being a macrotopos be characterized in terms of being a suitable localization of an $\infty$-category of accessible presheaves?
\item Are macropresentable $\infty$-categories and macrotopoi closed under the formation of slice categories, limits, etc.?
\end{enumerate}
\end{question}

\section{Effective epimorphisms of sheaves}
\label{appendix:effective_epimorphisms_of_sheaves}

In this appendix, we spell out some folklore results on effective epimorphisms of sheaves for which we have not been able to find a suitable reference.

Let $\mathcal{C}$ be a small $\infty$-category, and let $h \colon \mathcal{C} \to \presheaves(\mathcal{C})$ be the Yoneda embedding into the $\infty$-category of presheaves of anima on $\mathcal{C}$. A sieve on $S \in \mathcal{C}$ is an equivalence class of $(-1)$-truncated maps $j \colon U \to h(S)$ in $\presheaves(\mathcal{C})$. It determines and is determined by the full subcategory $\mathcal{C}_{/S}(U) \subset \mathcal{C}_{/S}$ spanned by the maps $g \colon S' \to S$ with the property that $h(g)$ factors through $j$.

\begin{example}If $f \colon T \to S$ is a map in $\mathcal{C}$, then~\cite[Proposition~6.2.3.4]{luriehtt} shows that the colimit of its \v{C}ech nerve provides the essentially unique factorization
$$\xymatrix{
{ h(T) } \ar[r]^-{i_0} &
{ U \simeq \varinjlim_{\Delta^{\op}} h(T)^{\times_{h(S)}[-]} } \ar[r]^-{j} &
{ h(S) } \cr
}$$
of $h(f)$ as the composition of a $(-1)$-connected map and a $(-1)$-truncated map. Therefore, a map $g \colon S' \to S$ in $\mathcal{C}$ factors through $f$ if and only if the map $h(g)$ in $\presheaves(\mathcal{C})$ factors through $j$. In particular, the full subcategory of $\mathcal{C}_{/S}(U) \subset \mathcal{C}_{/S}$ spanned by the maps with this property is closed under coproducts.
\end{example}

\begin{example}If $(f_i \colon T_i \to S)_{i \in I}$ is a family of maps in $\mathcal{C}$, then there is a smallest sieve $j \colon U \to h(S)$ with the property that $h(f_i)$ factors through $j$ for all $i \in I$. We say that $j \colon U \to h(S)$ is the sieve generated by $(f_i \colon T_i \to S)_{i \in I}$.
\end{example}

\begin{definition}
\label{definition:grothendieck_topology}
Let $\mathcal{C}$ be a small $\infty$-category. 
A Grothendieck topology on $\mathcal{C}$ is an assignment to every $S$ in $\mathcal{C}$ of a set $J(S)$ of sieves on $S$, called the covering sieves on $S$, such that the following hold:
\begin{enumerate}
\item[(1)]The class of $\id_{h(S)} \colon h(S) \to h(S)$ is a covering sieve for every $S$ in $\mathcal{C}$.
\item[(2)]If $j \colon U \to h(S)$ is a covering sieve on $S$, and if $f \colon T \to S$ is any map in $\mathcal{C}$, then the base-change $j' \colon U' \to h(T)$ of $j$ along $h(f)$ is a covering sieve on $T$.
\item[(3)]Let $j \colon U \to h(S)$ and $k \colon V \to h(S)$ be sieves on $S$. If $k$ is a covering sieve, and if for every $f \colon T \to S$ with the property that $h(f)$ factors through $k$, the base-change $j' \colon U' \to h(T)$ of $j$ along $h(f)$ is a covering sieve on $T$, then $j$ is a covering sieve on $S$.
\end{enumerate}
A family $(f_i \colon T_i \to S)_{i \in I}$ of maps in $\mathcal{C}$ is a covering family, if the sieve $j \colon U \to h(S)$ that it generates is a covering sieve. A Grothendieck topology is finitary, if every covering sieve contains a finite covering family.
\end{definition}

\begin{remark}
\label{remark:covering_sieve}
Let $\mathcal{C}$ be a small $\infty$-category equipped with a Grothendieck topology. We claim that if $j \colon U \to h(S)$ is a sieve on $S$, which contains a covering sieve in the sense that there exists a covering sieve $k \colon V \to h(S)$ such that $k$ factors through $j$, then $j$ is a covering sieve. Indeed, if $f \colon T \to S$ has the property that $h(f)$ factors through $k$, then $h(f)$ also factors through $j$, and hence, the base-change $j' \colon U' \to h(T)$ of $j$ along $h(f)$ is a $(-1)$-truncated map, which admits a section. But then $j'$ is equivalent to the identity map, and hence, is a covering sieve by~(1), so we conclude from~(3) that $j$ is a covering sieve, as claimed.
\end{remark}

\begin{definition}
\label{definition:local_epimorphism}
Let $\mathcal{C}$ is a small $\infty$-category equipped with a Grothendieck topology, let $f \colon Y \to X$ be map of presheaves of anima on $\mathcal{C}$, and let $j \colon U \to X$ be the canonical map from the colimit of its \v{C}ech nerve. The map $f \colon Y \to X$ is a local epimorphism if the base-change $j_S \colon Y_S \to h(S)$ of $j$ along any map $\eta \colon h(S) \to X$ with $S$ in $\mathcal{C}$ is a covering sieve.
\end{definition}

\begin{proposition}
\label{proposition:characterization_of_effective_epimorphisms_of_sheaves}
Let $\mathcal{C}$ be a small $\infty$-category equipped with a Grothendieck topology, and let $f \colon Y \to X$ be a map of presheaves of anima on $\mathcal{C}$. The following are equivalent:
\begin{enumerate}
\item[{\rm (1)}]The map $L(f) \colon L(Y) \to L(X)$ in $\sheaves(\mathcal{C})$ is an effective epimorphism.
\item[{\rm (2)}]The map $f \colon Y \to X$ is a local epimorphism.
\item[{\rm (3)}]For every $S$ in $\mathcal{C}$ and every $S$-valued point $\eta \in X(S)$, there exists a diagram
$$\xymatrix@C=10mm{
{ \coprod_{i}h(T_i) } \ar[r]^-{\tilde{\eta}} \ar[d]^-{\sum_{i}h(g_i)} &
{ Y } \ar[d]^-{f} \cr
{ h(S) } \ar[r]^-{\eta} &
{ X } \cr
}$$
with $(g_i \colon T_i \to S)_{i \in I}$ a covering family. If the Grothendieck topology is finitary, then said covering family can be taken to be finite.
\end{enumerate}
\end{proposition}

\begin{proof}We recall from~\cite[Proposition~6.2.3.4]{luriehtt} that
$$\xymatrix{
{ Y } \ar[r]^-{i_0} &
{ U \simeq \varinjlim_{\Delta^{\op}} Y^{\times_X[-]} } \ar[r]^-{j} &
{ X } \cr
}$$
is the essentially unique factorization of $f$ as the composition of a $(-1)$-connected map and a $(-1)$-truncated map. Now, to prove that~(1) and~(2) are equivalent, let us say that $\eta \colon X' \to X$ is \emph{good}, if the base-change $j' \colon U' \to X'$ of $j$ along $\eta$ is a covering sieve, or equivalently, if $L(j') \colon L(U') \to L(X')$ is an equivalence. So~(1) is the statement that the identity map $\id_X \colon X \to X$ is good, whereas~(2) is the statement that every map $\eta \colon S \to X$ from a representable presheaf is good. Now, if $\eta \colon X' \to X$ is good, and if $g \colon X'' \to X'$ is any map, then $\eta g \colon X'' \to X$ is good, so~(1) implies~(2). And the full subcategory of $\presheaves(\mathcal{C})_{/X}$ spanned by the good maps $\eta \colon X' \to X$ is closed under colimits, so~(2) implies~(1).

Finally, to prove that~(2) and~(3) are equivalent, we note that~(2) is the statement that for every $S$ in $\mathcal{C}$ and every $\eta \in X(S)$, the base-change $j_S$ of $j$ along $\eta$ is a covering sieve. But the sieve $j_S$ is a covering sieve if and only if there exists a covering family $(g_i \colon T_i \to S)_{i \in I}$ and a diagram
$$\xymatrix@C=10mm{
{ \coprod_{i \in I} h(T_i) } \ar[r]^-{\tilde{\eta}} \ar[d]^-{\sum_i h(g_i)} &
{ U } \ar[d]^-{j} \cr
{ h(S) } \ar[r]^-{\eta} &
{ X. } \cr
}$$
Moreover, this diagram exists if and only if the diagram in~(3) exists, so we have proved that~(2) and~(3) are indeed equivalent.
\end{proof}

If a covering family $(g_i \colon T_i \to S)_{i \in I}$ consists of a single map $g \colon T \to S$, then we will also say that the map $g \colon T \to S$ is a covering.

\begin{corollary}
\label{corollary:characterization_of_effective_epimorphisms_of_sheaves}
Let $\mathcal{C}$ be a small $\infty$-category equipped with a finitary Grothendieck topology, and suppose that $\mathcal{C}$ admits finite coproducts. A map $f \colon Y \to X$ of sheaves of anima on $\mathcal{C}$ is an effective epimorphism if and only if for every $S$ in $\mathcal{C}$ and every map $\eta \colon h(S) \to X$, there exists a diagram
$$\xymatrix@C=10mm{
{ h(T) } \ar[r]^-{\tilde{\eta}} \ar[d]^-{h(g)} &
{ Y } \ar[d]^-{f} \cr
{ h(S) } \ar[r]^-{\eta} &
{ X } \cr
}$$
with $g \colon T \to S$ a covering.
\end{corollary}

\begin{proof} \cref{proposition:characterization_of_effective_epimorphisms_of_sheaves} show that $f \colon Y \to X$ is an effective epimorphism of sheaves on $\mathcal{C}$ if and only if for every $\eta \colon h(S) \to X$ with $S$ in $\mathcal{C}$, there exists a diagram
$$\xymatrix@C=11mm{
{ \coprod_{i \in I}h(T_i) } \ar[r]^-{\tilde{\eta}} \ar[d]^-{\sum_ih(g_i)} &
{ Y } \ar[d]^-{f} \cr
{ h(S) } \ar[r]^-{\eta} &
{ X } \cr
}$$
with $(g_i \colon T_i \to S)_{i \in I}$ a finite covering family. But $Y$ preserves finite products, so the map $\tilde{\eta}$ factors essentially uniquely through the canonical map
$$\xymatrix{
{ \coprod_{i \in I} h(T_i) } \ar[r] &
{ h(\coprod_{i \in I} T_i), } \cr
}$$
and since $(g_i \colon T_i \to S$) is a covering family, we conclude from \cref{remark:covering_sieve} that the induced map $g \colon T \simeq \coprod_{i \in I} T_i \to S$ is a covering.
\end{proof}

\providecommand{\bysame}{\leavevmode\hbox to3em{\hrulefill}\thinspace}
\providecommand{\MR}{\relax\ifhmode\unskip\space\fi MR }
\providecommand{\MRhref}[2]{%
  \href{http://www.ams.org/mathscinet-getitem?mr=#1}{#2}
}
\providecommand{\href}[2]{#2}

\end{document}